\documentclass{article}
\usepackage{amsmath,amssymb,amsfonts,amsthm, bm}
\usepackage{float, caption, subcaption, wrapfig, color, verbatim, multirow, algpseudocode, soul, cancel, graphicx}
\usepackage[top=2cm, bottom=2cm, left=2cm, right=2cm]{geometry} 
\newtheorem{remark1}{Remark}[section]
\newtheorem{thm}{Theorem}
\newtheorem{statement1}[thm]{Statement}

\providecommand{\keywords}[1]
{
  \small	
  \textbf{\textit{Keywords---}} #1
}

\title{ Computer-assisted global analysis for vibro-impact dynamics:\\ a reduced smooth maps approach
\thanks{ This work was supported by the  National Science Foundation (USA) under collaborative grants No. CMMI-2009329 and CMMI 2009270, and Engineering and Physical Sciences Research Council (UK) under grant
EPSRC EP/V034391/1.}}

\author{Lanjing Bao\thanks{Department of Mathematics and Statistics, Georgia State University, P.O. Box 4110, Atlanta, Georgia, 30302-410, USA} ({\tt lbao1@student.gsu.edu})
	\and Rachel Kuske\thanks{School of Mathematics, Georgia Institute of Technology,
			686 Cherry Street Atlanta, GA 30332-0160 USA} ({\tt rachel@math.gatech.edu})
	\and Daniil Yurchenko\thanks{ Institute for Sound and Vibration Research, University of Southampton, Southampton SO17 1BJ, UK
			} ({\tt d.yurchenko@soton.ac.uk})
	\and 
		Igor Belykh\thanks{Department of Mathematics and Statistics and Neuroscience Institute, Georgia State University, P.O. Box 4110, Atlanta, Georgia, 30302-410, USA} ({\tt ibelykh@gsu.edu})
		}

\begin{document}

\maketitle
\begin{abstract}
 We present a novel approach for studying the global dynamics of a vibro-impact pair, that is, a ball moving in a harmonically forced capsule.  Motivated by a specific context of vibro-impact energy harvesting, we develop the method with broader non-smooth systems in mind. The traditional maps between impacts of the ball with the capsule are implicit and transcendental and, therefore, not amenable to global analysis.  Nevertheless, we exploit the impacts as useful non-smooth features to select appropriate return maps that provide a path for studying global behavior.  This choice yields a computationally efficient framework for constructing return maps on short-time realizations from the state space of possible initial conditions rather than via long-time simulations often used to generate more traditional maps.  The different dynamics in sub-regions in the state space yield a small collection of reduced polynomial approximations. Combined into a piecewise composite map, these capture transient and attracting behaviors and reproduce bifurcation sequences of the full system. Further ``separable'' reductions of the composite map provide insight into both transient and global dynamics. This composite map is valuable for cobweb analysis, which opens the door to computer-assisted global analysis and is realized via conservative auxiliary maps based on the extreme bounds of the maps in each subregion. We study the global dynamics of energetically favorable states and illustrate the potential of this approach in broader classes of dynamics.

\end{abstract}



\section{Introduction} \label{S1intro}
The prevalence of non-smooth dynamics, characterized by switches, impacts, sliding, and other abrupt alterations in behavior, permeates various fields, including physics, biology, and engineering \cite{Andronov59,Fil88,BernChamp2007}.  Non-smooth dynamical models are essential for understanding phenomena such as body component interactions with non-smooth contacts, impacts, friction, and switching in mechanical systems \cite{FeigBern99,SHAW201983,leine2013dynamics,belykh2021emergence}, and relay systems, switched power converters, and packet-switched networks in electrical and control engineering \cite{FeigBern99,di1998switchings,benadero2003two,hasler2013dynamics}.
In the life sciences, non-smooth dynamics are evident in diverse systems such as gene regulatory networks \cite{polynikis2009comparing,acary2014numerical} and pulse-coupled neurons \cite{ermentrout2019recent}.
While piecewise smooth, non-smooth, and vibro-impact dynamical systems represent vast research fields in nonlinear science, they have historically received far less attention than their smooth counterparts. In recent decades, increased efforts have pursued a comprehensive understanding of non-smooth bifurcations and related nonlinearities (see extensive reviews \cite{BernChamp2007,Je18b,Je20,beyond2023belykh} and references therein). Non-smooth systems and the vibro-impact systems we study here fall into the larger class of hybrid systems, whose breadth is reflected in combinations of discrete and continuous components with complementary features \cite{Clerk2020_hybrid},  dynamics obtained from combined models and measured or experimental data \cite{Sieber_2008}, embedded control systems  \cite{branicky2005introduction},  and perhaps even by the Wikipedia description of systems that  ``can both flow and jump'' \cite{hybrid-wiki}.

Vibro-impact (VI) systems constitute a distinct class of dynamical systems where impacts substantially influence the nonlinear behavior. Typical classes of VI systems include a forced mass and one or more stationary rigid barriers or, alternatively, a pair of moving impacting masses, each of which may be subject to external forcing. Classic examples include balls bouncing on moving surfaces \cite{luo1996dynamics,leine2013dynamics,leine2012global}, pendulums impacting barriers \cite{shaw1989transition,di2008bifurcations,tang2019periodic}, and VI pairs composed of two oscillating masses that impact each other \cite {Luo2013}. Generally, both masses in the VI pair may undergo forcing, complemented by elastic or inelastic impacts. A canonical VI pair, considered in this paper, consists of a forced capsule, with an inner mass moving freely within a cavity of a given length and impacting the ends of the capsule. This concept has been explored as an effective vibration mitigation alternative to linear tuned mass dampers or continuous nonlinear dampers \cite{Yu2014, Thota2006, ZhangFu2019, Mason2011, LiVINES2022, Li2017, DANKOWICZ2005, LUO2008}.  Recently, a VI pair was proposed as an energy harvesting mechanism, where the impacts between the inner mass and the capsule deform flexible dielectric polymer membranes on the capsule ends \cite{Yurchenko2017}. These membranes serve as capacitors, as the impacts deform them and change their capacitance, thus enabling energy harvesting \cite{lai2018energy}.  Previously, VI pairs have been studied by approximate methods, including averaging, multiple scales, and complexification averaging \cite{Dimentberg2004, Ibrahim2009, Li2017, Vakakis_rev22}, but with limited applicability to non-smooth systems with impacts. 

Recently, VI pair systems have been studied precisely using maps,  combining the system's motion between the impacts and the impact conditions \cite{Guo2012,Guo2012n,Fu2015,ZhangFu2019,luo2010switching}. The semi-analytical solution of these exact equations can provide exhaustive information regarding the bifurcation structure and local stability of different types of motion. In the case when the smaller mass is negligible relative to the larger one, between impacts this two-degree-of-freedom system can be reduced to a single differential equation for the relative displacement of the two masses \cite{Serdukova2019,  Luo2013}, used to explore, e.g., the interplay between classical and grazing bifurcations \cite{Serdukova2023} and comparisons of the influence of instantaneous and compliant impact conditions \cite{Costa2022}.
In settings where the smaller mass is non-negligible,  such as in targeted energy transfer, exact maps for the full system allow bifurcation analyses over a large range of parameters for modes with efficient energy transfer and their loss of stability to inefficient alternating chatter behaviors \cite{Kumar2024}.

These previous map-based results are primarily based on linear stability analyses, leaving a critical gap in analyzing the global, possibly chaotic dynamics of VI systems due to severe limitations of the existing global stability methods in handling impacts. 
One limiting factor in pursuing existing approaches of global analyses for the forced VI pair is that it is non-autonomous. As a result, the maps traditionally used in its bifurcation analyses yield analytically intractable coupled transcendental maps for the system response and impact time,  preventing explicit expressions for the system's state that one would normally use to study global dynamics. This fact motivates the new approach we propose here.

In a broader context, global stability approaches for non-autonomous, non-smooth systems are few and far between, which appears to be true also for those hybrid systems described with both continuous and discrete dynamics \cite{kowalczyk2012modelling,belykh2021emergence,CHEN2021_hybrid,Glendinning_2011,Chawla_2022}.
 One notable example is an extension of the Lyapunov function method to prove the global stability of the equilibrium state of a non-autonomous bouncing ball \cite{leine2012global}. In this setting, the Lyapunov-type method involves non-autonomous measure differential inclusions and constructs a decreasing step function above an oscillating Lyapunov function. However, its application to non-trivial dynamics of VI pairs with two-sided impacts seems elusive. Another notable sample is an averaging Lyapunov function approach developed to prove global convergence to absorbing domains of non-trivial attractors in non-smooth dynamical systems with a non-autonomous stochastic switching parameter rule \cite{hasler2013dynamics}. However, this approach is irrelevant for non-autonomous VI systems as it is based on knowledge of the averaged autonomous system's attractor. 
 Recently, a computer-assisted proof of chaos in piecewise linear maps was obtained by explicitly constructing trapping regions and invariant cones based on word sets representing the dynamics symbolically \cite{Simpson2023}.
An area-preserving map-based analysis for the global behavior of the VI pair's rare, restricted behavior was proposed in \cite{Cao2022}.
Yet, to date, there appear to be no global analyses relevant to applications such as energy harvesting, for which the VI pair dynamics of interest include sustained sequences of regular impacts on both barriers at the capsule ends, observed over a large range of parameters. Then, we are faced with the challenge of global analyses of behavior with at least two (alternating) impacts per forcing cycle.  This feature is in contrast with other studies of impacting systems that may consider the transition between no impacts and a single impact  \cite{Nordmark1991}, repeated impacts on a single barrier \cite{Simpson2018}, or the global attraction of a solution without impacts \cite{leine2012global}.

In this paper, we present a novel computer-assisted approach for studying the global dynamics of the VI pair, that is, a ball moving in a harmonically forced capsule. Motivated to develop an analytical global analysis for this system, we prioritize approaches that include explicit expressions wherever possible. As mentioned above, the repeated impacts of the ball with the capsule yield transcendental maps that are analytically intractable within existing global analyses. Yet, we exploit them as useful non-smooth features in constructing novel two-dimensional (2D) return maps that separate families of important sequences in the VI-pair dynamics. These families are used to characterize global dynamics and can be related to bifurcations of the VI pair.  Computationally efficient short-time realizations of these return maps divide the state space according to the different dynamics of these families. Our definition of return maps does not fall into standard choices for maps, such as Poincar\'{e}, stroboscopic, all impacts, or all returns to a particular state \cite{Luo2013, Nordmark1991, Pavlov2004, Simp20}. Instead, it divides the return maps based on the sequence of impacts that do or do not occur before the system returns to a particular impacting state. 
This innovative perspective is valuable for efficiently partitioning the state space into a few regions corresponding to distinct surfaces formed by maps from different families of the sequences of impacts. Identifying the regions with potentially attracting and transient behavior is straightforward by inspecting the surfaces' geometry and gradients relative to the diagonals in the phase planes for impact velocity and impact phase.  
As a result, we can focus a more detailed analysis on smaller regions with potentially attractive behavior. 
These computationally-realized return maps could be directly used for purely numerical yet efficient cobweb analysis of the system's global behavior. However, toward our goal of performing computer-assisted analysis and explicitly characterizing the system's global dynamical properties, we go one step further and define reduced polynomial approximations for the maps in each region.

Combining these polynomials into a piecewise smooth composite map, we demonstrate that it captures transient behaviors throughout the state space while reproducing the attracting behaviors. Furthermore, it reproduces an important sequence of period-doubling bifurcations and (apparently) chaotic behavior compared with the bifurcation sequences of the exact systems. 
In constructing the composite map, we find that in some regions with strongly transient dynamics, we can reduce the 2D return maps to a pair of 1D return maps without sacrificing the integrity of the attracting dynamics. While not a necessary step, these types of ``separable'' components of the composite map provide transparency for the overall dynamics.
Furthermore, this composite map derived from the non-smooth VI dynamics is remarkably valuable for cobweb analysis, as it is based on simple return maps corresponding to impacts on one end of the capsule rather than on compositions of map sequences.
Specifically, the separable representations of the 2D map are convenient for visualizations within this cobweb phase analysis that captures the different attracting behaviors for different parameter regimes.

Notably, this cobweb analysis motivates a valuable definition of auxiliary maps on the regions identified within the construction of the composite map, once the transient and attracting characteristics have been identified.  For regions with attracting dynamics, the auxiliary map is conservatively based on the extreme bounds on the map for each region and thus can be used to bound the attracting domain. A key feature of the auxiliary maps is that they simplify the 2D return maps into a set of 1D equations using the bounds for each region. Then, a cobweb phase space analysis is used to explore the system's long-term dynamics, and yields a limiting period-two cycle that bounds the attracting domain that contains all the system's non-trivial attractors.  With the auxiliary maps based on the polynomial approximations, we can obtain analytical expressions for the impact velocity map and, thus, for the attracting domain. Repeated application of the auxiliary maps, each with updated bounds obtained from the previous application, yields tighter bounds for the attracting domain.
  
 We outline the process of generating the approximate composite map in terms of a general algorithm adaptable for other non-smooth dynamical systems. A key step in the algorithm includes identifying families of short sequences of impacts that give the building blocks for the return maps. The resulting division of the state space is relatively simple and computationally efficient compared to, e.g., the identification of basins of attraction, which require long time computations to find complex regions for dynamics sensitive to initial conditions. Likewise, flow-defined Poincar\'{e} maps for the global dynamics of periodic and chaotic systems, derived from long-time simulations over the entire state space,  are often piecewise smooth even though they originate from a smooth dynamical system. Geometrical piecewise smooth Lorenz maps \cite{afraimovic1977origin,robinson1989homoclinic,guckenheimer1983nonlinear} representing the smooth chaotic dynamics of the Lorenz system are notable examples. Our approximate composite map constructed for only short-time realizations of the VI pair is conceptually different from classical piecewise smooth maps with regular and chaotic dynamics appearing in various biological, social science, and engineering applications \cite{Nordmark1997,avrutin2019continuous,zhusubaliyev2003bifurcations,belykh2011belykh,coombes1999mode,Gl14b,DeGl21}. However, it can still be interpreted as a geometrical model of the VI pair as it depicts the dynamics and bifurcations remarkably well and derives from a polynomial approximation of the state space partitions. The combination of the geometric interpretation and the polynomial approximation facilitates our goal of obtaining analytical results for the global dynamics directly related to the physical model. These results are in contrast to local analyses and computational studies of higher dimensional maps \cite{Pavlov2004,Sahari2021}. 

In this first development of the approach, we focus on parameter regimes for behaviors that drive favorable energy output in a VI pair-based energy harvesting device, behaviors with alternating impacts on either end of the capsule. The impact velocity and phase may repeat periodically with period $n{\cal T}$, where ${\cal T}$ is the period of the forcing, or the states may have apparently chaotic behavior within the alternating behavior. 
Besides its physical relevance, this choice of parameters facilitates a relatively straightforward presentation of the approach while exploring several types of non-trivial dynamics. Nevertheless, as discussed further in the conclusions, we expect that foundational concepts in this analysis are adaptable to other (more complex) sequences of impacts.
 
The remainder of the paper is organized as follows. Section 2 gives details of the VI pair model, including the transcendental form of the maps \cite{Serdukova2021,Serdukova2023} that motivates the computer-assisted analysis of global dynamics. Section 3 provides the return maps that form the building blocks of the computer-assisted approach, illustrating their key properties.  Section 4 provides the general algorithm for constructing a composite map realized for the VI pair by approximating the return maps with explicit piecewise polynomial maps over relevant regions that comprise the state space.  Section 5  compares the trajectories generated using the exact and composite maps in the state space and the phase plane. Section 6 develops an auxiliary map based on the composite map to identify the globally attracting dynamics and the corresponding domain for three qualitatively different types of the VI pair system behavior.
Section 7 contains conclusions and a brief illustration of the relevance of the approach for a VI pair-based energy harvesting device with stochastic forcing. Finally, Appendix A provides additional details on the construction of the return map. The supplementary material contains the exact map derivation and demonstrates its analytical intractability. It also contains the coefficients of the polynomials used in the composite map. Supplementary videos provide additional visualizations for constructing and iterating the composite map.

\section{The Model}\label{S2model}
The model takes the form of the canonical impact pair, comprised of an externally forced capsule with a freely moving ball inside. The friction between the ball and the capsule is neglected, so the ball's motion is driven purely by gravity and impacts one of the membranes on the capsule's ends.

One application based on the impact pair is a nonlinear vibro-impact energy harvesting device.  Each end of the capsule is closed by a membrane of dielectric (DE) polymer material with compliant electrodes \cite{Yurchenko2017}. 
The deformation of such a DE membrane is the vibro-impact energy harvesting device's primary means of energy generation. When the ball collides with the membrane, this action changes the ball's trajectory and deforms the membrane. The DE membrane's physical property, being a variable capacitance capacitor, allows the change of its capacitance when it is deformed; meanwhile, a bias voltage is applied when the deformation reaches its maximum state. After the collision, an extra voltage charge is harvested, and the membrane returns to its undeformed state.

\begin{figure}[htbp]
    \begin{subfigure}[T]{0.55\textwidth}
        \centering
        \caption{}
        \includegraphics[width=\linewidth]{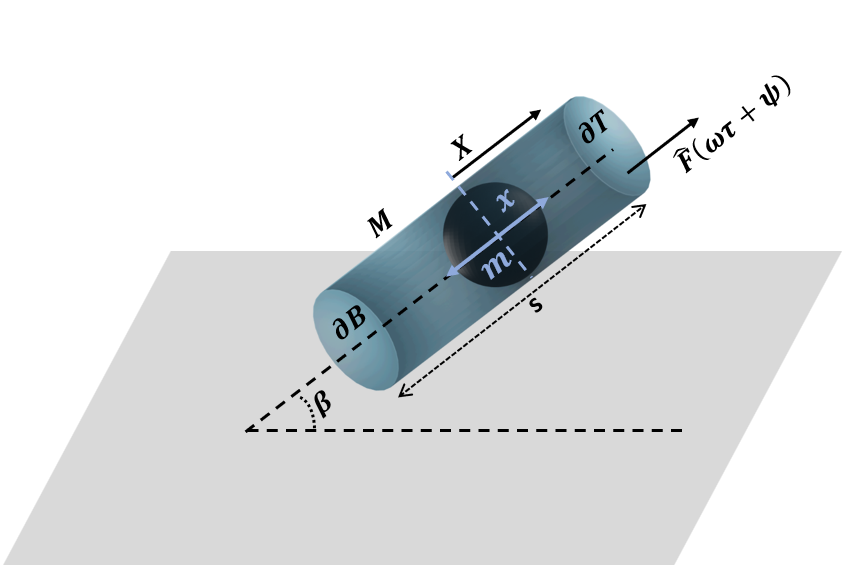}
    \end{subfigure}
    \begin{subfigure}[T]{0.35\textwidth}
        \centering
        \begin{subfigure}{\textwidth}
        \caption{}
        \includegraphics[width=\linewidth]{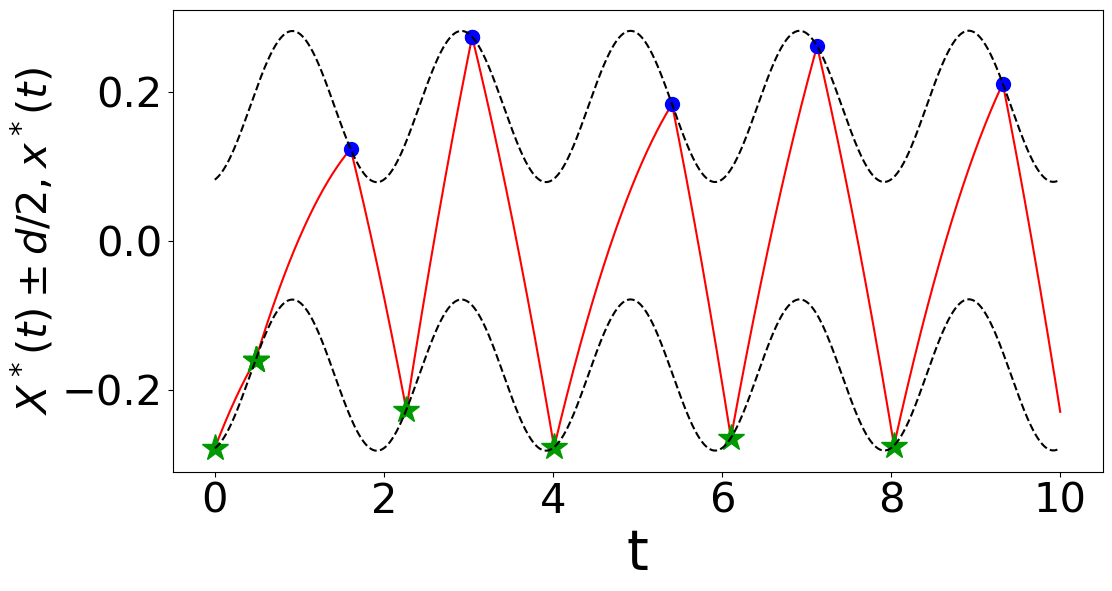}
        \end{subfigure}
        \begin{subfigure}{\textwidth}
        \caption{}
     \includegraphics[width=\linewidth]{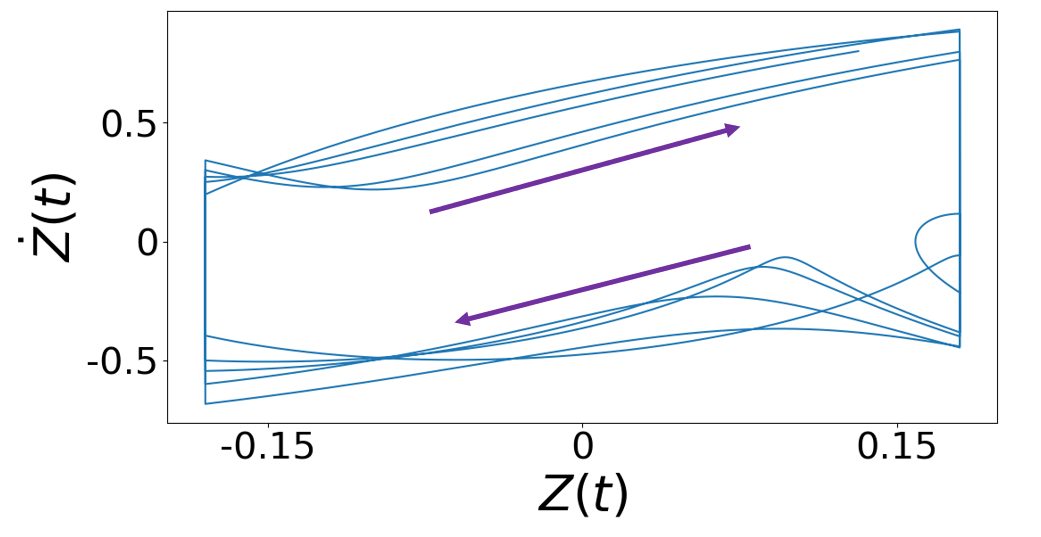}
        \end{subfigure}
    \end{subfigure}
    \caption{(a): Illustration of the VI pair:  A ball moves freely within a harmonically forced capsule enclosed by deformable membranes on both ends. The capsule is positioned with an angle $\beta$ relative to the horizontal plane and is excited by an external harmonic excitation $\hat{F}(\omega \tau + \psi)$. The mass, length of the capsule, and mass of the ball are $M, s$, and $m$, respectively. (b):  The two dashed black lines represent the displacement of the top and bottom membranes, $X(t)^*\pm d/2$. The green stars and blue dots indicate the impacts at $\partial B$ and $\partial T$, respectively. The red solid lines connect each impact at $\partial T$ and $\partial B$, representing the estimated ball movement between each impact. (c): Phase plane in terms of relative variables. The relative displacement $Z(t)$ oscillates between $- d/2$ and $ d/2$, and the relative velocity $\dot{Z}(t)$ has a sign change at each impact. The arrows indicate the direction of time. Parameters: $d=0.35$, $\dot{Z}_0 = 0.43$ and $\psi_0=0.26$.}
    \label{model:abc}
\end{figure}

The schematic for the VI pair is given in Fig. \ref{model:abc}(a). Neglecting the friction, the system is driven by forces generated at impact, gravity, and external harmonic excitation $\hat{F}(\omega \tau + \psi)$ with period $2 \pi/\omega$. Using Newton's Second Law of Motion, the model is described by the following differential equations:
\begin{align}
    \frac{d^2 X}{d\tau^2}&=\frac{\hat{F}(\omega \tau + \psi)}{M}, \label{cyl_dim}\\
    \frac{d^2x}{d \tau^2}& = -g\sin\beta,   \label{ball_dim}
\end{align}
where $X(\tau)$ and $x(\tau)$ are the dependent variables for the absolute displacement for the capsule and the ball, respectively. In addition, $M$ and $m$ are the mass of the capsule and the ball, respectively. Also, $\psi$ is the general phase of the forcing.

Treating the impact time as negligible compared to other time scales in the model, we use an instantaneous impact model given by
\begin{eqnarray}\label{simpact_dim}
\left(\frac{d  x}{d \tau}\right)^+=- r \left(\frac{d  x}{d \tau}\right)^- + (1+r) \left(\frac{d  X}{d \tau}\right).
\end{eqnarray}
Note that this is a reduced model based on the condition $M\gg m$,  as discussed in detail in \cite{Serdukova2021}.
The superscripts $+$ and $-$ signify the state of the ball after and before the impact, respectively.
The parameter $r$ is the restitution coefficient, which is a quantitative measure of the membrane's elasticity. The range of $r$ is $[0,1]$ with $r=1$ being perfectly elastic and $r=0$ being inelastic. In this paper, we consider moderate elasticity $r=0.5$.
Additionally, in \eqref{simpact_dim}, we do not distinguish the states before and after the impact for the capsule $dX/d\tau$ because the mass of the ball ($M\gg m$) is negligible and does not change the state of the capsule at impact.

To focus on the system's dependence on key parameters, we first non-dimensionalize the system. Following  \cite{Serdukova2021}, the dimensionless variables $X^*(t), \dot{X}^*(t), t$ are the following:
\begin{eqnarray}\label{nondXYt}
X(\tau)=\frac{\parallel \hat{F} \parallel \pi^2}{M \omega^2} \cdot X^*(t), \ \  \frac{d {X}}{d \tau}=\frac{\parallel \hat{F} \parallel\pi}{M \omega} \cdot \dot{X}^*(t), \ \  \tau=\frac{\pi}{\omega} \cdot t \, ,
\end{eqnarray}
where  $\parallel \hat{F} \parallel$ is an appropriately defined norm of the strength of the forcing $\hat{F}$. Here, we also use Newton's dot notation for differentiation when the derivative is calculated with respect to dimensionless time $t$. 

In addition to non-dimensionalization, relative variables are used to focus on the system dynamics as a whole rather than the separate motion of the ball and capsule.   Using the variables $X^*$, the relative displacement $Z(t)$ and relative velocity $\dot{Z}(t)$  are given in the dimensionless form: 
\begin{eqnarray}\label{dimlessrel}
& & {Z}={X}^*-{x}^*, \qquad \dot{Z}=\dot{X}^*-\dot{x}^*,\nonumber\\
& &\ddot{Z}=\ddot{X}^*-\ddot{x}^*= F(\pi t+ \psi)+ \dfrac{M g \sin \beta}{\parallel \hat{F} \parallel}, \label{eq_relative}
\end{eqnarray}
where the non-dimensional forcing $F(\pi t+\psi) = \frac{\hat{F}(\omega \tau+\psi)}{||\hat{F}||}$ has the unit norm, i.e. $\parallel F \parallel = 1$. For convenience, we define $\bar{g}=\dfrac{M g \sin \beta}{\parallel \hat{F} \parallel}$.

 Since we want to evaluate the system from one impact to the next, the system's state at each impact is particularly important. Combining conditions (\ref{nondXYt}), (\ref{dimlessrel}), the impact condition (\ref{simpact_dim}) can be rewritten using  $Z$ and $\dot{Z}$. For the $j^{\text{th}}$ impact occurring at time $t=t_j$, 
\begin{align}
&Z_j=X^*(t_j)-x^*(t_j)=\pm \frac{d}{2}, \ \
\mbox{ for } x \in \partial B \; (\partial T) \mbox{ the sign is } +(-), \nonumber \\
&\dot{Z}_j^+=-r \dot{Z}_j^-\, , \label{eq_impact} \qquad
d=\frac{s M \omega^2}{\parallel \hat{F} \parallel \pi^2}.
\end{align}
The notations $\partial B$ and $\partial T$ denote the bottom and top membranes, respectively. The parameter $d$ is the dimensionless length of the system, used throughout this paper as the bifurcation parameter. In contrast to the actual length of the capsule $s$, $d$ varies with multiple factors, including the device length ($s$), mass ($M$), angular velocity of the external force ($\omega$), and forcing strength ($\parallel \hat{F} \parallel$). As illustrated in Fig. \ref{model:abc}(b),(c), the relative position of the system is bounded, $Z(t)\in [-d/2, d/2]$.  At the impacts, which is when $Z_j = \pm d/2 $, the relative velocity $\dot{Z}_j$ changes sign: when the impact is on $\partial B$ ($Z_j = d/2$), $\dot{Z}$ changes from positive to negative; when the impact is on $\partial T$ ($Z_j = -d/2$), $\dot{Z}$ switches from negative to positive. Since the displacement at each impact is known, either $Z_j=d/2$ or $Z_j=-d/2$, the relative velocity and time $(\dot{Z}_j, t_j)$ describe the system state at the $j^{\rm th}$ impact.

We summarize results from \cite{Serdukova2021} for calculating the exact maps for $(\dot{Z}_j, t_j)$ between two consecutive impacts. Between the impact at $t_j$ and the next impact at $t_{j+1}$, the relative velocity and displacement can be derived by integrating \eqref{dimlessrel} for $t\in (t_j, t_{j+1})$ and applying \eqref{eq_impact}:
\begin{align}\label{relvar}
    \dot{Z}(t) &= -r \dot{Z}^-_j +\Bar{g}\cdot (t-t_j)+F_1(t)-F_1(t_j),\nonumber\\
    Z(t) &= Z^+_j 
    -r \dot{Z}^-_j \cdot (t-t_j)+\frac{\Bar{g}}{2}\cdot (t-t_j)^2 + F_2(t) - F_2(t_j) - F_1(t_j)\cdot(t-t_j),
\end{align}
 where $F_1(t) = \int F(\pi t+\psi)\; dt$ and $F_2(t) = \int F_1(t)\; dt$. At the $j^{\text{th}}$ impact, $Z_j^+=Z_j^-$. Therefore, the superscripts in ${\dot{Z} }^{\pm}$ are omitted, since \eqref{relvar} are in terms $Z^- \text{ and } \dot{Z}^-$ only. Using the equations (\ref{relvar}), there are four basic nonlinear maps $P_{BB}, P_{BT}, P_{TB}, P_{TT}$ corresponding to motion between consecutive impacts, in terms of the four combinations of impact locations: $\partial B \rightarrow \partial B, \; \;\partial B \rightarrow \partial T, \; \; \partial T \rightarrow \partial B, \; \; \partial T \rightarrow \partial T$.  All four maps take the form
\begin{align}\label{basicmapdZ}
\dot{Z}_{j+1} &= -r \dot{Z}_j +\Bar{g}\cdot (t_{j+1}-t_j)+ F_1(t_{j+1})-F_1(t_j),\nonumber\\
\pm\frac{d}{2} &=  \pm\frac{d}{2} -r \dot{Z}_j \cdot (t_{j+1}-t_j)+\frac{\Bar{g}}{2}\cdot (t_{j+1}-t_j)^2 +F_2(t_{j+1}) - F_2(t_j) - F_1(t_j)\cdot(t_{j+1}-t_j).
\end{align}
The maps derived above are based on the system \eqref{relvar}, which gives the exact map when evaluated at impact times $t=t_j$; specifically,  $P_{\ell}: (\dot{Z}_j,t_j) \to (\dot{Z}_{j+1},t_{j+1})$ for  $\dot{Z}_j = \dot{Z}(t_j)$.
 Notice, the sign for  $\pm d/2$ is chosen depending on the impact locations of $Z_{j}, Z_{j+1}$, $+$ ($-$) for $\partial B$ $(\partial T$).

Ideally, we would like to transform \eqref{basicmapdZ} into closed-form expressions for $(\dot{Z}_{j+1}, t_{j+1})$ in terms of $(\dot{Z}_{j}, t_{j})$, which can be used to analyze stability and other   (global) dynamic properties of these maps and their compositions. Furthermore, if we wish to determine the map for the first return to $\partial B$ for sequences as shown in Fig. \ref{model:abc}(b),(c), we would seek the exact map for the impact sequence   $\partial B \to \partial T \to \partial B$, or for two consecutive impacts on $\partial B$, which we refer to as  BTB or BB motion, respectively. Here, we use the simpler case of BB motion to demonstrate the difficulties in deriving closed-form expressions for such sequences. The map $P_{\rm BB}$ is given by \eqref{basicmapdZ}, using $Z_{j+1}=Z_j = d/2$, we have
\begin{align}\label{PBBt}
\dot{Z}_{j+1} &= -r \dot{Z}_j +\Bar{g}\cdot (t_{j+1}-t_j)+ F_1(t_{j+1})-F_1(t_j),\nonumber\\
\frac{d}{2} &=  \frac{d}{2} -r \dot{Z}_j \cdot (t_{j+1}-t_j)+\frac{\Bar{g}}{2}\cdot (t_{j+1}-t_j)^2 +F_2(t_{j+1}) - F_2(t_j) - F_1(t_j)\cdot(t_{j+1}-t_j).
\end{align}
For concreteness, we take $F(\pi t+\psi) = \cos(\pi t+\psi)$. Then $F_1(t) = \frac{1}{\pi} \sin(\pi t+\psi)$ and $F_2(t) = -\frac{1}{\pi^2}\cos(\pi t+\psi)$. Substituting these into \eqref{PBBt} and solving for  $(\dot{Z}_{j+1},t_{j+1}) $, we have
\begin{align}\label{PBBdZj+1}
\dot{Z}_{j+1} &= -r \dot{Z}_j +\Bar{g}t_{j+1} - \Bar{g}t_j+ \frac{1}{\pi}\sin(\pi t_{j+1}+\psi) -\frac{1}{\pi}\sin(\pi t_{j}+\psi),\\
    0 &= -r \dot{Z}_{j} t_{j+1} + r \dot{Z}_{j}t_{j}+\frac{\Bar{g}}{2} t_{j+1}^2 -\Bar{g} t_{j+1} t_{j} + \frac{\Bar{g}}{2} t_{j}^2 
    -\frac{1}{\pi^2}\cos(\pi t_{j+1}+\psi) + \frac{1}{\pi^2} \cos(\pi t_j + \psi) \label{PBBt1}\\
    &  
    - \frac{1}{\pi}\sin(\pi t_j+\psi) t_{j+1} + \frac{1}{\pi}\sin(\pi t_j+\psi)t_{j}. \nonumber
\end{align}
In \eqref{PBBdZj+1}, $\dot{Z}_{j+1}$ is a function of $\dot{Z}_j, t_j$, as well as $t_{j+1}$, determined from \eqref{PBBt1}.   Sorting terms containing  $t_{j+1}$ to simplify \eqref{PBBt1} yields 
\begin{align}\label{PBBtj1}
    \frac{\Bar{g}}{2} t_{j+1}^2 - \Big( r\dot{Z}_j + \Bar{g} t_j + \frac{1}{\pi}\sin(\pi t_j+\psi)\Big) t_{j+1} +\Big( r \dot{Z}_{j}t_{j} + \frac{\Bar{g}}{2} t_{j}^2 + \frac{1}{\pi^2} \cos(\pi t_j + \psi) + \frac{t_j}{\pi}\sin(\pi t_j+\psi) \Big)&\nonumber\\
       = \frac{1}{\pi^2} \cos(\pi t_{j+1}+\psi).\qquad \qquad \qquad \qquad \qquad \qquad&
\end{align}
From these expressions, we see that it is impossible to get a closed form expression for the state  ($\dot{Z}_{j+1}, t_{j+1}$) from the state at the previous impact, $(\dot{Z}_{j}, t_{j})$. This motivates the new approach that we discuss in detail in Section \ref{S3return-map}.
 The formulation in \eqref{basicmapdZ} is useful when determining conditions for periodic solutions with a fixed number of impacts, and their local stability. For example, as in \cite{Serdukova2021}, a composition of a fixed number of maps provides the basis for previous analyses of periodic solutions, and the corresponding linear stability analysis provides information about whether the periodic solutions are stable under small perturbations. In this previous work, different types of motion were generally categorized as n:m/$p{\cal T}$, where n and m are the numbers of impacts on $\partial B$ and $\partial T$, respectively, ${\cal T}$ is the excitation period, and $p$ is an integer. Furthermore, the impact pair has been demonstrated to yield n:m/$p{\cal T}$ and n:m/$C$ behaviors, with $C$ indicating complex, aperiodic, or chaotic behavior.

\begin{figure}[H]
\centering
\begin{subfigure}[b]{0.45\textwidth}
    \centering
    \caption{}
\includegraphics[width=\textwidth]{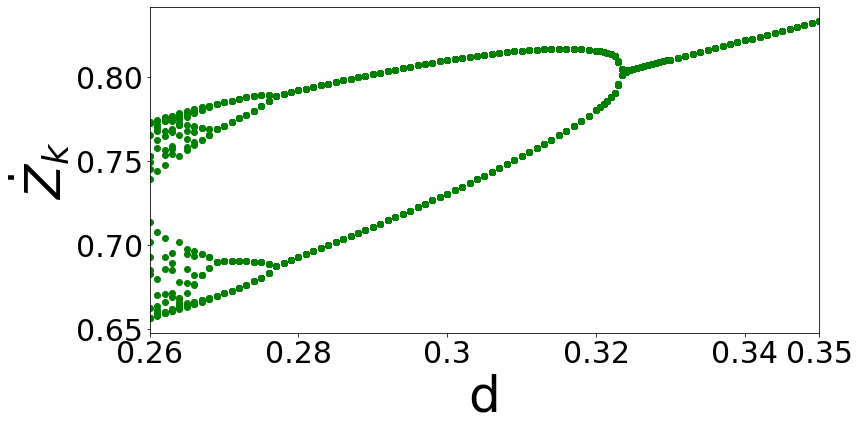}
\end{subfigure}
\hfill
\begin{subfigure}[b]{0.45\textwidth}
    \centering
    \caption{}
\includegraphics[width=\textwidth]{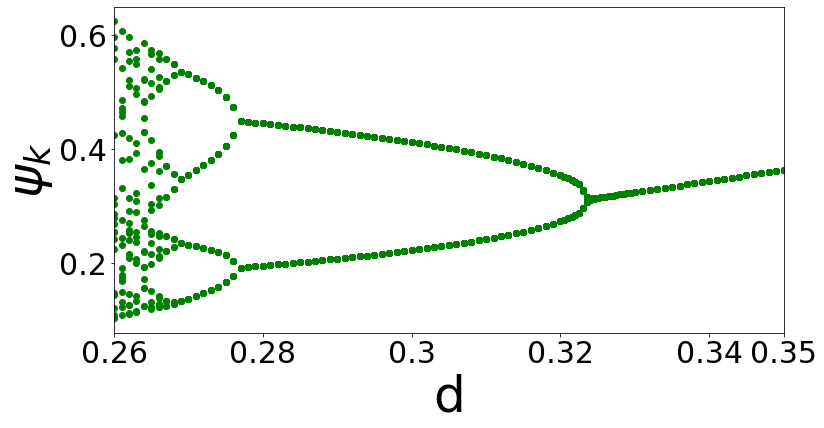}
\end{subfigure}
\caption{Bifurcation diagrams for $\dot{Z}_k$ and $\psi_k$ generated using the exact map from system \eqref{relvar}.
}\label{bifexact}
\end{figure}

In the remainder of this paper, we use $\psi_j =\mod(\pi t_j + \psi, 2\pi)$ rather than $t_j$ to quantify the impact timing within the forcing period of oscillation. Note that $\psi_j$ is  distinct from the general phase $\psi$  in the forcing term $F(\pi t+ \psi)$. This relative impact phase $\psi_j$ is more amenable than $t_j$ for considering transients and (quasi)-periodic behavior.
Figure~\ref{bifexact} shows the relative impact velocity $\dot{Z}_k$ and $\psi_k$ on $\partial B$, corresponding to a sequence of bifurcations with  1:1/${\cal T}$, 1:1/$p{\cal T}$ for $p$ an even integer, and 1:1/$C$ behavior over a range of the dimensionless length $d$. (Note: $\dot{Z}_k$ and $\psi_k$  on $\partial T$ not shown.)  We focus here on the parameters and the range of $d$ yielding 1:1-type behavior, with impacts alternating between $\partial B$ and $\partial T$ that is typically favorable for energy output and observed for the system \eqref{cyl_dim}-\eqref{simpact_dim}  over a large range of parameters  \cite{Serdukova2021, Serdukova2023}.

\begin{remark1}\label{num_bif}
 The numerical results in the bifurcation diagram (Fig. \ref{bifexact}) are generated by solving \eqref{cyl_dim}-\eqref{simpact_dim} over a long time, recording the limiting values for $\dot{Z}_k$ and $\psi_k$ on $\partial B$ for each value of $d$.  The attracting state then serves as the initial condition for the next value of $d$, using a continuation-type method with decreasing $d$.
Throughout this paper, the parameters used to generate the simulations are the following: $r=0.5$, $\lVert \hat{F} \rVert = 5$, $M =$ 124.5 g, $\omega=5\pi$, $\beta = \pi/3$, $g = 9.8$ m/${\rm s}^2$. Here, the non-dimensional parameter $d$ varies with the length of the capsule $s$, as given in  \eqref{eq_impact}. 
\end{remark1}

\section{Identification and visualization of the return maps} \label{S3return-map}

While the previous analyses capture the local stability of branches corresponding to periodic solutions, they do not provide information about the global attraction of this behavior or the potential for other attracting behavior.  In contrast, here, we seek to provide global stability results for the attraction of different types of solutions, including periodic,  nearly periodic, and chaotic behavior, as shown in Fig.~\ref{bifexact}. With that in mind, normally we would want to have the maps in an explicit form for the system state $(\dot{Z}_j,\psi_j)$. Equation \eqref{PBBtj1} has a solution if the quadratic function on the left-hand side (LHS) and the cosine function on the right-hand side (RHS) intersect, but it is impossible to get a closed form expression for $t_{j+1}$ or $\psi_{j+1}$ and similarly for $\dot{Z}_{j+1}$. Further details of the derivation of the equations for the maps can be found in Supplementary Section I.

For the BTB case, the same hurdle arises. In that case, the BTB motion is composed of maps $P_{TB}\circ P_{BT}$, and therefore a closed form first return map for $\partial B$ would require the composition of expressions for $(\dot{Z}_{j+1}, t_{j+1})$ and $(\dot{Z}_{j+2}, t_{j+2})$. The only difference in the equations for these quantities is the sign of $\pm d/2$ in \eqref{PBBt}, so the lack of closed-form expressions follows as in \eqref{PBBtj1}. Therefore, we propose a computer-assisted method to transform these non-autonomous, implicit maps into a composition of smooth maps using explicit polynomials. To achieve this, we define a novel type of return map that can be combined with phase plane analysis to identify regions of state space with potentially attracting or transient behavior.

 There are three key elements to our generalizable approach to the maps:
 \begin{enumerate}
     \item We exploit the non-smooth impact events in the dynamics, leading to the observation that any transient behavior can be broken down into a sequence of a small number of types of return maps to $\partial B$, as shown in Fig. \ref{model:abc}(b): those that impact $\partial T$ between sequential impacts on $\partial B$, and those that do not. 
 \item The second key element is the ability to approximate these return maps with polynomial functions, noting that there may be different choices for this approximation.
 \item We focus on families of return maps for which a valuable phase plane analysis follows naturally, in contrast to the maps between different impacts \eqref{relvar}-\eqref{basicmapdZ}. 
\end{enumerate}
\smallskip 

 With sequential impacts on $\partial B$ as a natural framework for defining the maps, we focus on the first return maps to $\partial B$ captured by $P_{\rm BTB}$ and $P_{\rm BB}$. Note that in order to capture all possible transients, one would normally have to consider sequences with multiple impacts on $\partial T$ before returning to $\partial B$, e.g., sequences such as $BTTB$, $BTTTB$, etc. While we could include these in our analysis, we note that for $\beta>0$ and for the values of the forcing amplitude $\hat{F}$ and restitution coefficient $r$ considered here, these sequences are generated in a limited range of larger initial  $\dot{Z}_k$ and nearly in-phase $\psi_k$. Furthermore, one can show by repeated applications of the maps that the larger values of $\dot{Z}_k$ can not be sustained for the given $\hat{F}$ and $r$ \cite{serdukova2023fundamental,DULIN2022}, so repeated impacts on $\partial T$ are highly transient for these parameters. Therefore, they play a negligible role in the global dynamics, particularly as we focus on potentially attracting regions. Some detailed comments on this are included under Remark \ref{range} below. 
 
 Our approach also allows for considering sequences such as BTBB and BBTB. As discussed in the Conclusions, these sequences correspond to grazing bifurcations to 2:1 solutions. Bifurcations to stable 2:1 behavior do not occur for the parameters considered here and can be demonstrated as transient, so they are not considered here.

\begin{remark1}\label{not_Poincare}
It is worth noting the distinction between this approach and that of a Poincar\'e map with $\partial B$ as the Poincar\'e section. Here, we divide the sequences that all return to  $\partial B$ into different families, depending on which other impacts occur before the system returns to $\partial B$, considering the maps for the different families separately. 
\end{remark1}
 
While above, we have used the subscripts $j$ and $k$ somewhat generically for impacts, for clarity with respect to the maps in \eqref{relvar}-\eqref{basicmapdZ}, we reserve the subscripts $j,j+1, \ldots$ for sequential impacts on either $\partial B$ or $\partial T$. Then, for the  sequential impacts on $\partial B$ only,  in the following we use the subscripts $k, k+1, \ldots$, so that for $k=j$ and $P_{\rm BTB}$ ($P_{\rm BB}$), the $(k+1)^{\rm th}$  impact on $\partial B$ corresponds to the $(j+2)^{\rm th}$ ( $(j+1)^{\rm th}$) impact when counting all impacts. That is, for  ${Z}_{j} \in \partial B$,
\begin{eqnarray} \label{exactmap}
 P_{\rm BTB}:& (\dot{Z}_j,\psi_j) \to &\{(\dot{Z}_{j+2},\psi_{j+2}) \ | \ {Z}_{j+1} \in \partial T, {Z}_{j+2} \in \partial B\}, \nonumber \\
 P_{\rm BB}:& (\dot{Z}_j,\psi_j) \to &\{(\dot{Z}_{j+1},\psi_{j+1}) |  \ {Z}_{j+1} \in \partial B\}.
 \end{eqnarray}
Note, for physical clarity, we have slightly abused notation in \eqref{exactmap}, using  $Z_j\in \partial B$ and $Z_j\in \partial T$ for impacts on either end of the capsule, in place of $Z_j = \pm d/2$ as discussed following \eqref{eq_impact}.

As illustrated in Fig. \ref{model:abc}(b),  the sequence length, for example, to (nearly) periodic behavior is not uniform over the space of initial conditions and cannot be anticipated {\it a priori}. The return map to $\partial B$ gives a flexible construction that can be applied over any length of the transient.
This framework is well-suited for capturing global dynamics through phase plane techniques and can also be applied in stochastic settings for the VI pair \cite{Kwanposter23}.  In identifying potentially attracting dynamics, we use projections of the return maps in the $\dot{Z}_{k} -\dot{Z}_{k+1}$ and $\psi_k - \psi_{k+1}$  phase planes, relative to the corresponding diagonals (see Section \ref{S3.1visualization}). The maps in \eqref{relvar}-\eqref{basicmapdZ} do not lend themselves to these goals, as these are not (necessarily) return maps. 

\begin{figure}[htbp]
\centering
 \begin{subfigure}[t]{0.5\textwidth}
\centering
\caption{}
\includegraphics[width=\textwidth]{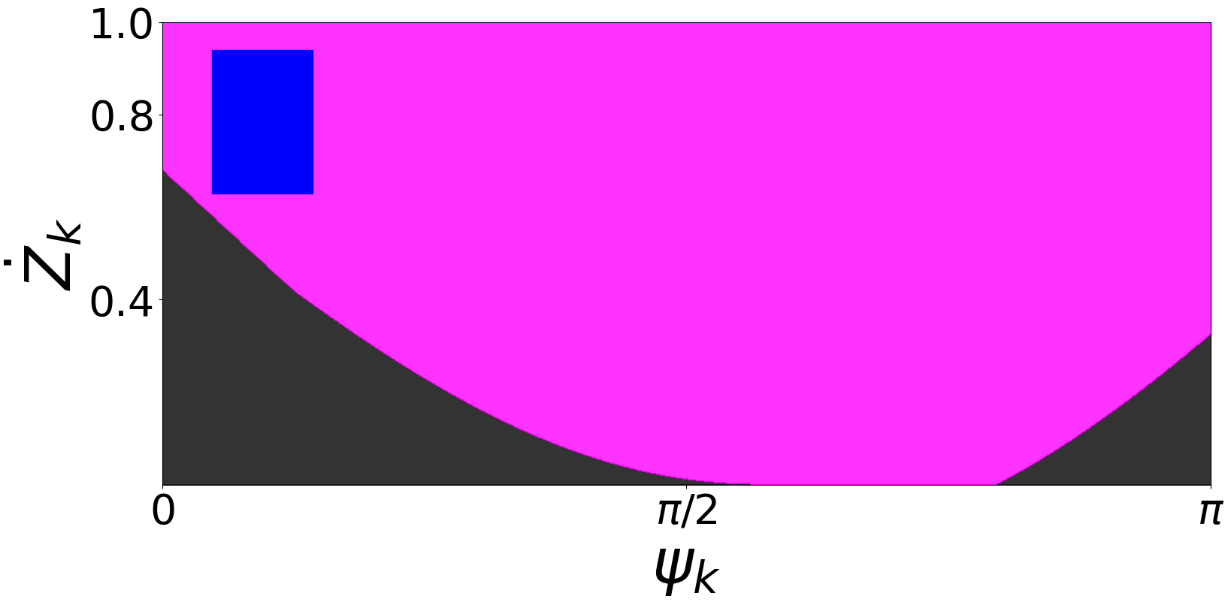}
 \end{subfigure}
     \hfill
 \begin{subfigure}[t]{0.4\textwidth}
\centering
\caption{}
\includegraphics[width=\textwidth]{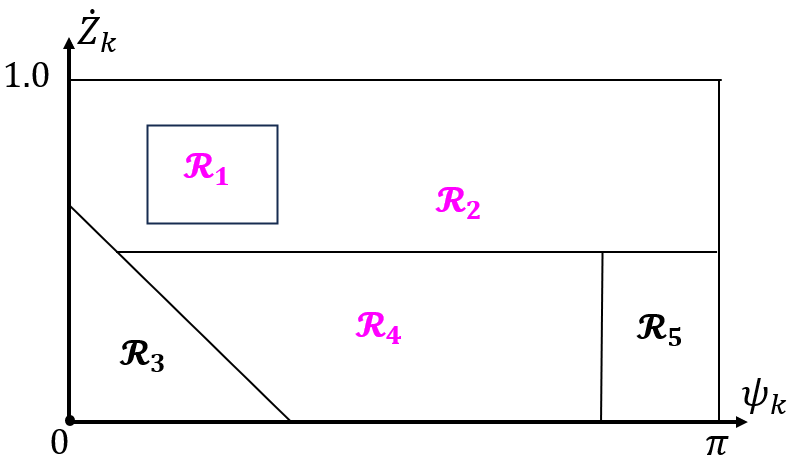}
 \end{subfigure}
\caption{(a): Using the building blocks in \eqref{exactmap}, the state space  $\dot{Z}_k-\psi_k$ can be partitioned based on two types of first return maps: $P_{\rm BB}$ (black regions) and $P_{\rm BTB}$ (magenta region). The blue square indicates the location of $\mathcal{R}_1$, a region within the $P_{\rm BTB}$ region that has special properties as studied in detail in Section \ref{S4algorithm}. (b): A  further partition of the state space into five regions, convenient for approximation as determined by the algorithm in Section \ref{S4algorithm}: Regions $\mathcal{R}_1$, $\mathcal{R}_2$, $\mathcal{R}_4$ divide the state space for the BTB motion, and Regions $\mathcal{R}_3$, $\mathcal{R}_5$ divide the state space for the BB motion. The partition in panel (b) shows an approximation to the exact solution in panel (a), so the dividing boundaries between regions do not match exactly those based on the exact map. The parameter used in (a) is $d=0.26$.}
\label{regions}
\end{figure}

 For the remainder of the paper, we track the first return maps for impact velocity and impact phase $(\dot{Z}_k,\psi_k)$ on $\partial B$, using the subscripts $k, k+1, \ldots$  to indicate sequential impacts on $\partial B$, composed of the building blocks in \eqref{exactmap}. Figure~\ref{regions} shows how the choice of these building blocks divides the state space for ($\dot{Z}_k$, $\psi_k$) by viewing this pair as the initial condition, which then yields one of these two return maps. Figure~\ref{regions}(a) shows how the ($\dot{Z}_k$, $\psi_k$) plane is divided by tracking the return maps. Figure \ref{regions}(b) illustrates a further division of the state space, necessary for applying straightforward polynomial approximations of the return maps, as discussed in the context of the full algorithm described in Section \ref{S4algorithm}.  Note that the building blocks \eqref{exactmap}  are analogous to short words in the symbolic representations used for piecewise linear maps in \cite{Simpson2023}, which form the basis for invariant cones and trapping regions. We note that the trapping region in \cite{Simpson2023} appears to be analogous to what we call the attracting domain in this paper, which is a compact region that attracts all non-trivial trajectories of the map.

\begin{remark1} \label{range}
   For the algorithm developed in this paper, we restrict our attention to the range of \, $ 0\leq\psi_k\leq \pi$, discussed further in the context of  Fig.~\ref{2Dpi-2pi} below.  As can be shown for the model \eqref{cyl_dim}-\eqref{simpact_dim} and the parameters considered in this paper, impacts with $\psi_k>\pi$ correspond to those where the ball and capsule are moving in the same direction, yielding smaller impact velocities and thus transient behavior in both $\psi_k$ and $\dot{Z}_k$ \cite{Serdukova2019}.  This point is discussed in Section \ref{S3.1visualization} below, in the context of projections of the 2D maps for $\dot{Z}_k, \psi_k$ into their corresponding phase planes. Likewise, for the parameter regimes considered in this paper, focusing on a range of $d$ with energetically favorable 1:1-type sequences of alternating impacts, the impact velocities in the range $\dot{Z}>1.0$ are transient.  Fig.~\ref{ICPP02pi} in Appendix \ref{Appendix-statespace} illustrates the additional regions with transient BTTB behavior, which can appear for $\dot{Z}>1.0$. While the approach proposed here can handle these values by including additional transient regions, for $\beta>0$ and the parameters considered here, these sequences are strongly transient and essentially negligible in the global behavior. Then for simplicity of exposition, we restrict our attention to $0\leq\psi_k\leq\pi$ and $0<\dot{Z} \leq 1.0$.   
\end{remark1}
  
\begin{figure}[htbp]
\centering
\includegraphics[width=\textwidth]{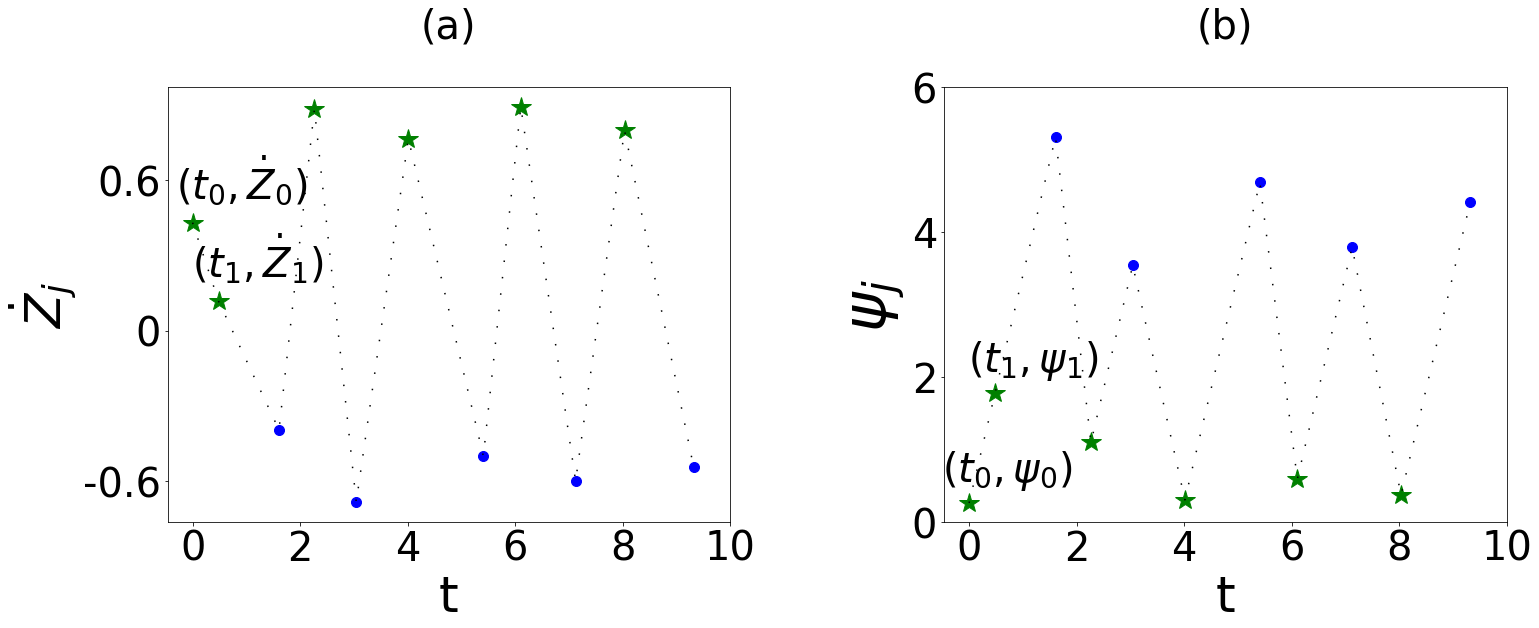}
\caption{ The values $(\dot{Z}_j,\psi_j)$ at impacts (both $\partial B$ (green stars) and $\partial T$ (blue circles)), starting with initial conditions $\dot{Z}_0 = 0.43$ and $\psi_0=0.26$ with $d=0.35$. Note that the location of the impact determines the sign of the relative velocity: $\dot{Z}_j>0$ for the impact on $\partial B$, and $\dot{Z}_j<0$ for  $\partial T$, and the dotted lines trace the order in which the impacts happen. In this paper, we focus on the return map for $\partial B$, denoted $(\dot{Z}_k, \psi_k)$.} 
\label{impact-time}
\end{figure}

Figure~\ref{impact-time} illustrates the reduction of our representation within the dynamics, focused on the impact velocity $\dot{Z}_j$ and phase $\psi_j$ on $\partial B$ (green stars),  in contrast to Fig.  \ref{model:abc}(b), which shows the exact behavior solution at and between the impact time.
The first return maps in \eqref{exactmap} are implicit in form and thus awkward to use directly in a global stability analysis.  Then, as a first step towards a more explicit approximation, we visualize the return maps in \eqref{exactmap}.

\subsection{Visualization of the maps and projections in the phase planes} \label{S3.1visualization}

\begin{figure}[htbp]
\centering
\includegraphics[width=0.9\textwidth]{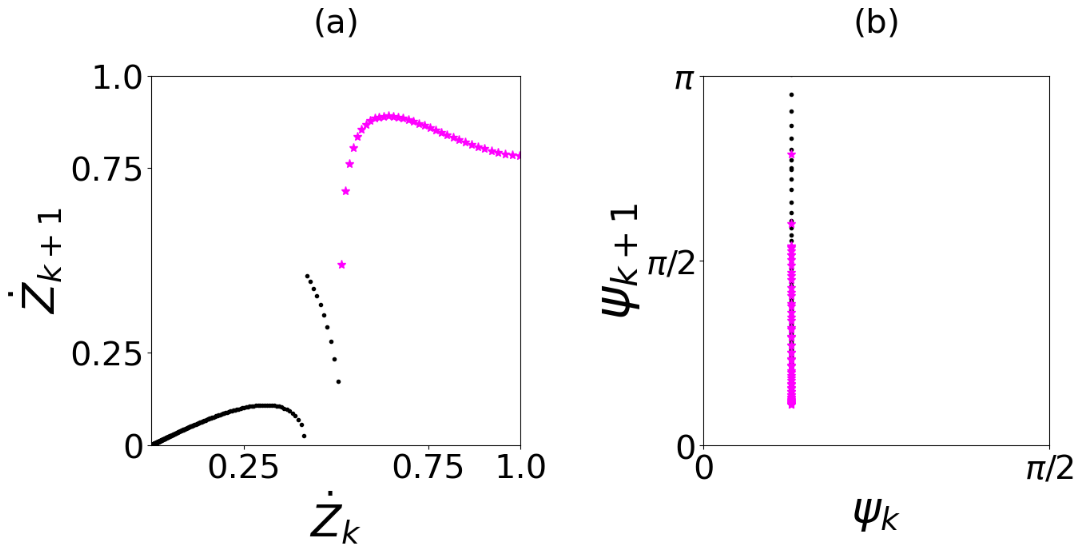}
\caption{Illustration of  $\dot{Z}_{k+1}$ and $\psi_{k+1}$, the first return maps on $\partial B$ using \eqref{exactmap} for fixed $\psi_k=0.4$ and sweeping through initial values $\dot{Z}_k \in (0,1.0)$ with $d=0.35$.  The magenta points correspond to the first returns via BTB type, and the black points represent the first returns of BB type. }
\label{2d-1phi}
\end{figure}

 Given that the return maps $P_{\rm BTB}, P_{\rm BB}$ are in terms of the  2D vector  ($\dot{Z}_{k},\psi_{k}$) we show two separate  surfaces for $\dot{Z}_{k+1}$ and $\psi_{k+1}$  generated by them. To build these up, we first show the maps projected in the phase planes $\dot{Z}_{k}- \dot{Z}_{k+1}$ and $ \psi_{k}-\psi_{k+1}$, for a fixed value of $0<\psi_k<\pi$, and sweeping through $\dot{Z}_k \in (0,1.0)$. In Fig. \ref{2d-1phi}(a), the resulting first return values $(\dot{Z}_{k+1}, \psi_{k+1})$ are sorted according to BTB and BB motion, as indicated by different colors.  In Fig. \ref{2d-1phi}(b), in this projection, these two types of behavior can interweave for a single value of $\psi_k$, as different values of $\dot{Z}_k$ yield a variety of $\psi_{k+1}$ that appear in both the $P_{\rm BTB}$ and $P_{\rm BB}$ return maps.
 
\begin{figure}[htbp]
    \centering
    \begin{subfigure}[t]{0.45\textwidth}
        \centering
    \caption{}
    \includegraphics[width=0.9\textwidth]{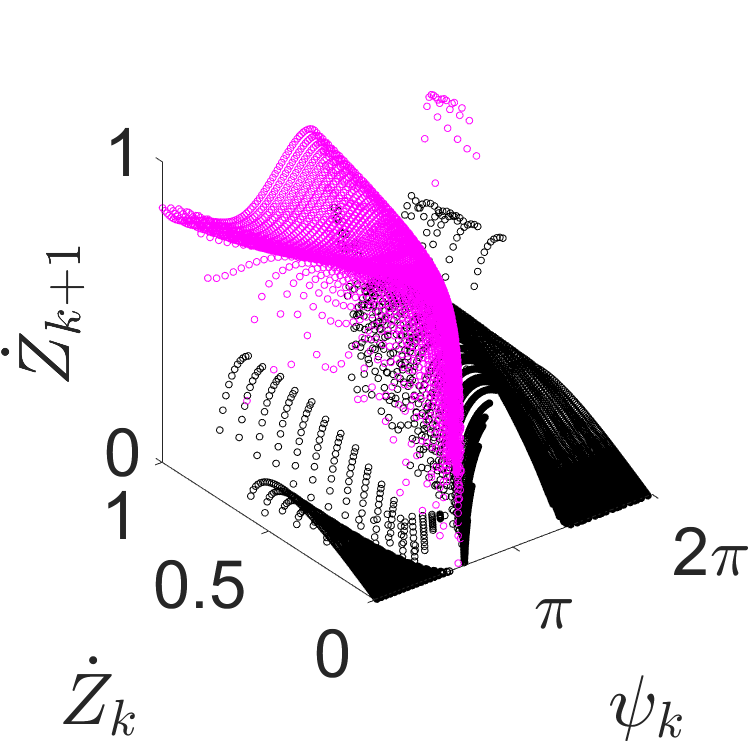}
    \end{subfigure}
    \hfill
    \begin{subfigure}[t]{0.45\textwidth}
        \centering
    \caption{}
    \includegraphics[width=0.9\textwidth]{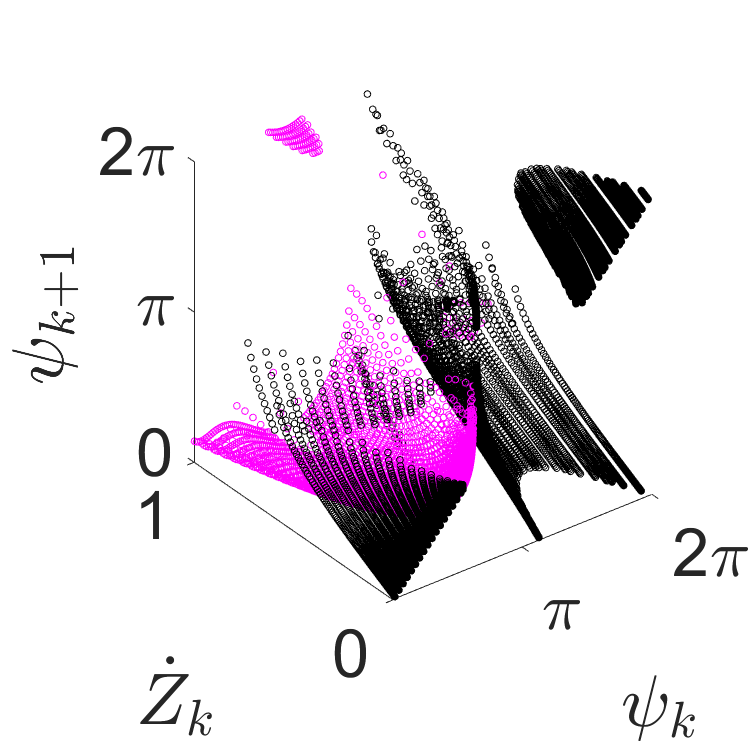}
    \end{subfigure}
\caption{Illustration of the 3D surfaces generated using the first return maps $P_{\rm BTB}$ (magenta) and $P_{\rm BB}$ (black) in \eqref{exactmap}, with $d=0.35$. Each initial condition pair $(\dot{Z}_k, \psi_k)$ has output $(\dot{Z}_{k+1}, \psi_{k+1})$, graphed on the vertical axes in panels (a) and (b), respectively. Supplementary Video 1 provides a $360^{\circ}$ rotating view of the surfaces. }
\label{3d}
\end{figure}

Repeating the application of the first return map \eqref{exactmap} over the range of initial phase values $\psi_k$ yields the surface visualized in Fig.~\ref{3d}, over a range of initial values in the horizontal $\dot{Z}_k-\psi_k$ plane. For $P_{\rm BB}$, shown by the black points,  in general  small values of $\dot{Z}_k$ (approximately $\dot{Z}_k<0.55$) map into small values of $\dot{Z}_{k+1}$, while $\psi_{k+1}$ tends towards values either near $0$ or above $2$. In the case of $P_{\rm BTB}$, shown by magenta points, larger $\dot{Z}_k$ map into larger values of $\dot{Z}_{k+1}$, with the corresponding $\psi_{k+1}$ spread out between $0$ and $\pi$. The visualization of the return maps $P_{\rm BB}$ and $P_{\rm BTB}$ indicates a few features that are important in approximating these surfaces with polynomial maps.  Not only are the surfaces disconnected, but the surfaces have dramatically different gradients corresponding to different regions in the $Z_k-\psi_k$ state space, which leads to the partitioning as shown in Fig. \ref{regions}(b). These regions are identified as part of the algorithm for approximating the surfaces, as discussed in detail in Section \ref{S4algorithm}.

\begin{figure}[htbp]
\centering
\includegraphics[width=0.8\textwidth]{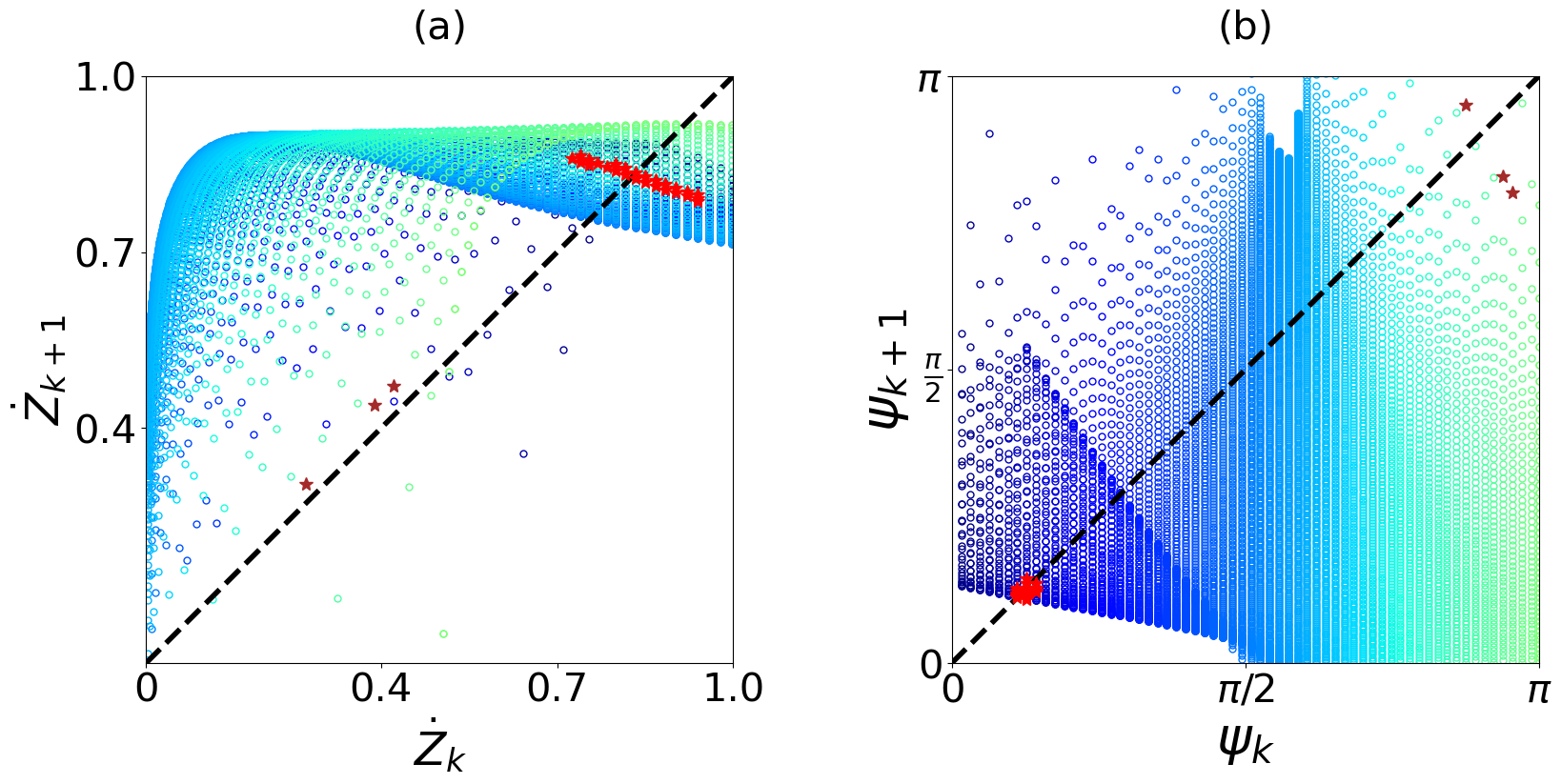}
\includegraphics[width=0.8\textwidth]{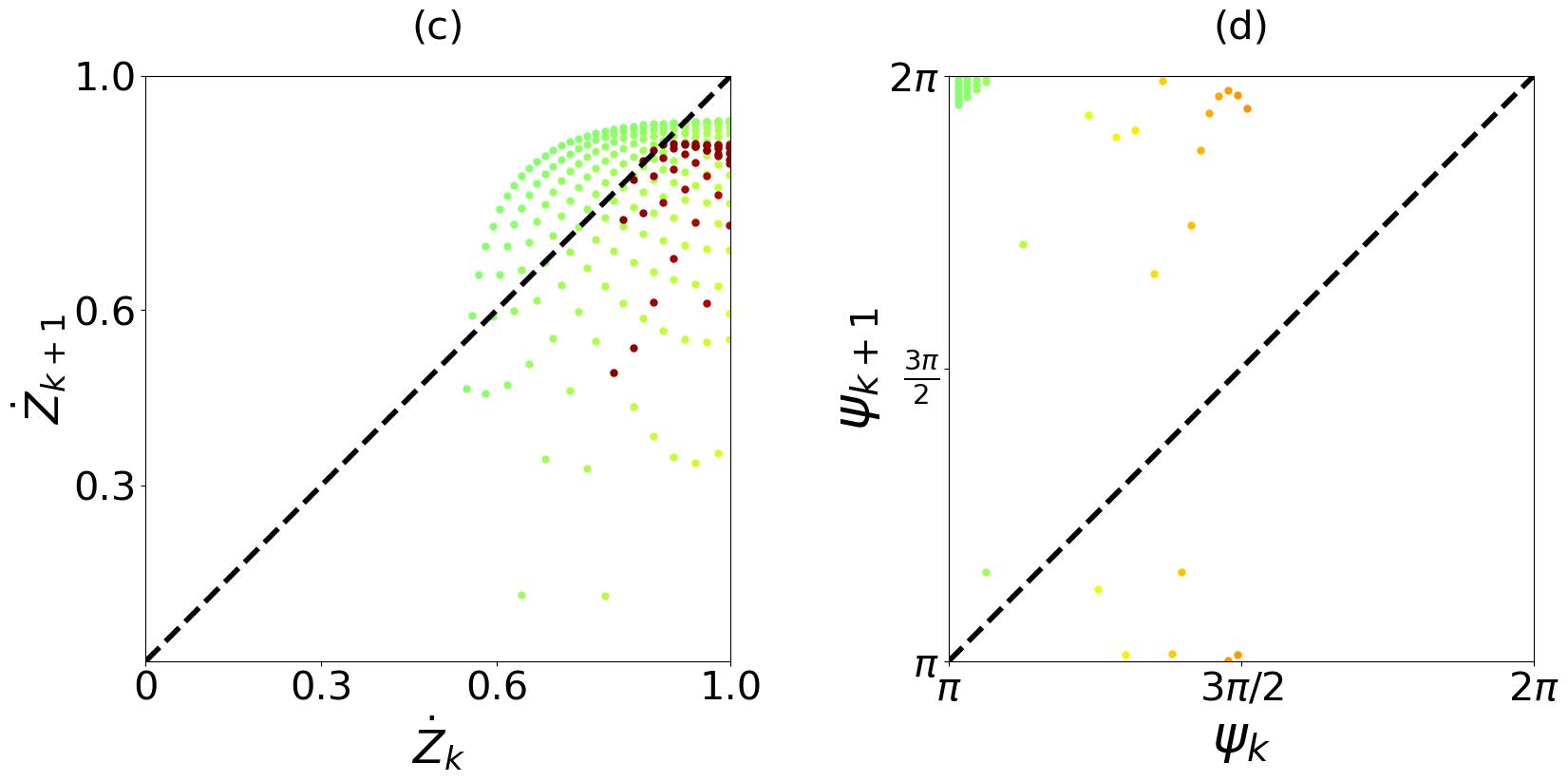}
\caption{The 2D projection of  the $P_{BTB}$ maps (magenta surfaces in Fig.~\ref{3d}) on the phase plane $\dot{Z}_{k}-\dot{Z}_{k+1}$ and $\psi_{k}-\psi_{k+1}$. For all panels, $d=0.35$. Different colors correspond to the maps for different values of $\psi_k$. Panels (a)-(b) show results for initial condition $\psi_k\in [0, \pi]$. Stars show cases where both maps take values near the diagonals in both phase planes; for red stars the slopes of the surfaces are smaller, suggesting potentially attracting dynamics near these values for $\psi_k<\pi/2$, while for brown stars the surfaces have steep slopes suggesting transient dynamics for these values when $\psi_k>\pi/2$. Panels (c)-(d) show results for initial condition $\psi_k\in [\pi, 2\pi]$; there are no cases where both maps take values near the diagonals, indicating transient dynamics over this range.}
\label{2Dpi-2pi}
\end{figure}

Comparison of the return maps with the diagonals in the $\dot{Z}_{k} - \dot{Z}_{k+1}$ and $\psi_k - \psi_{k+1}$ phase planes is achieved via projections of the return map surfaces on the phase planes, as shown in Fig. \ref{2Dpi-2pi} and in Appendix \ref{AppendixPrjt}, Fig. \ref{2D0-2pi}.  This projection is valuable as we identify potential regions for attracting and transient behaviors, following from comparisons of the map surfaces with the diagonals in the phase planes.
  
\begin{remark1}\label{diagonals}
To see the significance of the diagonals in the phase planes, recall the classic example of the logistic map $x_{n+1}=rx_n (1-x_n)$. The dynamics of the logistic map vary with the parameter $r$, directly related to the slope of the map $rx_n (1-x_n)$ at the fixed point $x_n^* =rx_n^* (1-x_n^*)$, which by definition is at the intersection of the phase plane diagonal and the map.  The fixed point is an attractor (repeller) if the absolute value of its slope is less than (greater than) 1.  This fact motivates us to look for potential attracting or transient dynamics by studying the intersection between the projections of the maps $\dot{Z}_{k+1}(\dot{Z}_k,\psi_k)$ and $\psi_{k+1}(\dot{Z}_k,\psi_k)$ in their respective $\dot{Z}_k-\dot{Z}_{k+1}$ and $\psi_k-\psi_{k+1}$ phase planes and the diagonals in those phase planes.
\end{remark1}

 Figure \ref{2Dpi-2pi} illustrates this comparison for the surfaces in the BTB region with the diagonals in the phase planes. There the surfaces, projected into the phase planes, are shown with different colors corresponding to different values of $\psi_k$, i.e., a different color for each ``strand'' in the map for fixed $\psi_k$ and sweeping over $\dot{Z}_k$ as in Fig. \ref{2d-1phi}; e.g., the value of $\psi = \pi/2$ is medium blue in both panels. Together, these form the complete surfaces for $\dot{Z}_{k+1}$ and $\psi_{k+1}$ then projected into their respective phase planes. Then, we look for cases where the same color strands cross the diagonals in both phase planes. These indicate potential fixed point values of $(\dot{Z}_k,\psi_k)$. To identify regions that contain these potentially attracting values, we look for regions where both maps have points near the diagonals of $\psi_k-\psi_{k+1}$ and $\dot{Z}_k-\dot{Z}_{k+1}$.
 Figure \ref{2Dpi-2pi} (a)-(b) shows these values for the  $P_{BTB}$ map for $0<\psi<\pi$,  with these points marked near the diagonals in both phase planes. There are two clusters of these points: red for those with $\psi_k<\pi/2$ and brown otherwise. Section \ref{S4.2} Iteration 2, step iii) in the algorithm below discusses the specific criteria for defining values near the diagonals, which yields the special region marked in blue in Fig. \ref{regions}(a) as a potentially attracting region.
Similarly,  Fig. \ref{2Dpi-2pi} (c)-(d) shows the $P_{BTB}$ maps for $\pi<\psi<2\pi$. While these points may satisfy the criteria for being near the diagonals, the steep slope of the curves forming the map for $\dot{Z}_{k+1}$ leads us to conclude that these points are not in a potentially attracting region, as repeated in Section \ref{S4.2} Iteration 2, step iii). 

Appendix \ref{AppendixPrjt}, Fig. \ref{2D0-2pi}, shows these comparisons of the maps projected into the phase planes for the $P_{BB}$ maps. The results of this comparison are discussed both there and in Section \ref{S4.2} Iteration 1, step iii), leading to conclusions about the transient nature of these regions. 
Section \ref{S4algorithm} further articulates these and other details in the application of the algorithm, combining visualizations of  Figs. \ref{3d}, \ref{2Dpi-2pi}, \ref{2D0-2pi},  and \ref{ICPP02pi} to give further insight into behavior on subdivisions of the return map surfaces together with approximating these surfaces with polynomials.

 \section{Composition of the Approximate Map}\label{S4algorithm}
 
We provide an algorithm for deriving a set of explicit piecewise polynomial maps $f_n$ and $g_n$  for each region ${\cal R}_n$ in the state space $\dot{Z}_k-\psi_k$, approximating the surfaces $\dot{Z}_{k+1}$ and $\psi_{k+1}$ as shown in Fig. \ref{3d}. The approximate return maps are given in terms of the variables $(v_k,\phi_k)$ that denote the approximate relative impact velocity on $\partial B$ and the corresponding impact phase, respectively, at the $k^{\rm th}$ return to $\partial B$.  We define the composite approximate map ${\cal M}$ that combines the continuous maps $f_n, g_n$ for the regions ${\cal R}_n$ in  Fig. \ref{regions}(b),  taking the form
\begin{eqnarray} \label{compositemapM}
(v_{k+1}, \phi_{k+1}) = {\cal M}(v_k,\phi_k) 
\equiv (f_n(v_{k},\phi_k),g_n(v_{k},\phi_k)), \mbox{ where } 
(v_k,\phi_k) \in {\cal R}_n.
\end{eqnarray} 

Given the complex nature of the surfaces for $\dot{Z}_{k+1}$ and $\psi_{k+1}$, the algorithm for constructing the maps $(f_n, g_n)$,  leads to refining the regions shown in Fig. \ref{regions}(a), resulting in the regions ${\cal R}_n$ for $n=1,2,3,4,5$ in Fig.~\ref{regions}(b).

\begin{figure}[htbp]
\centering
\begin{subfigure}[b]{0.45\textwidth}
\centering
\caption{}
\includegraphics[width=\textwidth]{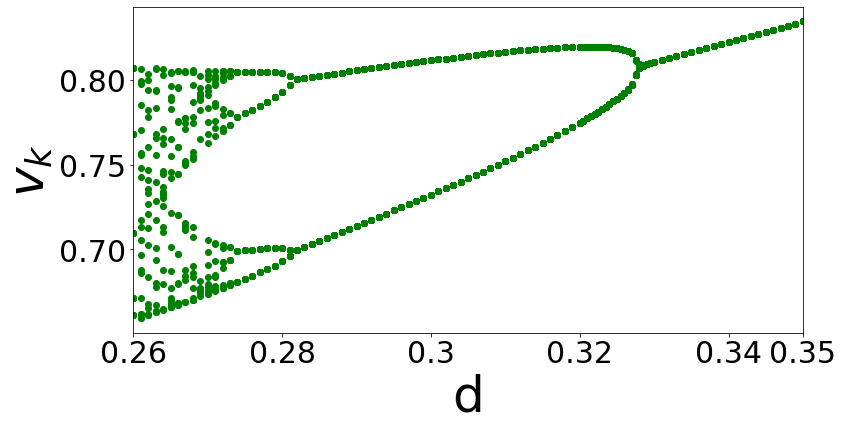}
     \end{subfigure}
     \hfill
     \begin{subfigure}[b]{0.45\textwidth}
         \centering
         \caption{}
\includegraphics[width=\textwidth]{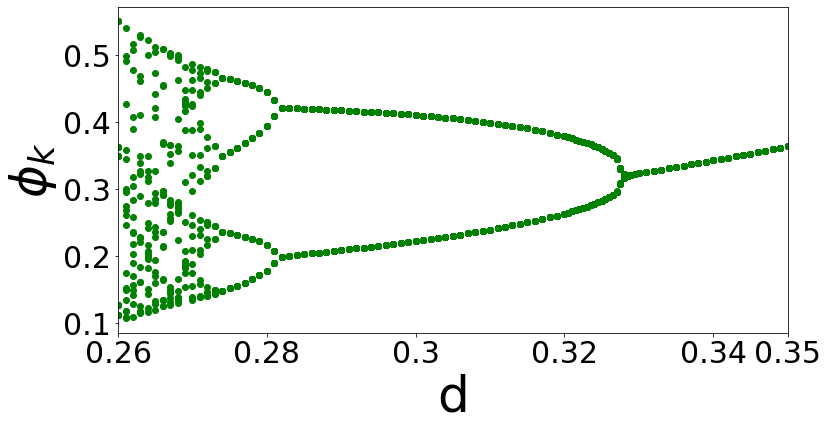}
     \end{subfigure}
\caption{Bifurcation diagrams generated using the composite approximate map $\mathcal{M}$, defined in \eqref{compositemapM} and Appendix \ref{algo_code}, with coefficients given in Supplementary Section II. The bifurcation structure obtained using $\mathcal{M}$ reproduces remarkably well that obtained for the exact map \eqref{relvar}-\eqref{basicmapdZ} presented in Fig. \ref{bifexact}.}
\label{bifapprox}
\end{figure}

Before constructing ${\cal M}$ in \eqref{compositemapM} (derived below, with specifics given in Appendix \ref{algo_code}), we give a first illustration that it captures the critical features of \eqref{relvar}-\eqref{basicmapdZ}  in the parameter range of interest, using it to obtain the bifurcation diagram analogous to Fig. \ref{bifexact}. Figure~\ref{bifapprox} shows the results for $v_k,\phi_k$ vs. $d$, generated using ${\mathcal M}$  via the continuation-type method described in Remark \ref{num_bif}. Comparing with the corresponding bifurcation diagram for the exact map in Fig. \ref{bifexact}, we see that the results from ${\mathcal M}$ capture a number of features of the original system, including $d$ values for the period-doubling bifurcations, the attracting values of $v_k$ and $\phi_k$ for the different branches, and the approximate range of values of $v_k$ and $\phi_k$ for the chaotic behavior obtained for smaller $d$ in the range shown in Figs. \ref{bifexact}, \ref{bifapprox}.

\smallskip

\subsection{General Algorithm: Construction of the composite map $\cal{M}$} \label{S4.1}

\begin{figure}[htbp]
\centering
\includegraphics[width=\textwidth]{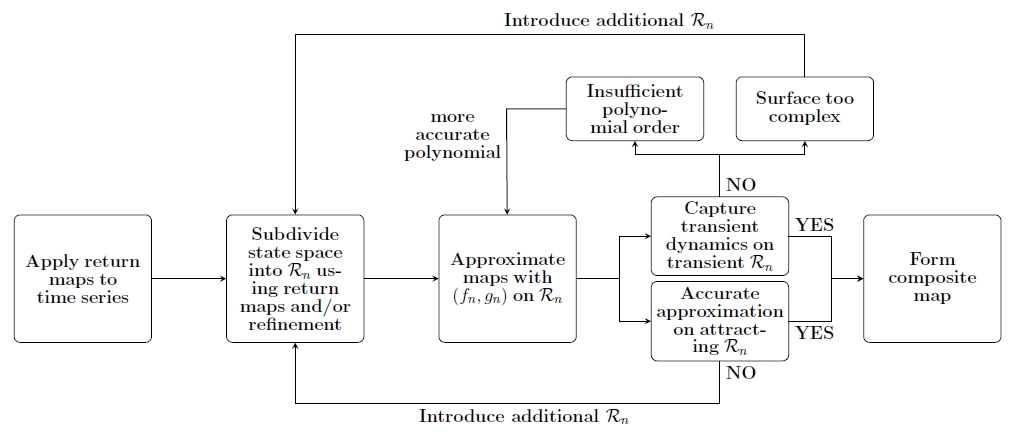}
\caption{Illustration of the general algorithm for constructing the composite map. }
\label{flochart}
\end{figure}

Illustrated in Fig. \ref{flochart}, the general algorithm consists of three main activities: identifying an initial partition of the state space based on the return map building blocks, iterating on approximations of the return maps on these regions, and including updates of the regions as necessary to improve the critical approximations.\\

\noindent
\underline{\bf Initialize: steps 0)-ii):}  Partition state space for the definition of the composite map.\\
0). Choose a state as the basis for return behavior.\\
i). Generate surfaces $(\dot{Z}_{k+1},\psi_{k+1})$ corresponding to the different families of the first return maps for this state;\\
ii). Partition regions in the state space based on these 
 different families of first returns.  Label these regions as  ${\cal R}_{n.1}$, denoting Region $n$ defined on  iteration 1.\\
 ~\\
\underline{\bf Iterate on steps iii)-vi)} until appropriate fit for surfaces corresponding to first return map for all regions ${\cal R}_{n.m}$, region $n$ on $m^{\rm th}$ iteration.\\
iii). Identify potential regions of attraction or transient behavior, e.g. via comparison of the surfaces with diagonals in the phase planes;
iv). Choose an appropriate order of polynomial fit for the surface(s) in each region, via testing different orders of polynomials and, depending on the resolution needed, to identify $f_n$ and $g_n$ for each ${\cal R}_{n.m}$.\\
v). If the fit of the polynomial is unsatisfactory, adjust 
 the size of the regions and/or locate new regions for additional partitions.\\
vi). Optional reduction: for regions that yield immediate transitions to other regions, replace with appropriate resetting conditions.\\
~\\
\underline{\bf Finalize}\\
vii). Once the polynomial approximations are defined for maps for all regions, finalize definitions of regions, labeled as ${\cal R}_{n}$, dropping the $.m$ label, together with their corresponding maps $(f_n,g_n)$. This final step includes a definition of the range for each map, as discussed further in the demonstration in Section \ref{S5validation}. \\

Steps iii)-vi) depend on the computer-assisted analysis of several different features of the first return map surfaces found in ii). In iii) we compute quantities relevant to the dynamics and geometric characteristics of the maps as we make comparisons with the corresponding phase planes (see Remark \ref{diagonals}). Keeping the user-defined constants to a minimum, we must define a level of accuracy and order of the polynomial when fitting any of the surfaces in step iv). There is also some flexibility in the size of the region defined for the potentially attracting region(s), as used in step v). 
In Section \ref{S4.2}, we illustrate the implementation of these steps and parameter choices in the concrete context of \eqref{cyl_dim}-\eqref{simpact_dim} and the corresponding non-dimensional form \eqref{eq_impact}.  There, we also highlight the points about computational efficiency and adjustments of any user-defined parameters related to the accuracy of the polynomial approximations.

\begin{remark1}\label{separable}
 As demonstrated below, in certain regions ${\cal R}_n$ where the shape of the map clearly indicates transient dynamics, we look for a simple approximation that takes the form of a single variable polynomial for each of the variables of interest, e.g., $v_{k+1}= f_n(v_k)$ and $\phi_{k+1}= g_n(\phi_k)$.  We refer to these as separable maps since we approximate the 2D map for $(v_k,\phi_k)$ with two 1D maps that each depend on a single variable. Such an approximation supports a cleaner visualization in the phase plane by simplifying the details of the transient behavior while approximating it as dictated by the shape of the exact map.
\end{remark1}

\subsection{\bf Algorithm implementation: a composite map for the VI pair model}\label{S4.2}

We apply the general algorithm outlined above  - Initialize, Iterate, and Finalize - to identify appropriate partitions of the state space and the approximations for the return maps on these regions for the non-dimensionalized VI pair model as in  \eqref{relvar}.  Here, we present this application step-by-step, with the specific details of the composite map ${\mathcal M}$ given in Appendix \ref{algo_code}. \\
~\\

\noindent
\underline{ \bf Initialize} the partition of the state space.\\
0). Choose $Z\in \partial B$ as the state for the basis of the first return maps.\\
i). Generate surfaces $\dot{Z}_{k+1}$ and $\psi_{k+1}$ for BTB and BB behavior as first return maps \eqref{basicmapdZ} over the range of possible initial conditions in the state space ($\dot{Z}_{k}, \psi_{k}$) (see, e.g., Fig. \ref{regions}(a)). \\
ii). Partition the state space into regions ${\cal R}_{n.1} $ according to these building blocks: BTB and BB: ${\cal R}_{1.1}$ corresponds to  BTB,  ${\cal R}_{3.1}$ corresponds to BB behavior for smaller $\psi_k$, and ${\cal R}_{5.1}$  corresponds to BB behavior with larger $\psi_k$.\\

\noindent
\underline{\bf Iteration 1:} steps iii)-vi)
\renewcommand{\labelenumi}{\roman{enumi}}
\begin{enumerate}
\item[iii).] Identify regions of potential attraction and transients as follows.
\begin{itemize}
\item ${\cal R}_{1.1}$:  entire region of BTB behavior, including both transient regions and potential attracting dynamics near the diagonals in the $\dot{Z}_{k} - \dot{Z}_{k+1}$  and  $\psi_k - \psi_{k+1}$  planes.
\item  ${\cal R}_{3.1}$: The surfaces for BB behavior with sharp gradients in the map near the diagonals. Thus, transient BB behavior is expected (see Fig. \ref{2D0-2pi}).
\item  ${\cal R}_{5.1}$:  The surfaces for  BB behavior are away from the diagonal in the $\psi_k$-$\psi_{k+1}$ plane, thus transient BB behavior is expected (see Fig. \ref{2D0-2pi}).
\end{itemize}
\item[ iv).] Polynomial approximation of surfaces  for  $\dot{Z}_{k+1}$ and $\psi_{k+1}$ in ${\cal R}_{1.1}$,
${\cal R}_{3.1}$, and
${\cal R}_{5.1}$ (see Fig. \ref{3d}):
\begin{itemize}
\item ${\cal R}_{1.1}$, BTB behavior: The surface in this region is a combination of subregions where the surfaces for $\dot{Z}_{k+1}$ and $\psi_{k+1}$ have more gradual variation, contrasted with others with sharp gradients. Thus, an accurate polynomial fit is challenging, which also limits an accurate approximation of potentially attracting dynamics near the diagonals in the $\dot{Z}_{k} - \dot{Z}_{k+1}$  and  $\psi_k - \psi_{k+1}$ phase planes. This motivates a further partitioning of the BTB region, as described in step v).
\item ${\cal R}_{3.1}$, BB behavior: As can be observed in Fig. \ref{3d},  there are two disjoint surfaces for $\dot{Z}_{k+1}$.  One is a curved surface with sharp gradients for which we use fifth/fourth order polynomials in $v_k$/$\phi_k$ for the approximate map $(f_3, g_3)$ (see Appendix \ref{appendixR3}).
There is a second segment, nearly vertical in $\dot{Z}_{k+1}$, discussed in (vi) below. 
\item  ${\cal R}_{5.1}$:  As the surfaces for $\dot{Z}_{k+1}$ and $\psi_{k+1}$ in ${\cal R}_{5.1}$ are away from the diagonal, we use  a ``separable'' approximation, as discussed in Remark \ref{separable}.
See  Appendix \ref{AppendixR5}  for a discussion of the resulting approximate map 
$(f_5,g_5)$.
\end{itemize} 
\item[ v).] Update regions in terms of additional partitions for ${\cal R}_{1.1}.$ The different features of the $\dot{Z}_{k+1}$ and $\psi_{k+1}$ surfaces in ${\cal R}_{1.1}$ motivates sub-dividing into two regions:
 \begin{itemize}
\item ${\cal R}_{1.2}$: identify potentially attracting states, e.g. states for which the repeated images of the return map $P_{BTB}$ are near the diagonals in the $\dot{Z}_{k} - \dot{Z}_{k+1}$  and  $\psi_k - \psi_{k+1}$ phase planes. This choice of ${\cal R}_{1.2}$ limits to cases where the slopes of the surfaces near the diagonals are primarily small, e.g., less than unity for some values of $d$. The details for defining ${\cal R}_{1.2}$ are given in Iteration 2, step iii), including a quantitative characterization of proximity to the diagonals.
 \item ${\cal R}_{2.2}$: the remaining  states that produce clearly transient BTB behavior. This region includes sections of the $P_{BTB}$ map located away from the phase plane diagonals and sections near the diagonals with sharp gradients.
\end{itemize}
\item[vi).] From physical considerations, some maps are replaced with resetting functions and/or approximate maps in nearby regions.
\begin{itemize}
\item $\pi<\phi<2\pi$: The transient behavior for this range of $\phi_k$ is discussed in Remark \ref{range} above. Then, we employ the reset: $\phi_{k+1}=1.2$ and $v_{k+1}=v_k$ if $\phi_k>\pi$ or $\phi_k<0$ (see Appendix \ref{algo_code}). The results are not sensitive to the choice of the user-supplied reset value of $\phi_{k+1} = 1.2$. The shape of the surfaces in Fig. \ref{3d}, consistent with observations from other studies \cite{DULIN2022,serdukova2023fundamental}, suggests that the system moves towards values  $\phi<\pi/2$. At the same time, for the sake of generality, we want to choose a value in a transient region consistent with the definition of ${\cal R}_{2.2}$ obtained in Iteration 2 below.
\item  The nearly vertical surface in ${\cal R}_{3.1}$ mentioned above represents strongly transient behavior, consisting of transitions to BTB behavior or other states in ${\cal R}_3$. This transient behavior is captured by using equations \eqref{Apr3eqn} throughout ${\cal R}_{3.1}$, without approximating the vertical surface. Likewise, there is a small vertical section of the surface $\psi_{k+1}$ in
${\cal R}_{5.1}$, also discussed in Appendix \ref{AppendixR5}.
\end{itemize}
\end{enumerate}

\smallskip

\noindent
\underline{\bf Iteration 2:} steps iii)-vi)\\
 Iteration 2 is focused on the newly defined ${\cal R}_{1.2}$ and ${\cal R}_{2.2}$. 
 
\begin{enumerate}

\item[iii).] Considering attracting and transient BTB behavior:
\begin{itemize}
\item   To identify ${\cal R}_{1.2}$ as described in Iteration 1 step v),  we introduce a filter ${\cal R}_{1.2}(d)$  for a given $d$ that selects states ($\dot{Z}_{k}, \psi_{k}$) near the diagonals ($\dot{Z}_{k}, \psi_{k}$)  in the $\dot{Z}_{k} - \dot{Z}_{k+1}$  and  $\psi_k - \psi_{k+1}$ phase planes with images ($\dot{Z}_{k+1}, \psi_{k+1}$) from $P_{BTB}$  near the same diagonals. We then take the union of these regions to obtain a region valid for the full range of $d$ of interest. Then, $\mathcal{R}_{1.2}$ is given by 
\begin{align}
& \mathcal{R}_{1.2}(d)  = \Big\{(\dot{Z}_k, \psi_k)\; : \; \frac{1}{\delta} < \left|\frac{\psi_{k+1}}{\psi_k}\right|<\delta \;{\rm  and}\;\frac{1}{\delta} < \left|\frac{\dot{Z}_{k+1}}{\dot{Z}_k}\right|<\delta \;\;  \Big\}, \nonumber\\
&\mathcal{R}_{1.2} = \cup_{d\in [0.26, 0.35]}
\mathcal{R}_{1.2}(d).
    \label{r1uniond}
\end{align}
 Of course, the size of ${\cal R}_{1.2}$ depends on the choice of the user-defined parameter $\delta$, which characterizes proximity to the diagonals. As discussed further in Appendix \ref{appendix-R1}, the choice of $\delta$, together with the choice of polynomial order, influences the error of the approximation of the surface in the region ${\cal R}_{1.2}$. Figure \ref{3d-area} shows an example of the definition of ${\cal R}_{1.2}$.

\begin{figure}[htbp]
\centering
\begin{subfigure}[t]{0.45\textwidth}
    \centering
    \caption{}
    \includegraphics[width=\textwidth]{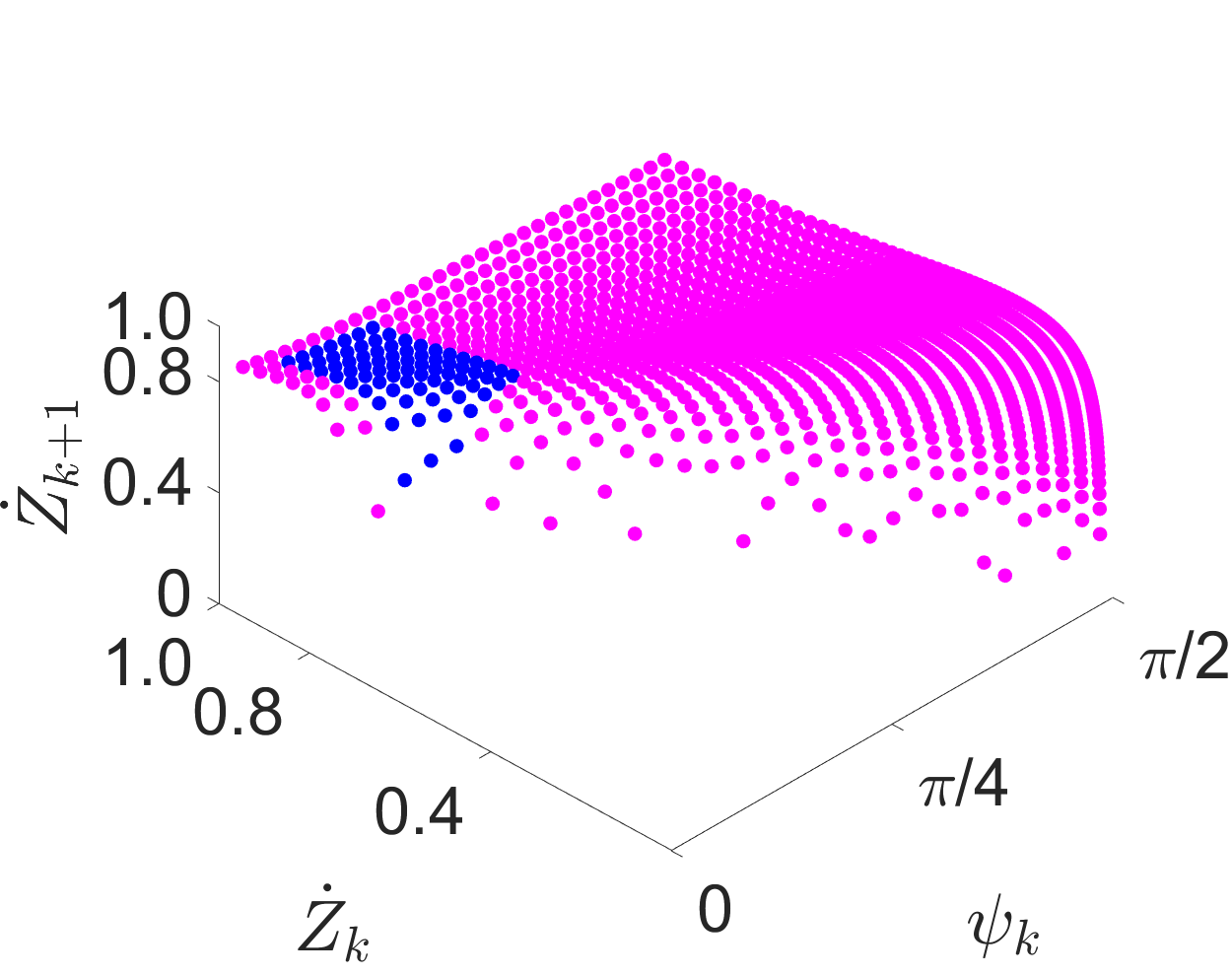}
\end{subfigure}
\hfill
\begin{subfigure}[t]{0.45\textwidth}
    \centering
    \caption{}
    \includegraphics[width=\textwidth]{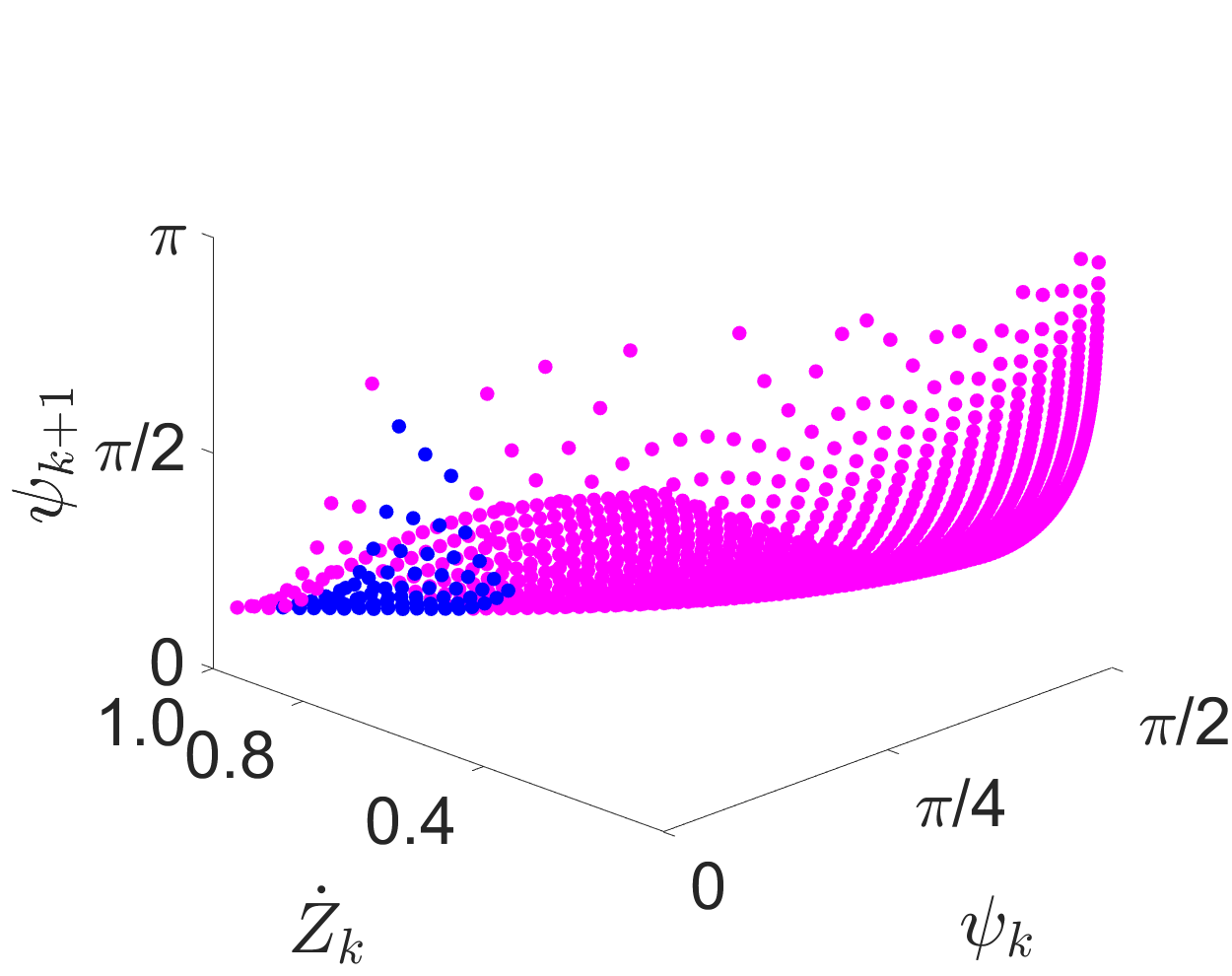}
\end{subfigure}
\caption{ The surface corresponding to $P_{\rm BTB}$ (magenta and blue combined), with $d=0.35$,  where  $\mathcal{R}_{1.2}$ (blue region), is obtained by using the filter \eqref{r1uniond} ($\delta=1.2$) to identify return maps located near diagonals in both the $\dot{Z}_{k+1}$ - $\dot{Z}_k$ and $\psi_{k+1}$ - $\psi_k$ phase planes.
}
\label{3d-area}
\end{figure}

 \item ${\cal R}_{2.2}$ is defined as the remainder of the BTB region, with transient behavior. This conclusion follows from Fig. \ref{2Dpi-2pi}, where the remainder of $(\dot{Z}_k,\phi_k)$ in the BTB region either do not correspond to points near the diagonals in both $\dot{Z}_{k+1}-\dot{Z}_k$ and $\psi_{k+1}-\psi_k$ phase planes, and/or are located on parts of the surfaces with steep slopes.
 \end{itemize}
\item[ iv).] Polynomial approximation of surfaces $\dot{Z}_{k+1}$ and $\psi_{k+1}.$
\begin{itemize}
 \item  ${\cal R}_{1.2}$: To capture subtle changes in the attracting behavior near the diagonals, the surfaces for $\dot{Z}_{k+1}$ and $\psi_{k+1}$ are approximated with polynomials of degree 3 in $v_k$ and degree 2 in $\phi_k$, 
\begin{align} 
  v_{k+1}(v_k, \phi_k) &= f_1(v_k,\phi_k) \nonumber\\
   &=b_0+b_1 \phi_k + b_2 v_k + b_3 \phi_k^2 + b_4 \phi_k v_k + b_5 v_k^2 + b_6 \phi_k^2 v_k + b_7 \phi_k v_k^2 + b_8 v_k^3,\label{cubic2dv} \\
        \phi_{k+1}(v_k, \phi_k) &= g_1(v_k,\phi_k)  \nonumber\\
     &=a_0+a_1 \phi_k + a_2 v_k + a_3 \phi_k^2 + a_4 \phi_k v_k + a_5 v_k^2  + a_6 \phi_k^2 v_k + a_7 \phi_k v_k^2 + a_8 v_k^3. \label{cubic2dp}
\end{align}

\item  ${\cal R}_{2.2}$:  We use a ``separable'' approximation (see Remark \ref{separable})  that takes the form
\begin{align}\label{r2map}
    v_{k+1}(v_k) &= f_2(v_k) = b_{20} v_k^5 + b_{21} v_k^4 + b_{22} v_k^3 + b_{23} v_k^2 + b_{24} v_k + b_{25},\nonumber\\
     \phi_{k+1}(\phi_k) &= g_2(\phi_k) = a_{20} \phi_k^5 + a_{21} \phi_k^4 + a_{22} \phi_k^3 + a_{23} \phi_k^2 + a_{24} \phi_k + a_{25}.
\end{align}
 Figure~\ref{r2approx}(a)-(c) shows  (green) curves representative of the transient behavior for this region, following from the shape of the surfaces for $\dot{Z}_{k+1}$ and $\psi_{k+1}$ shown in panel c) for  ${\cal R}_{2.2}$. The orange curves, showing the separable map in \eqref{r2map}, approximates this green curve. See further discussion in Appendix \ref{appR2}.
 \end{itemize}

 \begin{figure}[htbp]
 \centering
 \begin{subfigure}[t]{0.22\textwidth}
     \centering
     \subcaption{}
     \includegraphics[width=\textwidth]{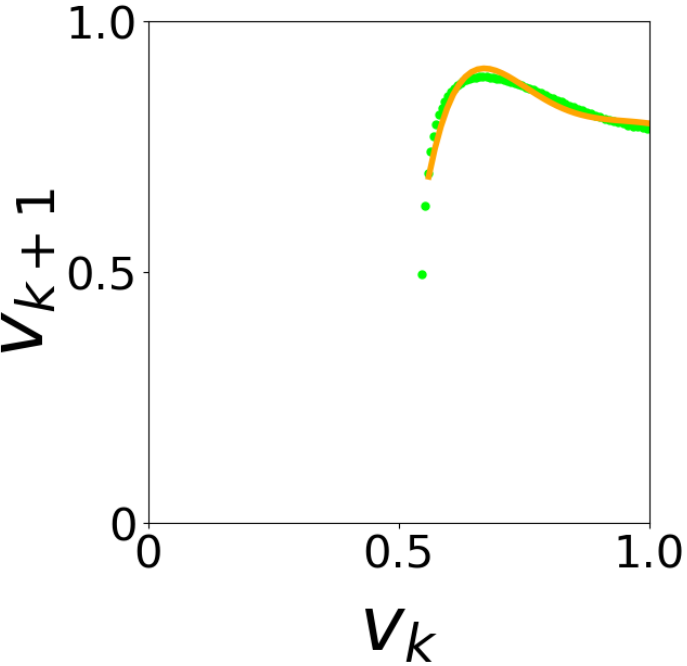}
 \end{subfigure}
 \hfill
 \begin{subfigure}[t]{0.21\textwidth}
     \centering
     \subcaption{}
     \includegraphics[width=\textwidth]{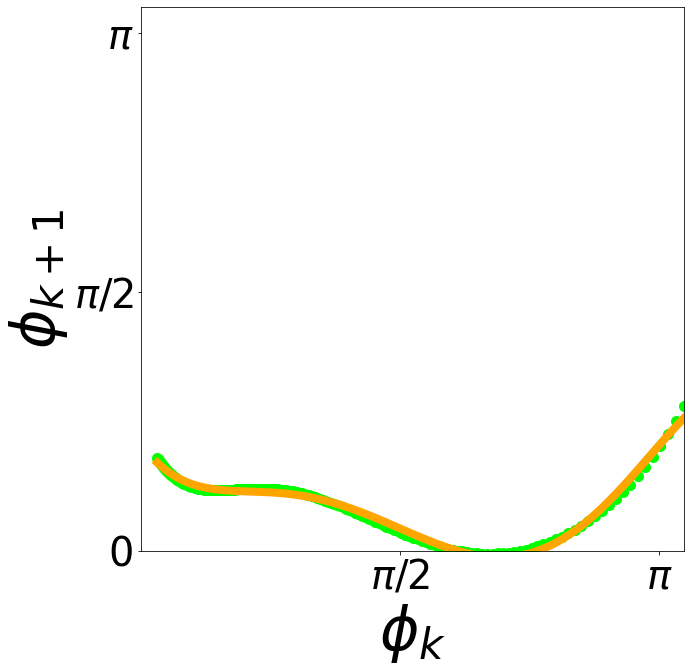}
 \end{subfigure}
  \hfill
 \begin{subfigure}[t]{0.5\textwidth}
     \centering
     \subcaption{}
     \includegraphics[width=\textwidth]{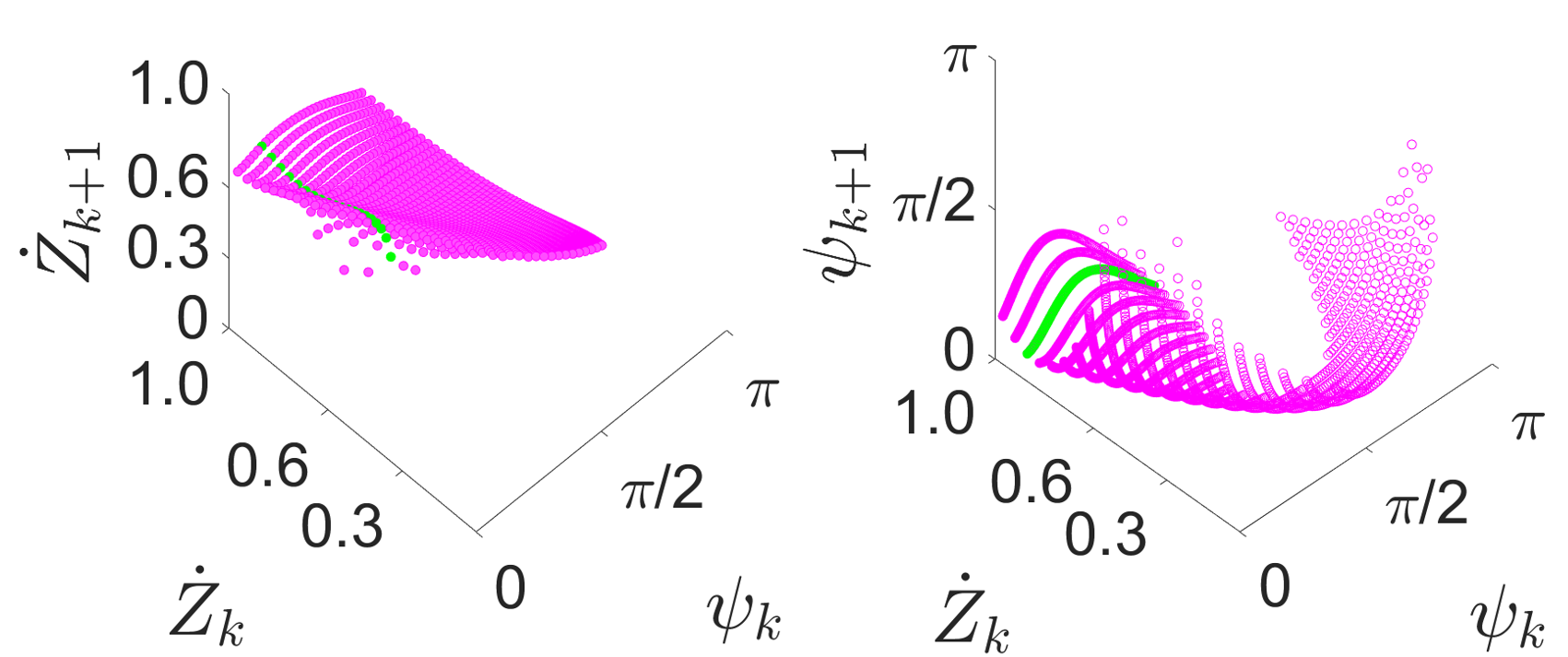}
 \end{subfigure}
 \vfill
 \begin{subfigure}[t]{0.22\textwidth}
     \centering
     \subcaption{}
     \includegraphics[width=\textwidth]{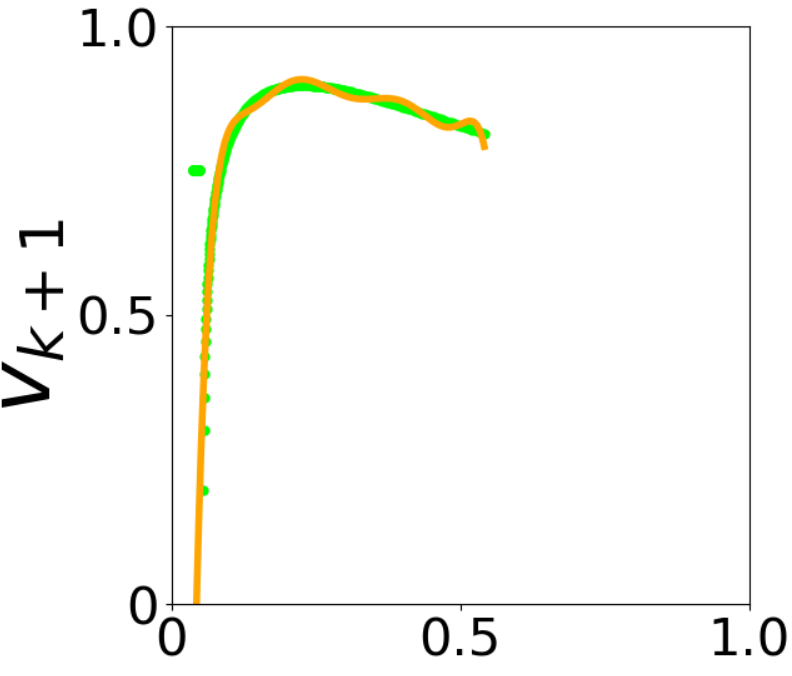}
 \end{subfigure}
 \hfill
 \begin{subfigure}[t]{0.22\textwidth}
     \centering
     \subcaption{}
     \includegraphics[width=\textwidth]{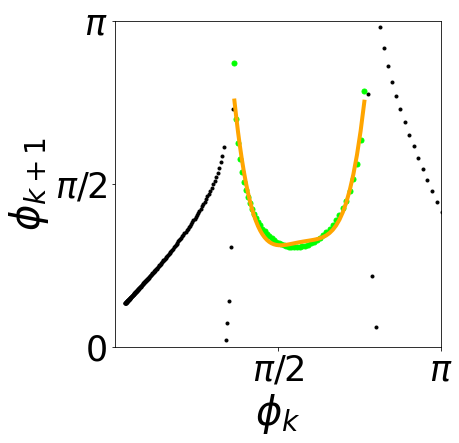}
 \end{subfigure}
 \hfill
 \begin{subfigure}[t]{0.5\textwidth}
     \centering
     \subcaption{}
     \includegraphics[width=\textwidth]{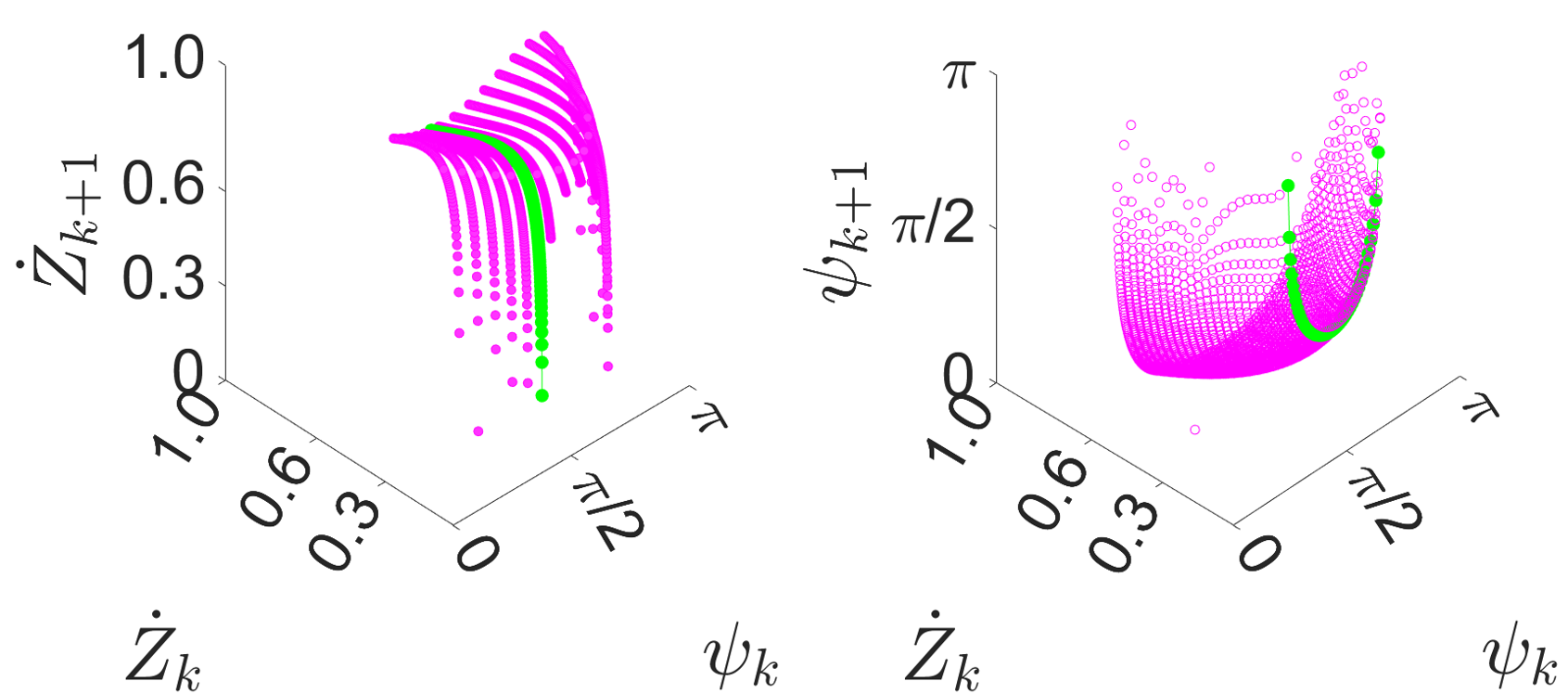}
 \end{subfigure}
 \caption{Illustration of the $P_{\rm BTB}$ surface (magenta surfaces in panels c,f) and its corresponding separable approximation (green and orange curves) for $\mathcal{R}_2$ (panels a, b, c) and $\mathcal{R}_4$ (panels d, e, f), with  $d=0.35$. Generated using the exact map \eqref{exactmap}, the green curves are chosen to represent the variation of the surface for fixed $\psi_k$ or $\dot{Z}_k$. Specifically,  for (c): $\psi_k=0.35$ (left) and $\dot{Z}_k = 0.85$ (right); for (f): $\psi_k = 1.35$ (left) $\dot{Z}_k = 0.12 $ (right). Panels (a)-(b) and (d)-(e) compare the green curves and the orange curves for the approximate separable map \eqref{r2map} in the phase planes. See Appendices \ref{appR2} and \ref{AppR4}  for details.}
 \label{r2approx}
\end{figure}
\item[ v).]Update regions/additional partitions for ${\cal R}_{2.2}$:   As seen from the curve shown in Fig. \ref{r2approx}, which forms the basis of the separable map, the map is not defined on smaller values of $\dot{Z}_k$ in ${\cal R}_{ 2.2}$. This suggests a further partition of ${\cal R}_{2.2}$ into ${\cal R}_{2.3}$ and ${\cal R}_{4.3}$, to capture all values of $\dot{Z}_{k+1}$, as described in Appendices \ref{appR2} and \ref{AppR4}. 
\item [ vi).] No further updates on this optional step.\\

\end{enumerate}
\begin{remark1}\label{R1R2overlap}
Here, we note that the individual curves 
$v_{k+1} = f_2(v_k)$ and $\phi_{k+1} = g_2(\phi_k)$ shown for ${\cal R}_{2.2}$ each overlap with the intervals for  $v_k$ and $\phi_k$ in ${\cal R}_{1.2}$. At first glance, this may seem to cause indeterminacy in the application of the map. In particular, since ${\cal R}_2$ surrounds ${\cal R}_1$, it is possible that one of $v_k$ or $\phi_k$ in ${\cal R}_{2.2}$ can take a value that also appears in the range for ${\cal R}_{1.2}$. However, for $(v_k,\phi_k)$  to be in ${\cal R}_{1.2}$,  both $v_k$ and $\phi_k$ must be in the intervals corresponding to ${\cal R}_{1.2}$. Then $(v_{k+1},\phi_{k+1}) = (f_1,g_1)$ as in \eqref{cubic2dv}-\eqref{cubic2dp},  and not the separable approximation $(f_2(v_k),g_2(\phi_k))$.   
\end{remark1}

We note that while the separable approximation requires some user choice of representative curves, this step is not necessary for determining the composite map ${\cal M}$. We include it here as it aids in visualizing the dynamics in cobweb phase portraits in Section \ref{S5validation}. Furthermore, the separable approximations inspire the auxiliary map applied in Section \ref{AuxiMapSection} to complete the global analysis.

\vskip .2cm

\noindent
\underline{\bf  Iteration 3:} steps iii)-vi)\\
This iteration focuses on ${\cal R}_{2.3}$ and  ${\cal R}_{4.3}$.\\
  \begin{enumerate}
\item[ iii).]  Considering transient dynamics for ${\cal R}_{4.3}$:  For values of small $v_k$ not covered by the approximate map \eqref{r2map} in ${\cal R}_{2.2}$, we consider surfaces as shown in Fig. \ref{r2approx}(f).
\item[ iv).]  Polynomial approximations of  ${\cal R}_{4.3}$:  Similar to the separable maps defined for ${\cal R}_{2.2}$, we use separable single variable approximations $(f_4,g_4)$ for the transient dynamics, given in equation \eqref{R4eqn} and shown in Fig. \ref{r2approx}(d) and \ref{r2approx}(e) .
\item[ v).] No additional partitions are needed.
\item[ vi). ] No further updates needed.\\
\underline{\bf Finalize}\\
vii) We finalize definitions of  the regions ${\cal R}_{n}$, $n=1,2,\ldots,5$ dropping the $.m$ label. The corresponding maps $(f_n,g_n)$ that define the composite map ${\cal M}$ are given in the detailed algorithm in Appendix \ref{algo_code}.
\end{enumerate}

We add some remarks about computational efficiency. In this framework, the main computation  identifies surfaces in regions based on short-time realizations of the impact pair over the state space of initial conditions. These contrast with long-time simulations over the entire state space traditionally used in deriving flow-defined Poincar\'{e} maps for global dynamics of limit-cycle or chaotic systems \cite{guckenheimer1983nonlinear} or for computing basins of attraction \cite{Rounak_basins2022,Kroetz_basins2018}.
 A second feature that contributes to efficiency is the comparison of projections of the surfaces with the diagonals in the phase planes, as in Figs. \ref{2Dpi-2pi} and \ref{2D0-2pi}. As we seek globally attracting dynamics, this division allows us to focus on accurate approximate maps in those regions with attracting dynamics, while relatively cheap approximations are sufficient for transient dynamics.

There are user-defined parameters - polynomial order, surface approximation accuracy, and surface values' proximity to the diagonals in the phase plane, as in determining  ${\cal R}_1$ in \eqref{r1uniond}-\eqref{cubic2dp} and Appendix \ref{appendix-R1}.  These require some iteration to improve the fit of the polynomial approximation to the surface for the region(s) containing the attracting dynamics, but this is not the main computational cost. Furthermore,  defining ${\cal R}_1$ as a union over the parameter $d$ of surfaces from \eqref{r1uniond} also aids in the efficiency of this fit. 

As discussed above, here we have made some additional parameter choices to apply separable maps for convenience in visualization, but they are not a necessary part of the algorithm.  Appendix \ref{algo_code} includes further discussion on values appearing in the algorithm.   
While certain aspects of the computation-based analysis do not rely on finding polynomial approximations for the return maps, we pursue them with the goal of explicit expressions for the global analysis obtained in Section \ref{S6.3GD}.

\section{Validation of the Composite Map}
\label{S5validation}

 In this section, the composite map ${\cal M}$ is validated using three distinct types of solutions, showing that it can reproduce the dynamics of different types of solutions.  The first type of solution is the fixed point of ${\cal M}$, which we call Case FP, corresponding to the 1:1/${\cal T}$ solution of the full system \eqref{cyl_dim}-\eqref{simpact_dim}. The second type is the period doubled case, i.e., the period-2 cycle of ${\cal M}$, called Case PD, corresponding to the 1:1/$2T$ behavior in the full system. Lastly, the chaotic dynamics of ${\cal M}$,  called Case CD, corresponds to the chaotic 1:1/$C$ behavior in the full system.  These different dynamics can be observed from the bifurcation diagrams in Figs. \ref{bifexact}, \ref{bifapprox} for $d = 0.35$, $d=0.30$, and $d=0.26$, respectively.

\begin{figure}[htbp]
 \centering
 \begin{subfigure}[b]{0.45\textwidth}  
     \centering
     \subcaption{}     \includegraphics[width=\textwidth]{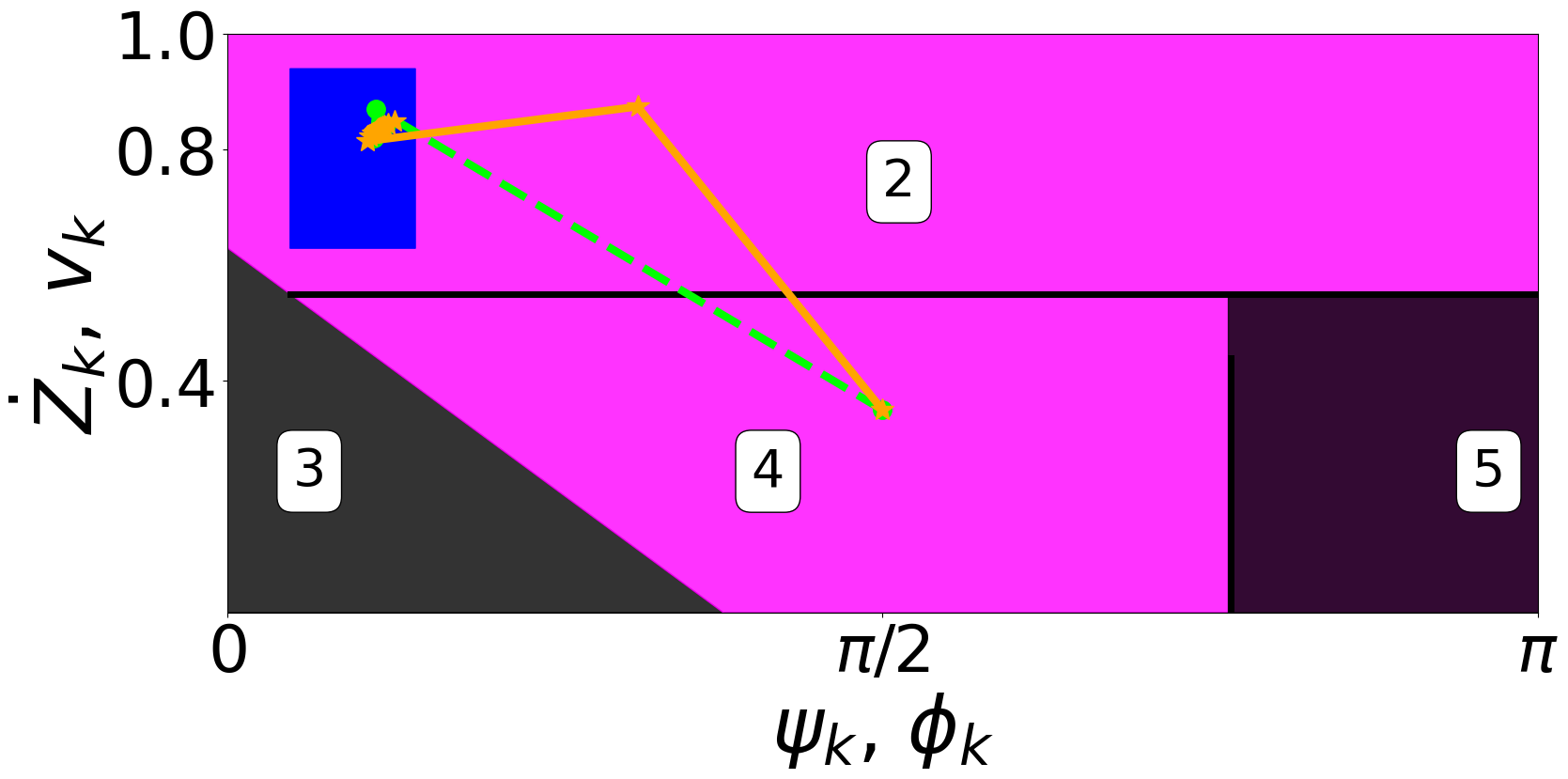}
 \end{subfigure}
 \hfill
 \begin{subfigure}[b]{0.45\textwidth}  
     \centering
     \subcaption{}
     \includegraphics[width=\textwidth]{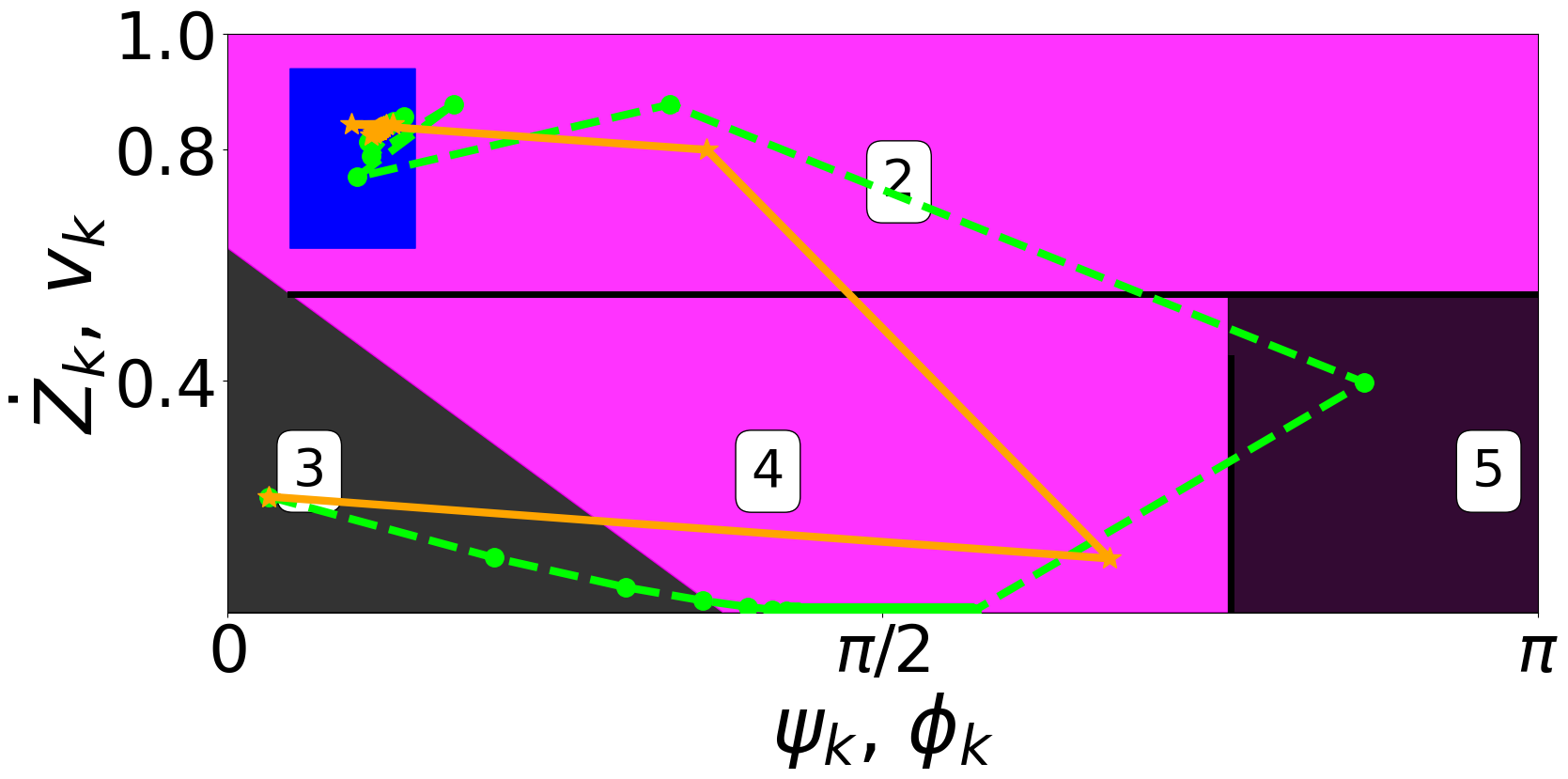}
 \end{subfigure}
 \vfill
 \begin{subfigure}[b]{0.45\textwidth}  
     \centering
     \subcaption{}
     \includegraphics[width=\textwidth]{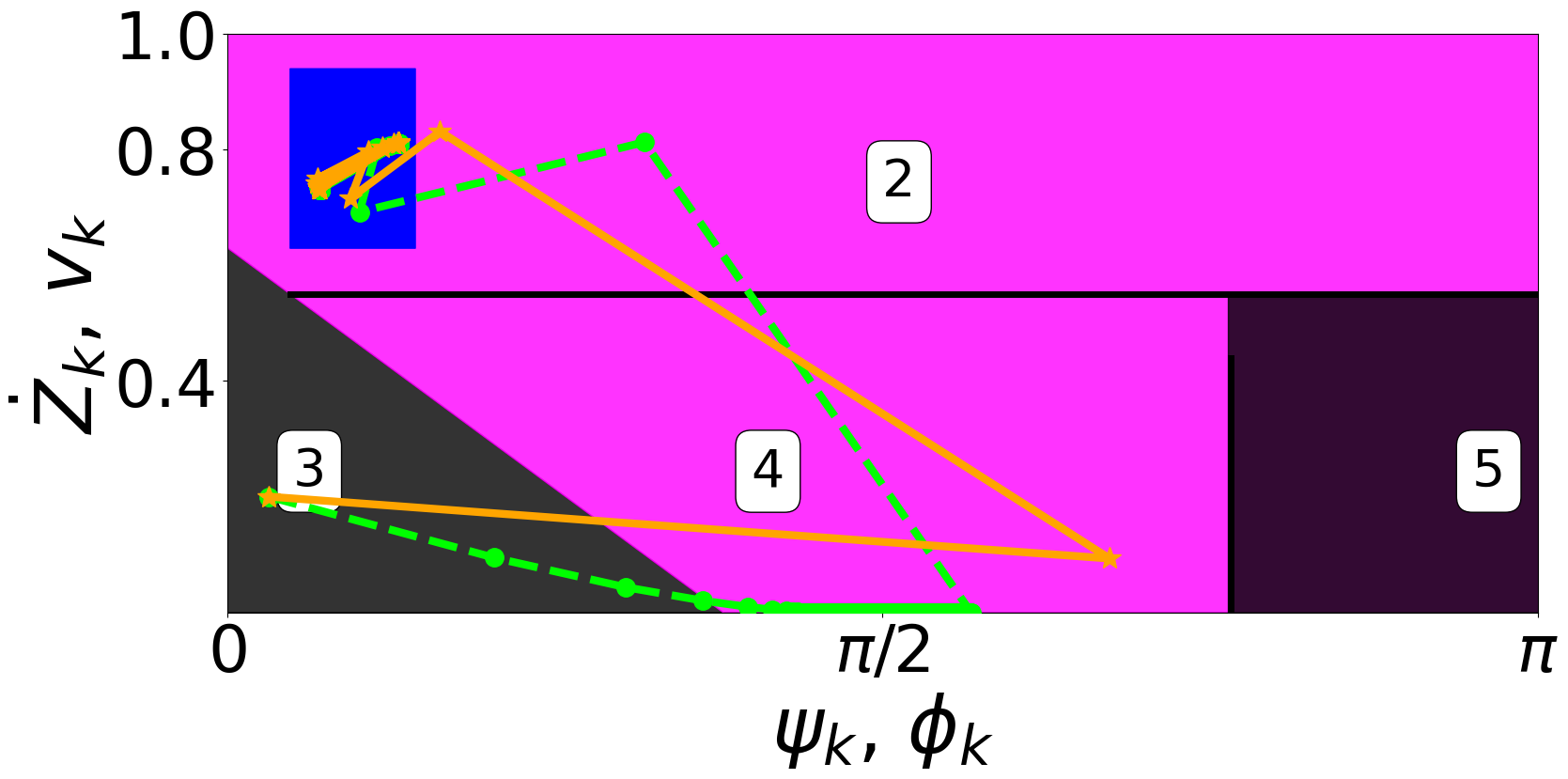}
 \end{subfigure}
 \hfill
 \begin{subfigure}[b]{0.45\textwidth}  
     \centering
     \subcaption{}
     \includegraphics[width=\textwidth]{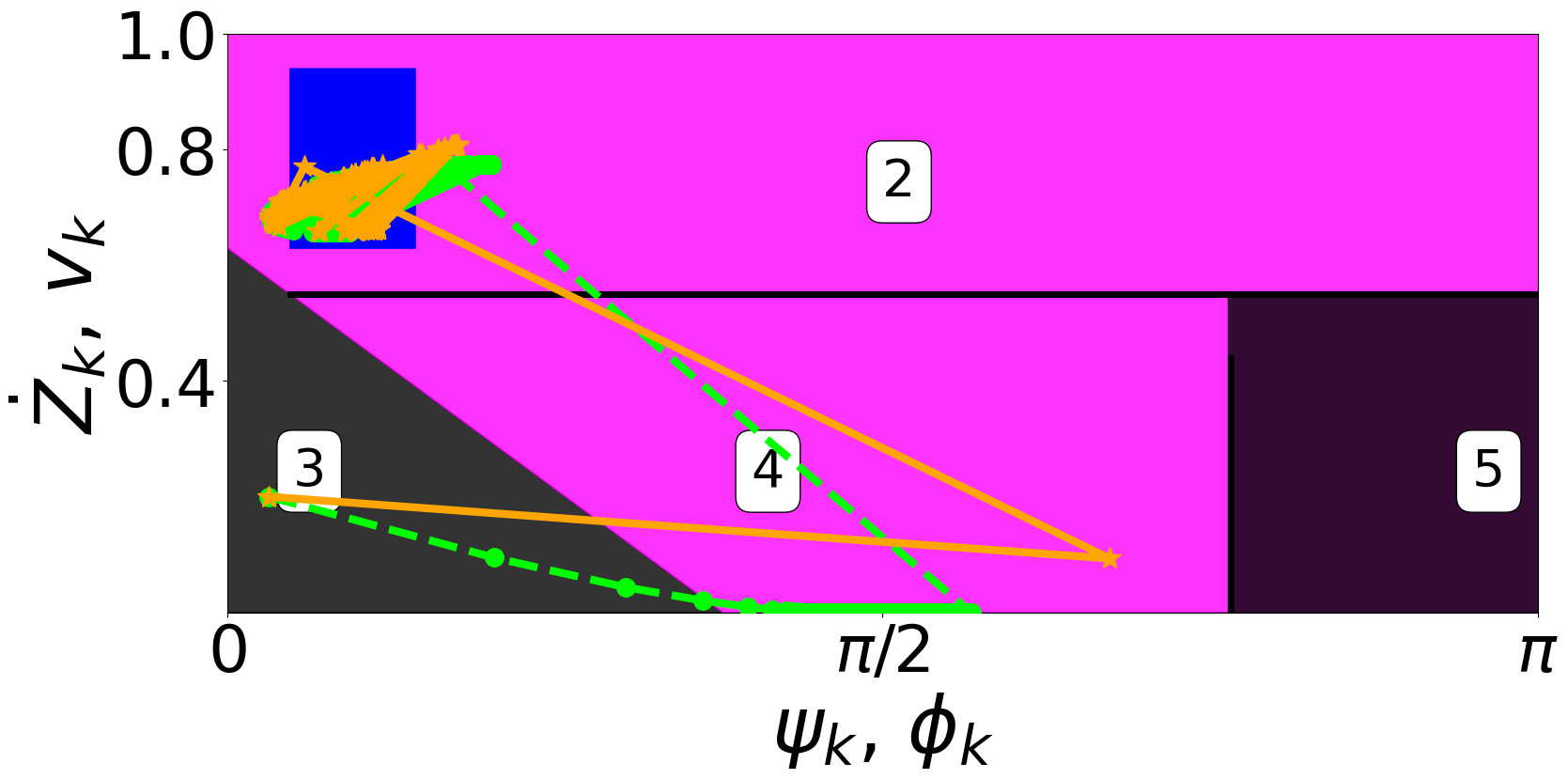}
 \end{subfigure}
 \caption{ 
 Comparison of trajectories in state space from the exact map \eqref{exactmap} (orange) and the composite map ${\cal M}$ \eqref{compositemapM} (green), superimposed on regions  $\mathcal{R}_n$ used in the definition of  ${\cal M}$ as specified in 
 Appendix \ref{algo_code}.  (a) and (b) correspond to Case FP,  also shown in cobweb phase portraits in Fig.~\ref{glosmallic}(a),(b); (c) corresponds to Case PD, also shown in Fig.~\ref{glosmallic}(c),(d); (d) corresponds to Case CD, also shown in  Fig.~\ref{glosmallic}(e),(f). 
 Parameters and initial conditions: (a) $d=0.35$,  $\phi_0=\pi/2, v_0=0.35$; (b) $d=0.35$, $\phi_0=0.1, v_0=0.2$; (c) $d=0.30$, $\phi_0=0.1, v_0=0.2$; (d) $d=0.26$, $\phi_0=0.1, v_0=0.2$. Here, we show representative results for initial conditions in the transient regions ${\cal R}_3$ and ${\cal R}_4$.}
 \label{test}
\end{figure}

Figure~\ref{test} shows the implementation of the composite map $\mathcal{M}$ (dashed green line), with the corresponding pseudocode given in  Appendix \ref{algo_code}.  Initial condition pairs $(v_k, \phi_k)$ are selected from transient regions ${\cal R}_3$ and ${\cal R}_4$ to demonstrate that ${\cal M}$ can reliably predict the long-term system behavior, reaching a potentially attracting region after traveling through other transient regions ${\cal R}_n$.  Similar results were obtained for other randomly selected initial pairs (not shown here).   Trajectories for ${\cal M}$ are plotted together with the trajectories generated with the exact map \eqref{exactmap} (solid orange line). Panels (a) and (b) correspond to Case FP. Panels (c) and (d) correspond to Case PD and Case CD, respectively. In all cases, both ${\cal M}$ and the exact map \eqref{exactmap} trajectories follow each other to reach the same attracting dynamics. Of course, the transient dynamics are not reproduced exactly, e.g., given the less accurate separable approximations used in ${\cal M}$ to facilitate visualization of the maps.

\begin{figure}[htbp]
\centering
\begin{subfigure}{0.78\textwidth}  
    \centering 
\includegraphics[width=\textwidth]{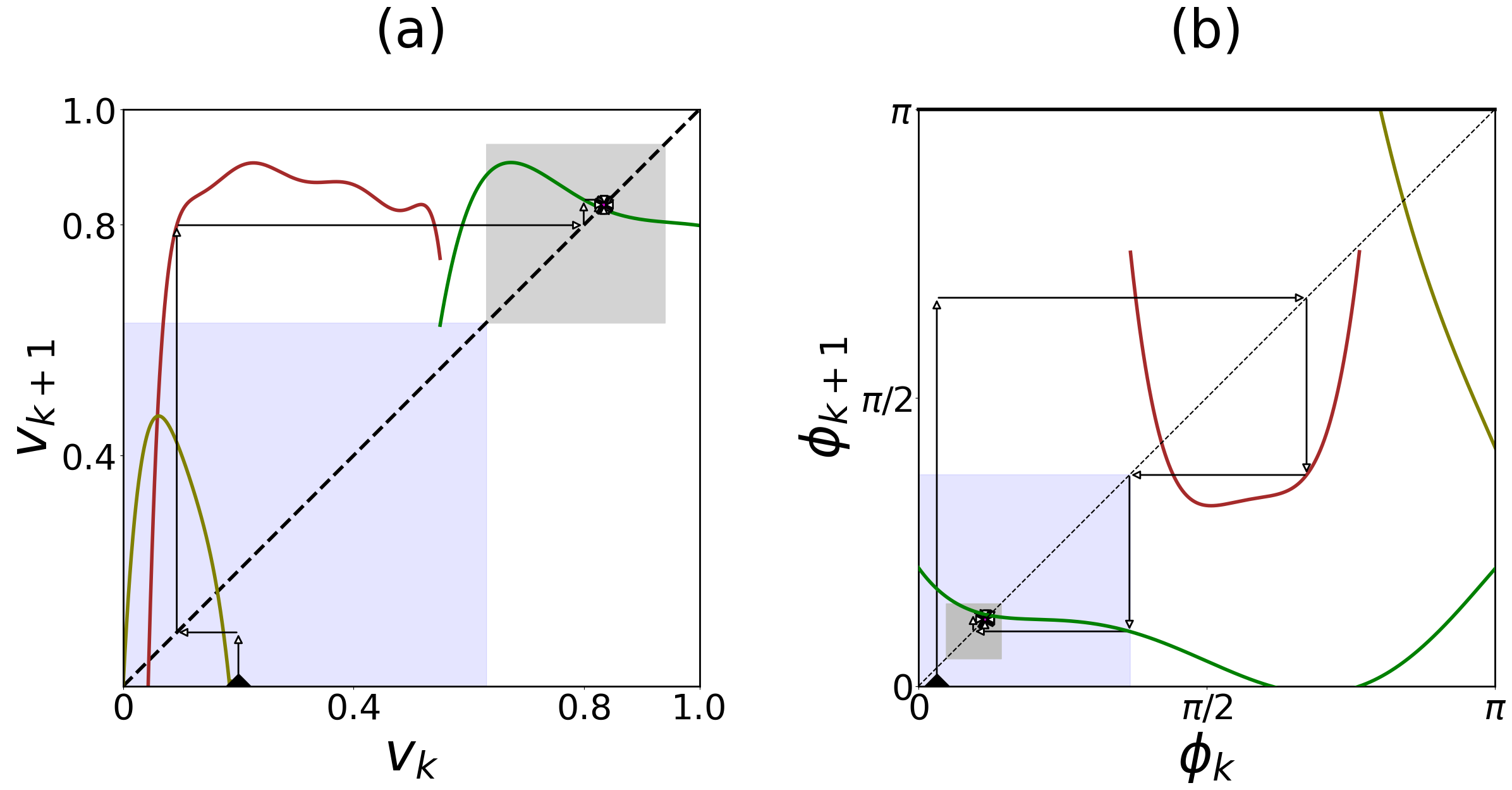}
\end{subfigure}
 \vfill
 \begin{subfigure}{0.78\textwidth} 
 \centering
\includegraphics[width=\textwidth]{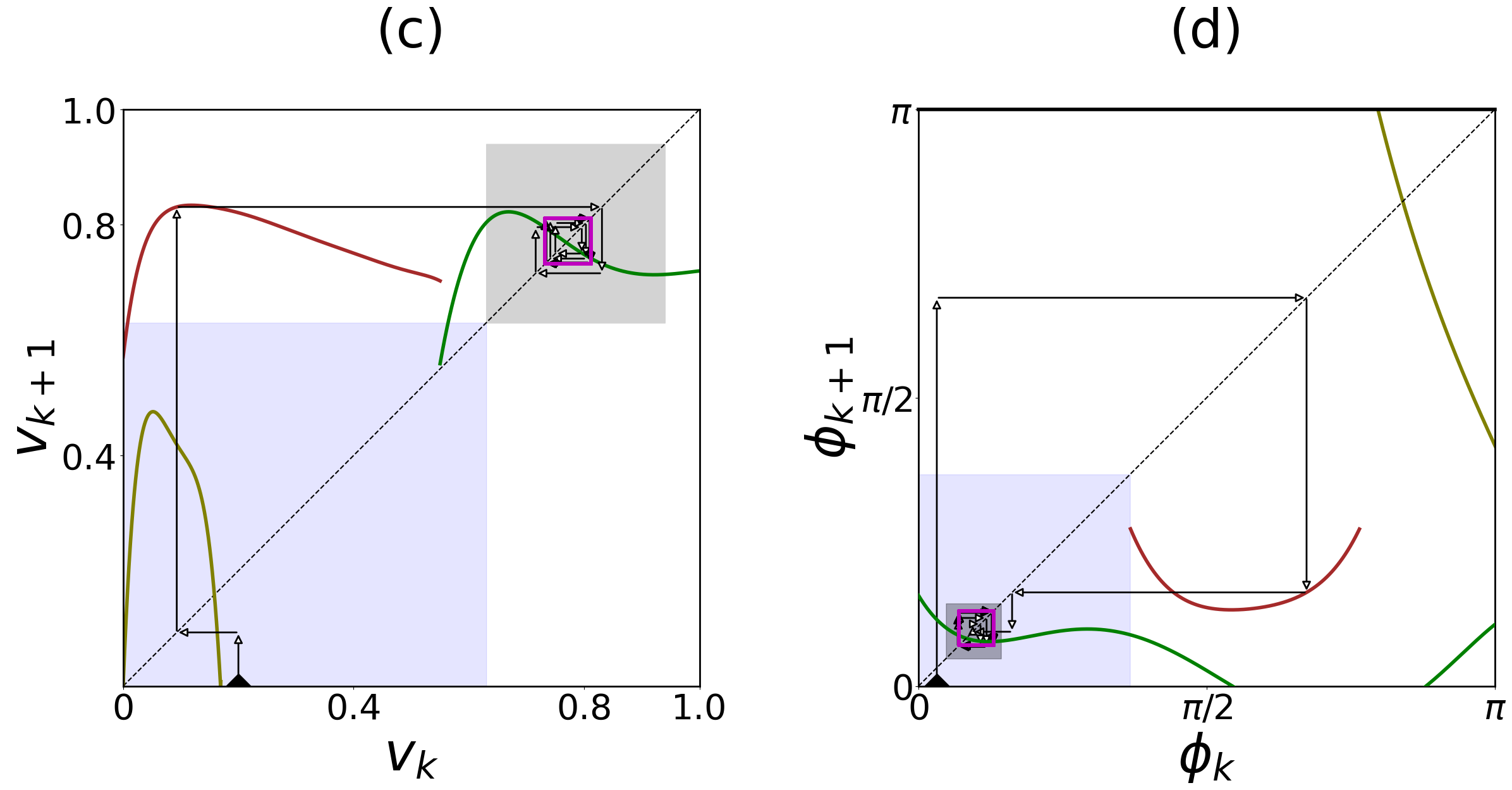}
 \end{subfigure}
 \vfill
 \begin{subfigure}{0.78\textwidth}
 \centering
\includegraphics[width=\textwidth]{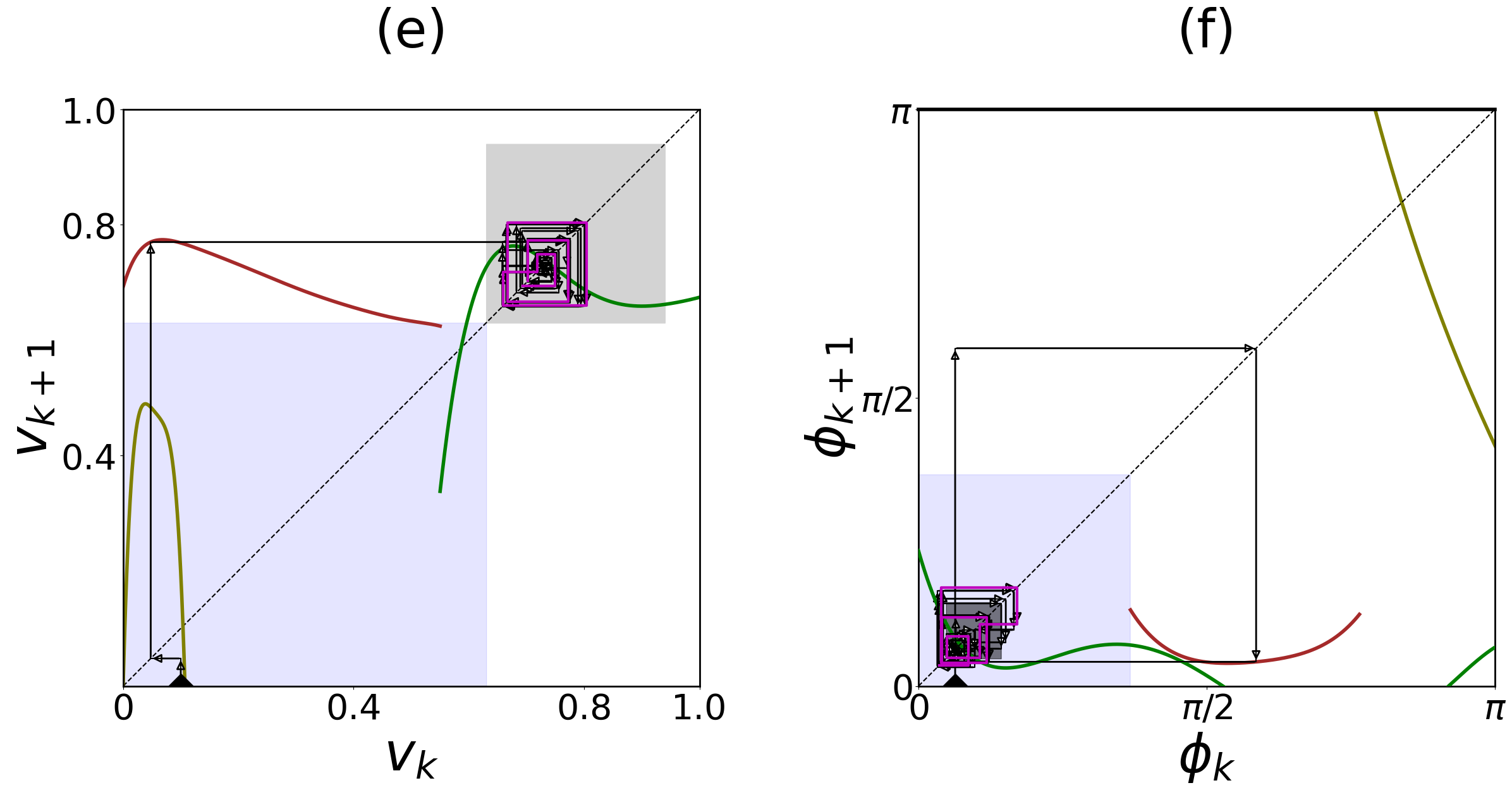}
 \end{subfigure}
\caption{ Application of  ${\cal M}$ \eqref{compositemapM} projected  on the $v_{k}$ and $\phi_k$ phase planes, with step  navigation for $(f_n,g_n)$ discussed in the text. 
Curves show (separable) maps for Regions $\mathcal{R}_2$ (green), $\mathcal{R}_4$ (red), and $\mathcal{R}_5$ (olive). Shaded regions are for approximate 2D maps for $\mathcal{R}_1$ (gray) and $\mathcal{R}_3$ (light blue), which can not be drawn in these projections. Black dashed lines show the respective diagonals.
Parameters: Case FP (a),(b): $d=0.35, v_0=0.2, \phi_0=0.1$; Case PD (c),(d):  $d=0.30, v_0=0.2, \phi_0=0.1$; Case CD (e),(f):  $d=0.26, v_0=0.1, \phi_0=0.2$. Supplementary Video 2 provides a step-by-step demonstration of using the overlapped curves in cobwebbing.
}
\label{glosmallic}
\end{figure}

Complementary to the validation of  ${\cal M}$ in Fig. \ref{test},  Fig.~\ref{glosmallic} demonstrates the attracting behavior in the projected $v_{k}-v_{k+1}$ and $\phi_{k}-\phi_{k+1}$ phase planes with initial conditions for small $v_k$ and $\phi_k$ ($v_0 = 0.2, \phi_0 = 0.1$). Repeated application of the composite map is demonstrated via cobweb phase portraits, indicating the steps toward the attracting behavior. The dynamic behavior is shown for the three types of solutions listed above.  In both Case FP and PD, the trajectories limit to values within $\mathcal{R}_1$ while in Case CD, the long-term trajectory takes values in $\mathcal{R}_1$ and $\mathcal{R}_2$. All of these are consistent with the bifurcation structure shown in Fig. \ref{bifapprox} (and in Fig. \ref{bifexact}).

For the cobweb analysis using the maps $(f_n,g_n)$ in the $v_{k}-v_{k+1}$ and $\phi_{k}-\phi_{k+1}$ phase planes shown in Fig.~\ref{glosmallic},  it is possible to visualize the curves for the maps in $\mathcal{R}_2$, $\mathcal{R}_4$, and $\mathcal{R}_5$, as we use separable (1D) approximations in those regions. For $\mathcal{R}_1$ and $\mathcal{R}_3$ we can not show a single curve in this projection, given the 2D polynomial map used in \eqref{cubic2dv}-\eqref{cubic2dp} and \eqref{Apr3eqn}, respectively. Instead, shaded regions show the range of $v_k$ and $\phi_k$ in $\mathcal{R}_1$  (gray) and $\mathcal{R}_3$ (light blue). Then, the cobweb steps in these regions follow the (surface) maps \eqref{cubic2dv}-\eqref{cubic2dp} and \eqref{Apr3eqn} for ${\cal R}_1$ and ${\cal R}_3$, respectively, for $(v_k,\phi_k)$ in these regions, even though specific curves are not shown.
Given the width of these shaded regions, it is possible to give a maximum and minimum for $v_{k+1}$ and $\phi_{k+1}$, which also motivates the auxiliary map defined and applied in Section \ref{AuxiMapSection} for $\mathcal{R}_1$.

We provide some navigation in order to trace the cobweb behavior for ${\mathcal M}$ as shown in Fig.~\ref{glosmallic}. A more detailed step-by-step navigation is provided in Appendix \ref{NaviFig13} and Supplementary Video 2. Since the panels show projections of the higher dimensional maps $(f_j,g_j)$ in the phase planes, there is an overlap in these projections, and thus, one must take care to use the correct map for $v_k$ and $\phi_k$ in the overlap region. Specifically, for each cobweb step, $(v_{k+1},\phi_{k+1})$ takes a value according to the map for the region that is common for both $(v_{k},\phi_{k})$. In all cases shown, the initial condition $(v_k,\phi_k)$ for $k=0$ takes small values in ${\cal R}_3$. 
We observe that ${\cal R}_3$, ${\cal R}_4$, and ${\cal R}_5$ overlap in the $v_k-v_{k+1}$ phase plane for these smaller values of $v_k$, while in the $\phi_k-\phi_{k+1}$ phase plane the curve for ${\cal R}_2$ and region $\mathcal{R}_3$ overlap for smaller $\phi_k$. Since ${\cal R}_3$ is the only region in common for $v_k$ and $\phi_k$ for these small values, we conclude that $(v_k,\phi_k)\in {\cal R}_3$,  and the first step follows $(v_{k+1}, \phi_{k+1}) = (f_3(v_k,\phi_k),g_3(v_k,\phi_k))$ in \eqref{Apr3eqn},  as shown in Fig.~\ref{glosmallic}. The result gives the new $(v_k, \phi_k)$, for the next step with $k=1$, for which $v_k$ remains small while $\phi_k$ has increased (before reaching the attracting dynamics observed for later steps in ${\cal R}_1$). 
Again ${\cal R}_3$, ${\cal R}_4$, and ${\cal R}_5$ overlap in the $v_k - v_{k+1}$ phase plane for these smaller values of $v_k$, while in the $\phi_k - \phi_{k+1}$ plane, $\phi_k$ takes a value corresponding to the range for ${\cal R}_4$ only. Then
$(v_{k+1},\phi_{k+1})$ follows the map $(f_4,g_4)$ for ${\cal R}_4$. Note that the approximate maps for ${\cal R}_3$ and ${\cal R}_5$ are not applied for $v_k$, even though $v_k$ takes values in their range, since  $\phi_k$ is not in either ${\cal R}_3$ or ${\cal R}_5$.
Eventually, for a larger $k>1$,  $v_k$ has increased to a range with an overlap between ${\cal R}_2$ and ${\cal R}_1$, while $\phi_k$ has decreased back to the region with overlap between ${\cal R}_2$, ${\cal R}_1$ and $\mathcal{R}_3$.  Then, the cobweb steps are governed by $(f_2,g_2)$  for   $(v_k,\phi_k)\in {\cal R}_2$, and by   $(f_1,g_1)$ in \eqref{cubic2dv}-\eqref{cubic2dp} for  $(v_k,\phi_k)\in\mathcal{R}_1$, as already discussed in Remark \ref{R1R2overlap}  about the overlap between the green curves and the grey shaded ${\cal R}_1$ region. From there, the dynamics are dictated by the attracting dynamics of ${\cal R}_1$ for panels (a),(b) and (c),(d) corresponding to Cases FP and PD, respectively. In panels (e)  and (f), the attracting chaotic dynamics for Case CD alternate between ${\cal R}_1$ and ${\cal R}_2$, as described in Remark \ref{R1R2overlap}.

\section{The Auxiliary Map Method for Global Dynamics} \label{AuxiMapSection}
It is worth noting that a computationally realized implicit composite map could have been employed up to this point, bypassing the need for polynomial approximations of the surfaces in Fig.~\ref{3d}. While such a map could still offer insights into the system's global dynamics, it would not allow the explicit computer-assisted analysis of the system's attracting domain. This limitation underscores the value of our explicit composite map, which is an analytical tool for deriving tight bounds on the size of the system's non-trivial attractors through analyzing the auxiliary maps, as demonstrated in Section \ref{S6.3GD}. 

The trajectories above indicate visually that Regions $\mathcal{R}_1$ and $\mathcal{R}_2$ contain an {attracting domain} that attracts all non-trivial trajectories in $\mathcal{R}_1$ and $\mathcal{R}_2$ for the considered range of parameter $d$. In particular, the magenta orbits in Fig. \ref{glosmallic} highlight the last 10\% of the cobweb trajectories, and the stable orbits correspond to the solution given by the composite map $\mathcal{M}$ \eqref{exactmap}. In Case FP, the solution is a fixed point, shown in panel (a)(b), and is contained in $\mathcal{R}_1$. In Case PD, the solution has period two and is also contained in $\mathcal{R}_1$, as shown in panel (c)(d). In Case CD, the solution is chaotic but is also contained within $\mathcal{R}_1$ and $\mathcal{R}_2$, shown in panel (e)(f). Therefore, the stable orbits shown in Fig. \ref{glosmallic} indicate the existence of an attracting domain. However, determining the bounds for the attracting domain as a subset of regions $\mathcal{R}_1$ and $\mathcal{R}_2$ directly from the 2D composite map is out of reach, especially for multi-period and chaotic dynamics.
In Fig.~\ref{glosmallic}, iterations of the closed-form composite map visualize the system's long-term behavior, with explicit curves shown only for regions ${\cal R}_2$, ${\cal R}_4$, and ${\cal R}_5$ when projected onto the $\dot{Z}_{k+1}-\dot{Z}_{k}$ and  $\phi_{k+1}-\phi_k$ planes. In contrast,  for  ${\cal R}_1$ and ${\cal R}_3$ the maps cannot be visualized under this projection, suggesting that an alternate approach is needed to capture global attraction using these cobweb phase portraits. The difference between the regions follows from the separable form of the maps in ${\cal R}_2$, ${\cal R}_4$, and ${\cal R}_5$, in contrast to the 2D maps of ${\cal R}_1$ and ${\cal R}_3$. This observation inspires the design of an auxiliary map, in which we dissect each 2D map into a pair of 1D maps based on the lower and upper bounds of the 2D map domain. This definition can then take advantage of the separable form and lead to bounds on the composite map's {attracting domain}. 
 
\subsection{Constructing the Auxiliary Maps}\label{S6.1auxiliary}

The auxiliary map is constructed using the bounds on the approximate maps $(f_n,g_n)$ for each Region $\mathcal{R}_n$, where  $(f_n,g_n)$  depends on both variables $v_k$ and $\phi_k$. In our case, these regions are $R_1$ and $R_3$.  We define the auxiliary maps in terms of the maxima and minima of $(f_n,g_n)$,  yielding the form for {the upper bound:} $\xi_{U} (v_k): v_k \to v_{k+1}$ and $\eta_{U} (\phi_k): \phi_k\to \phi_{k+1}$, and similarly for the lower bound. This decouples the two 2-D equations into two separable 1-D equations for each ${\mathcal R}_n$. The advantage of this formulation is its ability to track the dynamics of velocity $v_k$ and the phase $\phi_k$ separately, thus facilitating a 1D cobweb phase portrait for each.  At the same time, it captures the worst-case scenario and provides conservative bounds on the maximum and minimum range of $(f_n,g_n)$ at each iterate. Furthermore, we show that a repeated application of this auxiliary definition hones in on the attracting solutions or regions of the full map. While here we give the construction in terms of general region number $n$, we emphasize that below it is applied for $\mathcal{R}_1$ only, as we focus on the attracting behavior.

The construction of the auxiliary map begins with the bounds for $v_k$ and $\phi_k$ for a given ${\cal R}_n$: $v_k\in [v_0^{\rm min}, v_0^{\rm max}]$ and $\phi_k\in [\phi_0^{\rm min}, \phi_0^{\rm max}]$.  Then two curves $\xi_U(v_k)$  and $\xi_L(v_k)$ are determined for $v_{k+1}$  in terms of the $\max$ and $\min$ of $f_n$ over the range of possible $\phi_k$ values,
\begin{eqnarray}\label{auxv}
\xi_n^{(N)} =\left\{
\begin{array}{l} v_{k+1} = \xi_U^{(N)}(v_k), \;\;\;{\rm where}\;\; \xi_{U}^{(N)}  
:= \max \limits_{\phi \in {\cal A}_n^{(N)}} \{ f_n(v_k, \phi)\},\\
 v_{k+1} = \xi_{L}^{(N)}(v_k), \;\;\;{\rm where}\;\; \xi_{L}^{(N)}:=
  \min \limits_{\phi \in {\cal A}_n^{(N)}} \{ f_n(v_k, \phi)\}.
\end{array}\right. 
\end{eqnarray}
The auxiliary maps, $\xi_n^{(N)}$ and 
$\eta_n^{(N)}$ defined similarly below, alternate between the two curves in order to provide a sequence of bounds on the maximum and minimum range of $(f_n,g_n)$. Here, we use a generic initial $v_k$, with refinements discussed below and in  Section \ref{S6.2}. The superscript $N$ gives the index of updates of the auxiliary map after the first and subsequent applications, particularly valuable when the auxiliary map is contracting, as demonstrated below for the specific cases considered in Section \ref{S6.2}. Likewise, the auxiliary map $\eta_n^{(N)}$ is given in terms of  two  maps $\eta_U$, $\eta_L$ that bound  $\phi_{k+1}$ for $v_k\in [v_0^{\rm min}, v_0^{\rm max}]$:
\begin{eqnarray}\label{auxphi}
\eta_n^{(N)} =\left\{
\begin{array}{l}\phi_{k+1} = \eta_U^{(N)}(\phi_k), \;\;\;{\rm where}\;\; \eta_{\rm max}^{(N)} := \max_{v\in {\cal A}_n^{(N)}} \{ g_n( v, \phi_k)\}, \\
 \phi_{k+1} = \eta_L^{(N)}(\phi_k), \;\;\;{\rm where}\;\; \eta_{\rm min}^{(N)} := \min_{v\in {\cal A}_n^{(N)}} \{ g_n(v, \phi_k)\} .
 \end{array} \right. \label{auxphi}
\end{eqnarray}
To track the (possible) contraction of the region for each update, we define  ${\cal A}_n^{(N)}$ in \eqref{An}-\eqref{Bn} below.  There  ${\cal A}_n^{(N)} = {\cal R}_n$  for all $N$ if the region does not contract, while ${\cal A}_n^{(1)} =  {\cal R}_n$ and ${\cal A}_n^{(N)} \subseteq {\cal R}_n$ for $N>1$ for a contracting region,  updated as the auxiliary map is updated. For the system studied here, it is only for $n=1$ that ${\cal A}_n^{(N)}$ contracts.

We then write the full auxiliary map, replacing  $\mathcal{M}$\, \eqref {compositemapM}  with  $\mathcal{M}_{\cal A}^{(N)}$, which is composed of a combination of maps $(f_n,g_n)$ and $(\xi_n^{(N)},\eta_n^{(N)})$, with $ v_k, \phi_{k}$ corresponding to impact velocities on $ \partial B$ as in \eqref{compositemapM}. For our system it is only ${\cal A}_1^{(N)}$ that contracts as $N$ increases, so we define the full auxiliary map as
\begin{eqnarray}
 (v_{k+1}, \phi_{k+1}) & = & {\cal M}^{(N)}_{\mathcal A} (v_k,\phi_k) ,
 \nonumber\\
{\cal M}^{(N)}_{\mathcal A} (v_k,\phi_k)  &\equiv & \left\{
\begin{array}{ll}
(\xi_1^{(N)}(v_{k}),\eta_1^{(N)}(\phi_k))  & \mbox{ for } (v_k,\phi_k) \in {\cal A}_1^{(N)}, \\
(\xi_3^{(N)}(v_{k}),\eta_3^{(N)}(\phi_k))  & \mbox{ for } (v_k,\phi_k) \in {\cal R}_3, \\
(f_n(v_{k},\phi_k),g_n(v_{k},\phi_k))  &  \mbox{ for } (v_k,\phi_k) \in {\cal R}_n, n = 2,4,5 .
\end{array} \right.  \label{auxM}
\end{eqnarray}

 We define region $\mathcal {A}_1^{(N)}\subseteq \mathcal{R}_1 $ to allow a change in its size  over the  $N$ updates of the auxiliary construction,
\begin{eqnarray}
{\mathcal A}_1^{(N)} & = {\mathcal A}_{1v}^{(N)}\times{\mathcal A}_{1\phi}^{(N)} =  \left\{ 
\begin{array}{ll}
  \mathcal{R}_1 &  \mbox{ for } N=1,\\
 \mathcal{B}_1^{(N)} &  \mbox{ otherwise, }
 \end{array}
    \right. \label{An}\\
 \mathcal{B}_1^{(N)}  & =   {\cal B}^{(N)}_{1v} \times {\cal B}^{(N)}_{1\phi} \equiv
 [v_{\ell}^{\rm min},v_{\ell}^{\rm max}] \times
[\phi_{\ell}^{\rm min},\phi_{\ell}^{\rm max}]  \label{Bn}\\
 &  \quad 
\mbox{ for }  (v_{\ell},\phi_{\ell}) = 
\left({\cal M}_{\cal A}^{(N-1)}\right)^\ell (v_0,\phi_0) , \ \ell \gg 1. \nonumber
\end{eqnarray}
 Stated in words,  \eqref{An}-\eqref{Bn}  simply indicate that for the $N^{\rm th}$ ($N>1$) update of $(\xi_1^{(N)}(v_{k}),\eta_1^{(N)}(\phi_k))$,  the region ${\mathcal A}_1^{(N)}$ is updated to  the limiting range of $(v_k,\phi_k)$  obtained from a large number of iterations of $(\xi_1^{(N-1)}(v_{k}),\eta_1^{(N-1)}(\phi_k))$ using \eqref{auxv}-\eqref{auxphi}. Given the separable form of \eqref{auxv}-\eqref{auxphi}, both ${\cal A}$ and ${\cal B}$ are defined in terms of the ranges of $v$ and $\phi$.

The iteration of \eqref{auxv}-\eqref{auxphi} is particularly valuable for region(s) in which the dynamics are contracting since these iterations identify a relaxation within the extremes imposed by the defined maxima and minima, leading to an update of the region  ${\cal A}_1^{(N+1)}$ and ${\cal M}^{(N+1)}_{\mathcal A}$ as in \eqref{An}-\eqref{Bn}.  Here, we have proposed  \eqref{auxv}-\eqref{auxphi} starting from generic values of $(v_k,\phi_k)$.   In  Section \ref{S6.2}, we refine the iterative cobwebbing-type application of the auxiliary maps based on a choice of $(v_k,\phi_k)$ that ensures improvements within the worst-case scenario.
Then, repeated updates for increasing $N$ give conservative bounds on the limiting size of the attracting domain. 

\subsection{Application of the auxiliary map $\mathcal{M}_{\cal A}^{(N)}$}\label{S6.2}

In Section ~\ref{S5validation}, the application of $\mathcal{M}$ via cobweb phase portraits indicates that the {attracting dynamics} are concentrated in $\mathcal{R}_1$ for the larger values of $d$ considered in this study. Specifically,  in Fig.~\ref{glosmallic}, we see attracting solutions contained in $\mathcal{R}_1$ in Case FP and PD, while the trajectories oscillate between $\mathcal{R}_1$ and $\mathcal{R}_2$ in Case CD. 

While we could construct an auxiliary map in the setting where the dynamics revisit regions with transient dynamics (e.g., $\mathcal{R}_2$), this would require a different construction to be useful in demonstrating global stability; instead, the {attracting dynamics} suggest a more efficient approach. From Fig.~\ref{glosmallic},  the {attracting domain} covers values in $\mathcal{R}_1$ for Cases FP and PD, and in a region just outside of $\mathcal{R}_1$ for Case CD.  This suggests constructing the auxiliary map on a slightly expanded region $\mathcal{R}_1^+\supseteq \mathcal{R}_1$, noting that this does not reduce the accuracy of the approximation as it uses the more accurate $2D$ approximation over a larger region, reducing the region over which the separable approximation $(f_2,g_2)$ is used. Then we can expand the size of Region $\mathcal{R}_1$ to $\mathcal{R}_1^+$ sufficiently so that the long-term dynamics remain in $\mathcal{R}_1^+$ and $\mathcal{R}_1^+ \supseteq \mathcal{R}_1$, and here we consider the auxiliary map for $\mathcal{R}_1^+$ only. 

The following are the ranges of the initial region $\mathcal{A}_1^{(1)}={\cal R}_1^+$ for the three cases, the fixed point (FP) case, the period-doubling (PD) case, and the chaotic dynamics (CD) case of the composite map $\mathcal{M}$:
\begin{eqnarray}
\text{{\bf Case FP:} \ \ }    {\cal R}_1^+ := &\  & \{ (v_k, \phi_k):\; v_k\in [0.7,1] \text{ and }  \phi_k\in [0.2, \pi/3]\} \label{R1+FP}\\
 \text{\bf Case PD:\ \ } {\cal R}_1^+ := &\  & \{ (v_k, \phi_k):\; v_k\in [0.65,1]
 \text{ and }  \phi_k\in [0.13, \pi/3]\} \label{R1+PD} \\
  \text{\bf Case CD:\ \ }   {\cal R}_1^+ := &\ & \{ (v_k, \phi_k):\; v_k\in [0.64,1]  \text{ and }   \phi_k\in [0.08, \pi/3]\}. \label{R1+CD}
\end{eqnarray}

Here, ${\cal R}_1^+$ is typically an over-estimate of the attracting domain, given that it is based on the approximation obtained by comparing the projection of the exact maps with diagonals in the phase planes shown in Section \ref{S3return-map}.

\begin{figure}[htbp]
    \centering
\includegraphics[width=0.65\textwidth]{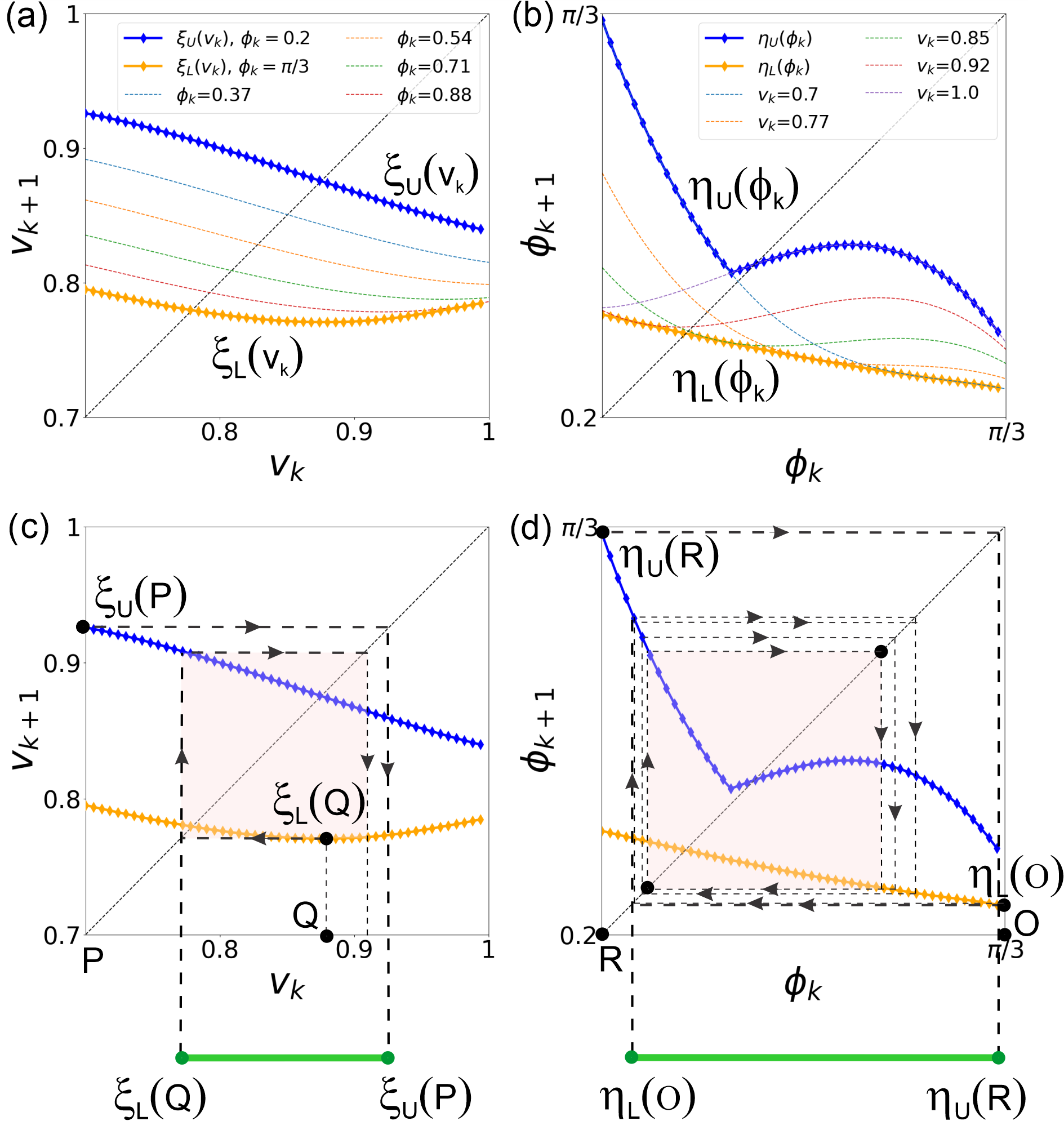}
\caption{(a)-(b): Visualization of 
the auxiliary maps $\xi_{U/L}^{(1)}$ and $\eta_{U/L}^{(1)}$ from \eqref{auxv} and \eqref{auxphi} 
for $\mathcal{R}_1^+$ and $d=0.35$. Here, we dropped the superscript on $\xi_{U/L}$ and $\eta_{U/L}$ for clarity of visualization.
  In (a)-(b) $f_1$ and $g_1$ are graphed for fixed $\phi_k$ and $v_k$ values, respectively, confined between their respective lower, $\xi_L(v_k)$ and $\eta_L(\phi_k)$ (orange diamonds) and upper, $\xi_U(v_k)$ and $\eta_U(\phi_k)$ (blue diamonds) bounds.
(c): The points $P$ and $Q$ are the location of the respective maximum and minimum of the maps $\xi_U$ (blue) and $\xi_L$ (orange) over the entire interval $I_{0v}^{(1)}$ for $v \in {\cal R}_1^+$. The images $\xi_U(P)$  and $\xi_L(Q)$ define the reduced interval $I_{1v}^{(1)}$ (green). The black-dotted line indicates the successive application of \eqref{WCScobweb-v}, each time applying  $\xi_U$  and $\xi_L$  to the points corresponding to their max and min, respectively, over the interval $I_{jv}^{(1)}$,
$j=1,2,\ldots$. 
In panel (d) iterations are similar to those in (c), shown for $\phi_k$, $\eta_{U/L}$  and the auxiliary map \eqref{WCScobweb-phi}. The first iterate applies \eqref{WCScobweb-phi} to the points $R$ and $O$, corresponding to the max and min of $\eta_{U}$ and $\eta_L$ on the full interval $I_{0\phi}^{(1)}$ for $\phi\in {\cal R}_1^+$, whose image yields the reduced interval $I_{1\phi}^{(1)}$ (green). The black-dotted line indicates successive application of \eqref{WCScobweb-phi} to the points corresponding to the max and min of $\eta_{U}$ and $\eta_{L}$, respectively, over the interval $I_{j\phi}^{(1)}$,
$j=1,2,\ldots$. }
\label{fig:illauximaps}
\end{figure}

By iterating the auxiliary maps \eqref{auxv}-\eqref{auxphi}, under a ``worst-case-scenario'' (WCS) cobwebbing  application described below, we can improve the lower and upper bounds for all trajectories of the composite map $\mathcal{M}$ \eqref {compositemapM}  within repeated updates for the bounds on the attracting domain.

{Figure~\ref{fig:illauximaps} (a)-(b)} illustrates the construction of $\xi_{U/L}^{(1)}$ and $\eta_{U/L}^{(1)}$ used in \eqref{auxv} and \eqref{auxphi} for Case FP, with ${\cal A}_1^{(1)}={\cal R}_1^+$ and $N=1$. In the $v_k - v_{k+1}$ phase plane,  the family of curves $f_1(v_k,\phi_k)$ do not cross each other, so $\xi_U^{(1)} :=   f_1(v_k, \min(\phi_k))$ and $\xi_L^{(1)} :=   f_1(v_k, \max(\phi_k))$ for $\phi_k\in [0.2,\pi/3]$ , thus yielding closed-form expressions for $\xi_{U/L}^{(1)}$ in terms of $f_1$. In contrast for $\phi_k$, the family of curves for $g_1(v_k, \phi_k)$ with fixed $v_k$ cross each other so that the envelope for $g_1$ is found computationally from the definition of $\eta_U^{(1)}$ and $\eta_L^{(1)}$ in \eqref{auxphi}.
Auxiliary maps for  $\mathcal{R}_3$  can also be constructed using the method described in Section \ref{S6.1auxiliary}. However, since ${\cal R}_3$ is a transient region, we do not pursue its construction here but focus on using the auxiliary map in $\mathcal{R}_1^+$.

 We break the WCS cobwebbing process into three steps.\\

{\bf Step 1).} We start by considering the evolution over one WCS iterate, using the lower and upper bounds $\xi_{U/L}^{(1)}$  and $\eta_{U/L}^{(1)}$ in \eqref{auxv}-\eqref{auxphi} for the $v$ and $\phi$ components.   We apply these maps to the maximum and minimum of the upper bound curve $\xi_U^{(1)}(v_k)$ and  lower bound curve $\xi_L^{(1)}(v_k)$, respectively, over the full range of $v$ in ${\cal A}_1^{(1)}$, i.e.,  $I_{0v}^{(1)} = [ v_{0}^{\rm min},\;v_{0}^{\rm max}]$,
\begin{eqnarray}
\xi_U^{(1)}(P) \mbox{ for } P= {\rm argmax}_{v}\ {\xi_U^{(1)}(v)},\qquad \xi_L^{(1)}(Q) \mbox{ for } Q=  {\rm argmin}_{v}\ {\xi_L^{(1)}(v)} \, . \label{WCS1}
 \end{eqnarray}
 Reflecting these through the diagonal gives us the first WCS iterates $v_1$, which define a new interval $I_{1v}^{(1)} = [ v_{1}^{\rm min},\;v_{1}^{\rm max}].$ 
 By the definition \eqref{WCS1}, the first-iterate images of any other point from $I_{0v}^{(1)}$ via the auxiliary maps also fall inside the updated interval $I_{1v}^{(1)}$.  This process and the analogous iteration for $\eta_U$ and $\eta_L$ are illustrated in Fig. \ref{fig:illauximaps}, including $I_{1v}^{(1)}$ and the analogous $I_{1\phi}^{(1)}$ shown in green. We notice that both of these intervals fall within the range of $v$ and $\phi$ for ${\cal A}_1^{(1)}={\cal R}_1^+$. \\

\begin{remark1}\label{WCS attracting-domain}
A single application of the WCS as in Step 1 to the initial region $ [ v_{0}^{\rm min},\;v_{0}^{\rm max}] \times [ \phi_{0}^{\rm min},\;\phi_{0}^{\rm max}] =\mathcal{R}_1^+$, gives bounds for the attracting domain, given by   $I_{1v}^{(1)} \times
I_{1 \phi}^{(1)}$, that is, ${v_1^{\min}<v_k<v_1^{\max}} $ and $ \phi_1^{\min}<\phi_k<\phi_1^{\max}.$  
\end{remark1}
This result from Step 1 represents a conservative bound for the region that attracts all trajectories of the 2D composite map ${\cal M}$, since the other regions ${\cal R}_\ell$,
$\ell \neq 1$ are demonstrated as transient, as discussed further following Statement \ref{attracting-domain-bounds} below. This observation motivates further iterations of this type,  seeking additional reductions of the attracting domain in  Steps 2 and 3. 

\begin{figure}[htbp]
    \centering
\includegraphics[width=0.7\textwidth]{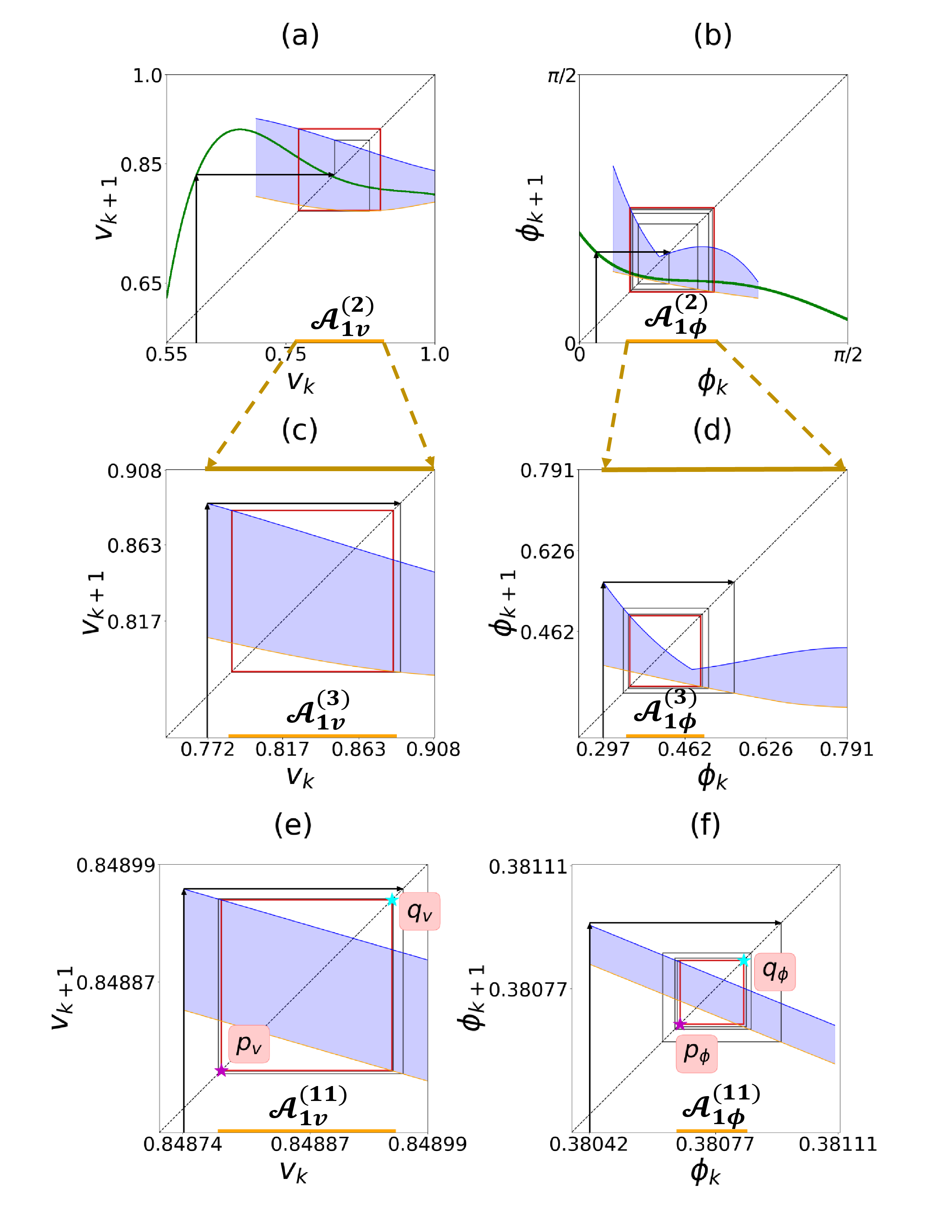}
 \caption{Illustration of the $1^{\rm st}, 2^{\rm nd}$, and $11^{\rm th}$ update of the auxiliary map ${\cal M}_{\cal A}^{(N)}$ \eqref{auxM} for Case FP ($d=0.35$). The blue and orange curves show $(\xi_U^{(N)}, \eta_{U}^{(N)})$ and $(\xi_L^{(N)}, \eta_L^{(N)})$, respectively for $\mathcal{R}_1^+$ and \eqref{auxv}-\eqref{auxphi}. The blue shaded areas between these curves  represent the possible values of $v_k$ and $\phi_k$ in ${\cal R}_1^+$.
 For each $N$, 400 steps are taken, and the last 40 steps are highlighted in red. (a)-(b): $N=1$ and ${\cal A}_1^{(1)} =\mathcal{R}_1^+$, defined in \eqref{R1+FP}. Generic cobwebbing (thin dark line) via \eqref{auxv}-\eqref{auxphi} with initial conditions ($v_0,\phi_0$) in ${\cal R}_2.$ The first few steps are governed by $(f_2,g_2)$ \eqref{r2map} (green line). In (a), the limiting behavior of \eqref{auxv} (red square) is a slight underestimate for size of ${\cal A}^{(2)}_{1v}$ (yellow) obtained via the WCS cobwebbing \eqref{WCScobweb-v}. In contrast, in (b), results from  \eqref{auxphi} and \eqref{WCScobweb-phi} are identical; (c)-(d): $N=2$, with the attracting domain from $N=1,$ $\mathcal{A}_{1v}^{(2)}$ and $\mathcal{A}_{1\phi}^{(2)}$ used as the initial domain size. 
(e)-(f): $N=11$, with the attracting domain from $N=10$ [not shown] used as the initial condition. The generic cobwebbing trajectory converges to a period-2 cycle ($p_v,q_v$ in (e) and $p_{\phi},q_{\phi}$ in (f)) that determines the size of the attracting domain: $\mathcal{A}_{1v}^{(11)}: v_k\in [0.8488,0.8490]$ and $\mathcal{A}_{1\phi}^{(11)}:\phi_k\in [0.3804,0.3811].$ Note its negligible size and the overall reduction from the original size $\mathcal{A}_{1}^{(1)}.$ Results from \eqref{auxv}-\eqref{auxphi} in  (c)-(f) yield the same attracting domain as the WCS cobwebbing \eqref{WCScobweb-v}-\eqref{WCScobweb-v}.}
 \label{cobweb12d35}
\end{figure}

{\bf Step 2).}   We repeat the procedure of Step 1 under the  WCS scenario,  applying \eqref{auxv}-\eqref{auxphi} to the max and min of the curves $\xi_U^{(1)}, \xi_L^{(1)}$  within consecutive intervals $I_{kv}^{(1)}$, $k=1, 2, \ldots$, as in \eqref{WCS1}. This repeated application can be expressed mathematically
\begin{equation}\label{WCScobweb-v}
\begin{array}{l} 
v_{k+1}^{\rm min} = \min \limits_{v_k \in I_{kv}^{(N)}} \{ \xi_L^{(N)}(v_k)\}, \\
v_{k+1}^{\rm max} = \max \limits_{v_k \in I_{kv}^{(N)}} \{ \xi_U^{(N)}(v_k)\}, \;\;\;{\rm where}\;\; I_{kv}^{(N)} = [ v_{k}^{\rm min},\;v_{k}^{\rm max} ], \qquad k=0,1,2,\ldots
\end{array}
\end{equation}
for $N=1$. These iterations are illustrated in Fig. \ref{fig:illauximaps}(c)-(d)  by the dotted black curves.
For example, since $\xi_U^{(1)}(v_k)$ is monotonically decreasing  on $I_{0v}^{(1)},$ while $\xi_L^{(1)}(v_k)$ is not, the point $P$ falls outside  $I_{1v}^{(1)},$  while $Q$ is inside.  Then, $\xi_L(Q),$ and $Q$ give the points to use in the application of \eqref{WCScobweb-v}. Then $I_{2v}^{(1)}= [v_{2}^{\rm min}=\xi_L^{(1)}(Q),\;v_{2}^{\rm max}=\xi_U^{(1)}(\xi_L^{(1)}(Q))]$. Note that the application of \eqref{WCScobweb-v} yields $v_1^{\rm min}$ and $v_2^{\rm min}$ that are identical to the previous iterate, which implies that $I_{2v}^{(1)}$  is the final estimate for the attracting domain for $v$, ${\cal B}^{(2)}_{1v}$,  starting with initial domain ${\cal R}_1^+$. Likewise iterating for $\phi_k$  from $I_{1\phi}^{(1)}$, using 
 \begin{equation}\label{WCScobweb-phi}
\begin{array}{l} 
\phi_{k+1}^{\rm min} = \min \limits_{\phi_k \in I_{k\phi}^{(N)}} \{ \eta_L^{(N)}(\phi_k)\}, \\
\phi_{k+1}^{\rm max} = \max \limits_{\phi_k \in I_{k\phi}^{(N)}} \{ \eta_U^{(N)}(\phi_k)\}, \;\;\;{\rm where}\;\; I_{k\phi}^{(N)} = [ \phi_{k}^{\rm min},\;\phi_{k}^{\rm max} ], \qquad k=0,1,2,\ldots
\end{array}
\end{equation}
for $N=1$, yields $I_{2\phi}^{(1)}= [\phi_{2}^{\rm min}=\eta_L^{(1)}(\eta_U^{(1)}(R)),\;\phi_{2}^{\rm max}=\eta_U^{(1)}(\eta_L^{(1)}(O))]\subset I_{1\phi}^{(1)}$, following from the shape of the curves $\eta_U^{(1)}$ and $\eta_L^{(1)}$. 
The repeated application of \eqref{WCScobweb-phi}  converges to its period-two cycle that yields the approximation for the size of the attracting domain for $\phi$ given by  ${\cal B}^{(1)}_{1\phi}$. Then, the reduction of ${\cal R}_1^+$ obtained through WCS cobwebbing is ${\cal B}_1^{(2)} ={\cal B}^{(2)}_{1v} \times {\cal B}^{(2)}_{1\phi}$. 
 
  The results of Step 2 are specific to starting with an initial domain ${\cal A}_{1}^{(1)} = {\cal R}_1^+ $ which defines $\xi_U^{(1)}$ and $\xi_L^{(1)}$,  used at each iteration in \eqref{WCScobweb-v} and \eqref{WCScobweb-phi}. Since the iterations yield ${\cal B}_{1}^{(2)}\subset {\cal R}_1^+$, this suggests that additional reductions for the bounds on the attracting domain may be obtained by updating the bounds $\xi_{U,L}^{(2)}$ and $\eta_{U,L}^{(2)}$, using \eqref{An}-\eqref{Bn}, i.e., using ${\cal A}_{1}^{(2)} = {\cal B}_{1}^{(2)} \subset {\cal R}_1^+$. 
 This leads to Step 3, which we write generically for the $N^{\rm th}$ update:\\
 
 {\bf Step 3).} Define an updated initial region  ${\cal A}_1^{(N)} = {\cal B}_1^{N}$, obtained via  \eqref{An}-\eqref{Bn}, with corresponding updates for $\xi_{U,L}^{(N)}$ and $\eta_{U,L}^{(N)}$. Then, repeated application of \eqref{WCScobweb-v} and \eqref{WCScobweb-phi} yields iterates that converge to a 2-cycle. The values of this 2-cycle then give new bounds on the attracting domain, denoted by $ {\cal B}_1^{N+1}$. 

In the remainder of this section Steps 1-3 are applied to the FP, PD, and CD cases, to illustrate the results of the WCS auxiliary map \eqref{WCScobweb-v}-\eqref{WCScobweb-phi}.  The implications for the attracting domain are discussed in Section \ref{S6.3GD}.   

\begin{figure}[htbp]
    \centering
    \includegraphics[width=\textwidth]{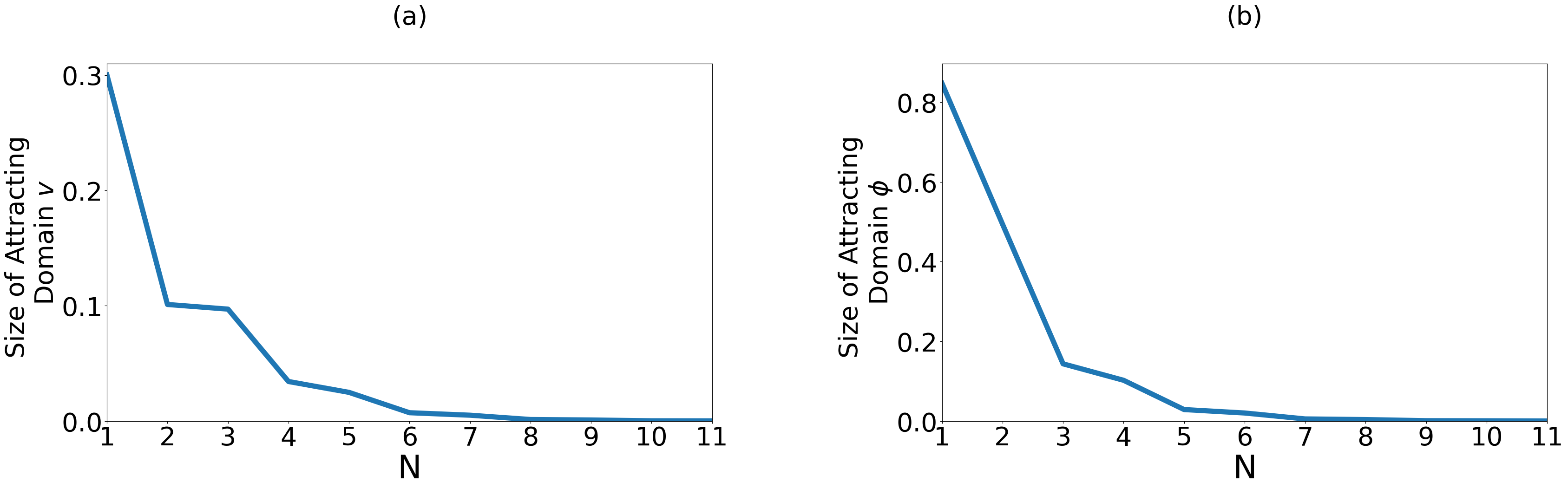}
 \caption{Illustration of the size of the domain ${\cal A}_N$ for each $N$, showing that the attracting domain size decreases monotonically for Case FP, reaching 0.000185 and 0.0001867 in the  $v_k,\phi_k$ directions, respectively. }
 \label{cobweb12d35_updates}
\end{figure}

\begin{figure}[htbp]
    \centering
\includegraphics[width=0.8\textwidth]{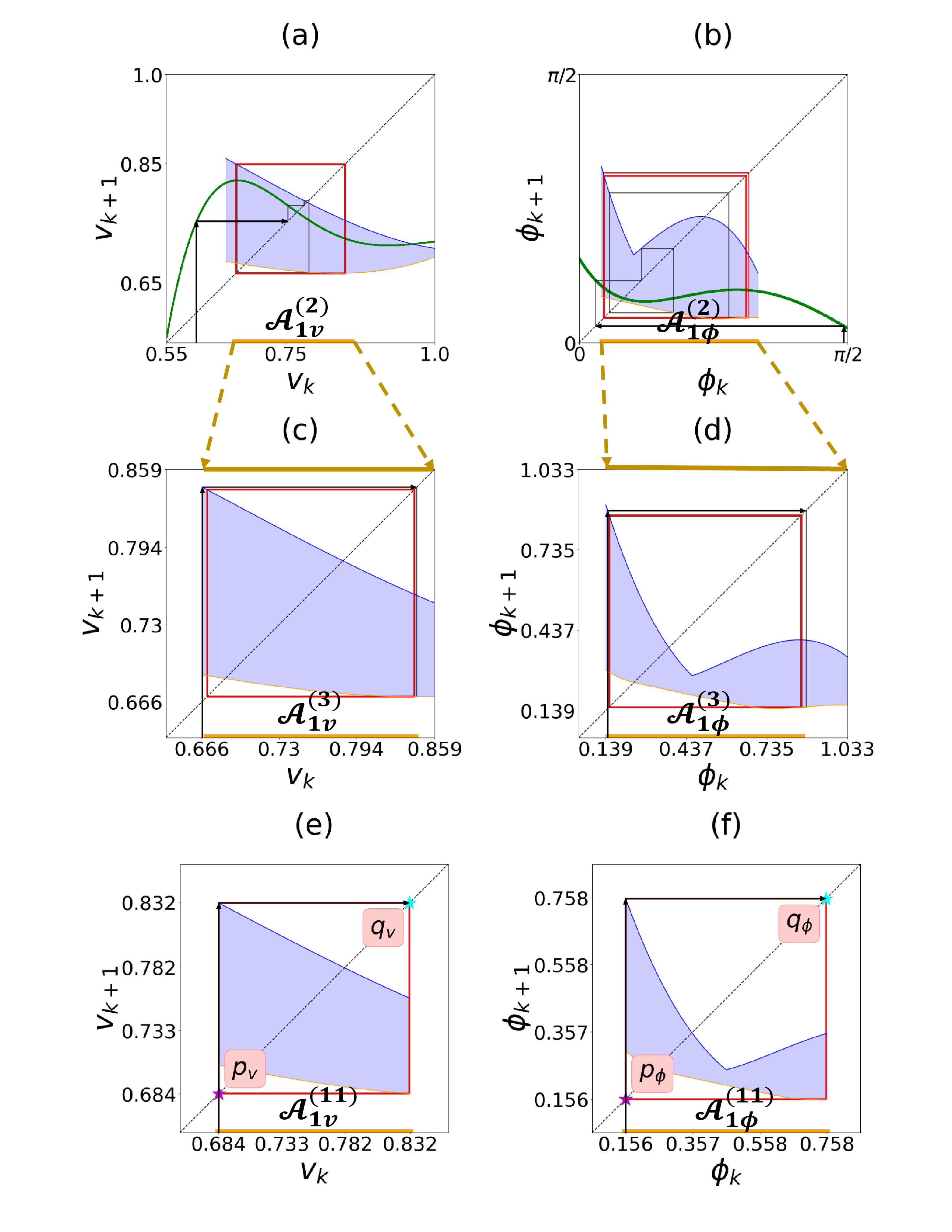}
\caption{ Illustration of the $1^{\rm st}, 2^{\rm nd}$, and $11^{\rm th}$ update of the auxiliary map ${\cal M}_{\cal A}^{(N)}$ \eqref{auxM}, for Case PD ($d=0.30$), using the same procedure as in Fig.~\ref{cobweb12d35}. 
(a)-(b): $N=1$ and $\mathcal{A}_1^{(1)}=\mathcal{R}_1^+$  \eqref{R1+PD}.
Generic cobwebbing (thin dark line) via \eqref{auxv}-\eqref{auxphi} with initial conditions ($v_0,\phi_0$) in ${\cal R}_2.$ Similar to Fig.~\ref{cobweb12d35}(a), in (a) the  limiting behavior via \eqref{auxv}-\eqref{auxphi} (red square) is a slight underestimate for the size of  the actual size of ${\cal A}^{(2)}_{1v}$ (yellow), obtained via the WCS cobwebbing \eqref{WCScobweb-v}. In (b), they are identical.
(c)-(d): $N=2$ and $\mathcal{A}_1^{(2)}: v_k\in [0.666,0.850]$ and $\phi_k\in [0.146,0.977].$ (e)-(f): $N=11$ and $\mathcal{A}_1^{(11)}: v_k\in [0.684,0.832]$ and $\phi_k\in [0.156,0.758]$, where the size of $\mathcal{A}_1^{(N)}$ for $N>1$ follows directly from the limiting (yellow) behavior in the $(N-1)^{\rm th}$ update (\eqref{An}-\eqref{Bn}). The stars with ($p_v,q_v$) and ($p_\phi,q_\phi$)  in panels (e) and (f) indicate the min and max of the period-2 cycle.}
\label{Auximapsd30}
\end{figure}

 While \eqref{WCScobweb-v}-\eqref{WCScobweb-phi} refines the generic \eqref{auxv}-\eqref{auxphi} with a WCS choice of $(v_k,\phi_k)$ on each iteration,   the two cobwebbing approaches  are equivalent in some cases.
 This property is determined by the shape of $\xi_{U,L}^{(N)}$ and $\eta_{U,L}^{(N)}$  (blue and orange curves in Fig. \ref{fig:illauximaps}). 
 For example,  if $v_{k+1}^{\rm min} >v_{k}^{\rm min}$ and $v_{k+1}^{\rm max} <v_{k}^{\rm max}$, as in the case of   monotonically decreasing $\xi_{U,L}^{(N)}$ and $\eta_{U,L}^{(N)}$, both yield the same result. 
This observation is useful, since cobwebbing based on the generic auxiliary map  \eqref{auxv}-\eqref{auxphi} is more straightforward to implement computationally since it does not restrict its application to the maximum or minimum on $\xi_{U,L}^{(N)}$ and $\eta_{U,L}^{(N)}$  curves, as in the WCS approach \eqref{WCScobweb-v} -\eqref{WCScobweb-phi}. However, in the WCS treatment of general functions $\xi_{U,L}^{(N)}$ and $\eta_{U,L}^{(N)}$, the bounds $v_{1}^{\rm min/max}$ and/or $\phi_{1}^{\rm min/max}$ might not be improved with further iterates of  \eqref{WCScobweb-v} and/or \eqref{WCScobweb-phi} for fixed $N$. In such cases  \eqref{auxv}-\eqref{auxphi} may underestimate the size of the attracting domain, since it does not restrict its application to the maximum or minimum on the upper or lower curves, which the WCS approach takes into account. For example, Figure~\ref{cobweb12d35}(a) shows results from  \eqref{auxv}-\eqref{auxphi}  iterating from a generic $v_k$, obtained from a random initial condition, with a limiting period-2 cycle in red.  This is slightly smaller than the limiting results from WCS cobwebbing \eqref{WCScobweb-v} and/or \eqref{WCScobweb-phi}   (yellow bar).   In contrast, in  Fig.~\ref{cobweb12d35}(b), the intervals are identical, due to the property that $\eta_L^{(1)}(\phi_k)$ and $\eta_U^{(1)}(\phi_k)$  are monotonically decreasing outside ${\cal A}^{(2)}_{1\phi}$, so that the WCS cobwebbing procedure repeatedly excludes the previous global maximum and minimum over ${\cal A}_{1\phi}^{(1)}$. As discussed further below, the generic approach achieves the same result as WCS in most of the FP, PD, and CD cases considered here, particularly when applying Step 3, that is, additional $N$ updates of the interval ${\cal A}_1^{(N)}$.  
 
Figure~\ref{cobweb12d35}(c)-(f) illustrates these updates of region ${\cal A}_1^{(N)}$ and  ${\cal M}_{\cal A}^{(N)}.$ Each row shows results for a different update, specifically for $N=2$ and $N=11$.  The red box highlights the last 10\% of the cobweb iterations,  indicating the limiting dynamics for ${\cal M}_{\cal A}^{(N)}$.  The size of the corresponding  {attracting domain (indicated by the yellow interval)} shrinks with $N$, and ${\cal A}_1^{(N)}$ for $N>1$ is updated accordingly, as in \eqref{An}-\eqref{Bn}.  For example, in Fig.~\ref{cobweb12d35}(a)(b), $\mathcal{A}_1^{\rm (1)}=\mathcal{R}_1^+$ for $N=1$ with $v_k\in [0.7, 1]$ and $\phi_k\in[0.2,\pi/3]$, and the limiting range shown by the yellow interval is  $v_k\in [0.771, 0.909]$ and $\phi_k\in[0.297,0.791]$. Continuing with this process for increasing $N$, Figs. \ref{cobweb12d35}(c),(f) and (e),(f) illustrate the smaller range of $v_k$ and $\phi_k$ given by $\xi_{U/L}^{(N)}$ and $\eta_{U/L}^{(N)}$, mirroring  the smaller size of ${\cal A}_1^{(N)}$ with increasing $N$. As shown at $N=11$, $\mathcal{A}_1^{\rm (11)}$ is significantly smaller than $\mathcal{A}_1^{\rm (1)}$. This contraction property with increasing $N$ is summarized in Fig. \ref{cobweb12d35_updates},  which shows how the length and width of the  {attracting domain} for  $v_k$ and $\phi_k$ decreases with increasing $N$. Thus, even though the max/min characteristics of the auxiliary map do not allow the limiting behavior of ${\cal M}_{\cal A}$ to be a fixed point, nevertheless, for Case FP, we see that region ${\cal A}_1^{(N)}$ shrinks to a negligible size for large $N$.

\begin{figure}[htbp]
    \centering
\includegraphics[width=\textwidth]{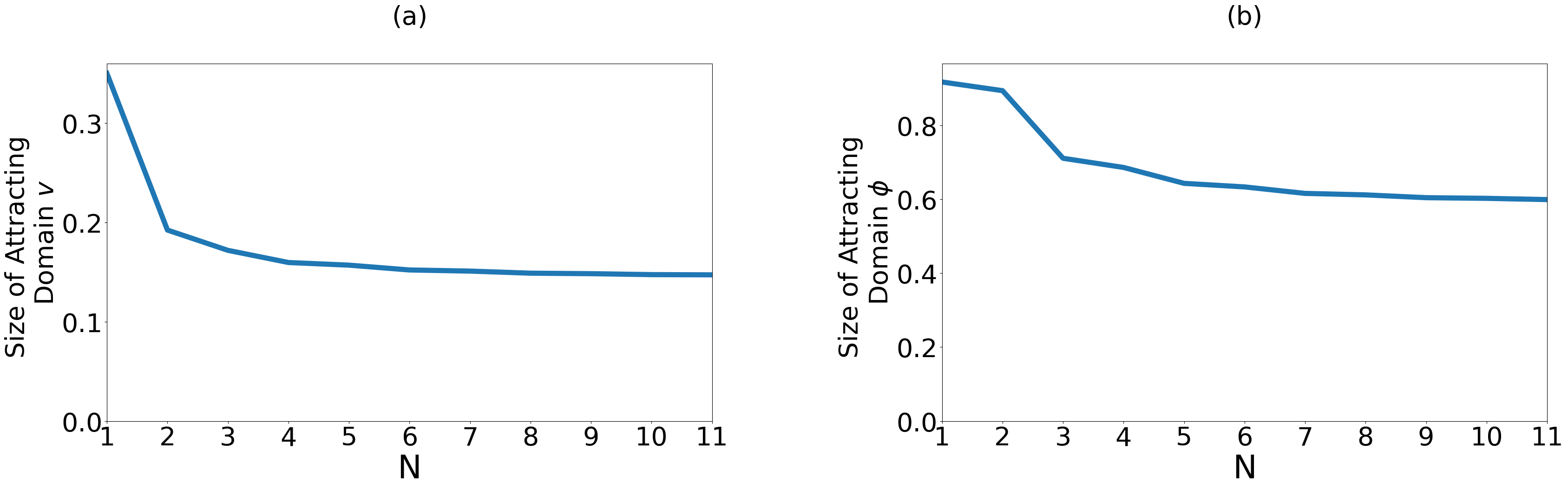}
\caption{ Illustration of the {attracting domain size} for case PD that decreases to a limiting size, with the final limiting size as 0.1472 and 0.5991 for $v$ and $\phi$, respectively. }
\label{Auximapsd30_updates}
\end{figure}

\begin{figure}[htbp]
    \centering
\begin{subfigure}{0.7\textwidth}
    \centering
\includegraphics[width=0.9\textwidth]{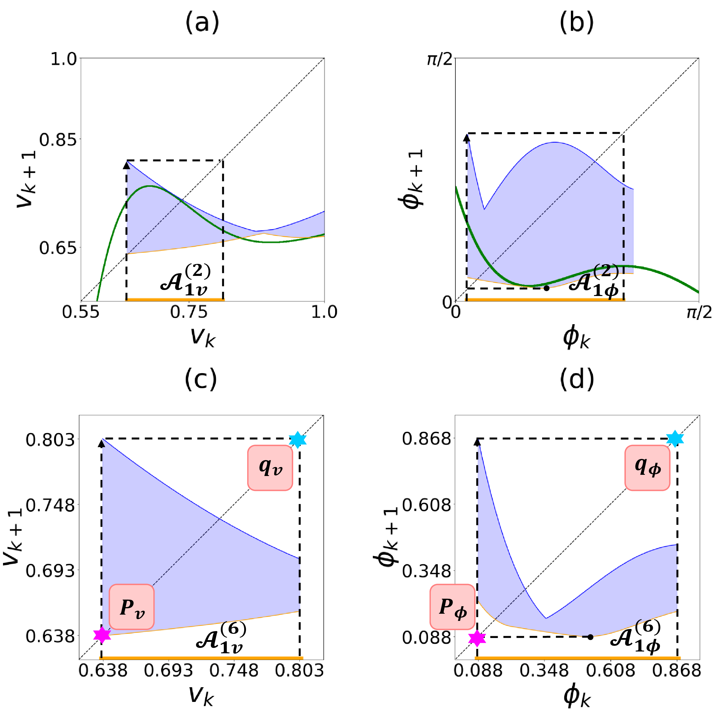}
\end{subfigure}
\vfill
\begin{subfigure}{0.9\textwidth}
    \centering
\includegraphics[width=0.8\textwidth]{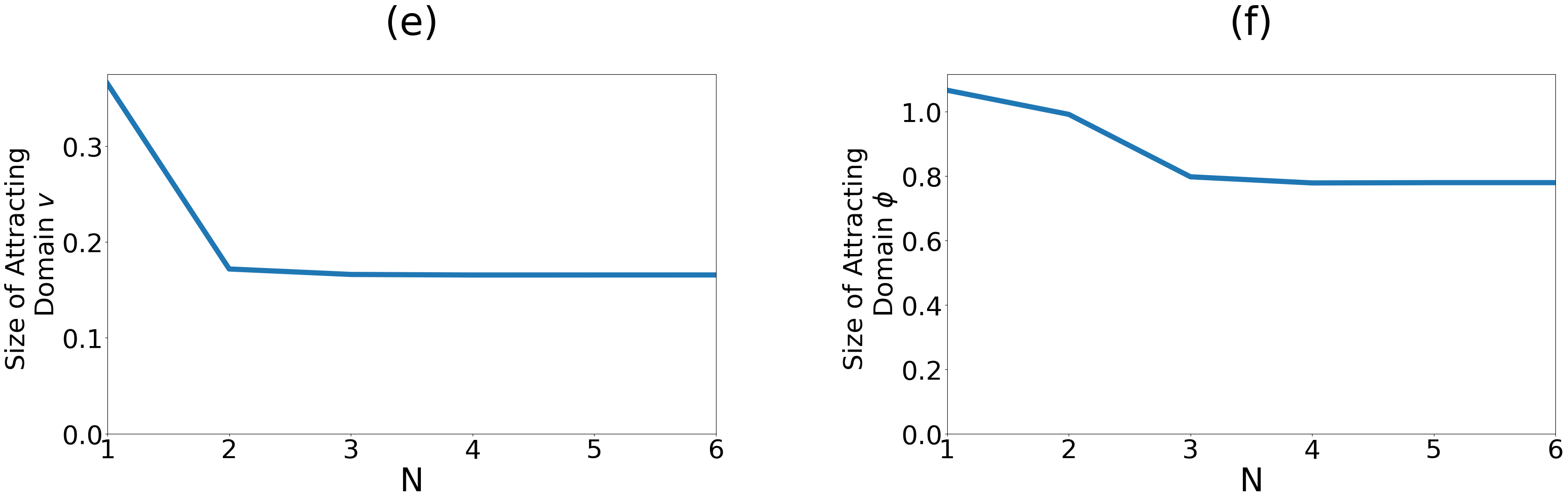}
\end{subfigure}
 \caption{ Illustration of the $1^{\rm st}$ and $6^{\rm th}$ update of the auxiliary map ${\cal M}_{\cal A}^{(N)}$ \eqref{auxM}, for $d=0.26$, corresponding to Case CD, using the same procedure as in Fig.~\ref{cobweb12d35}. {In panels (a)-(b), the attracting domain (dashed line) is calculated using one iterate of the WCS cobwebbing \eqref{WCScobweb-v}-\eqref{WCScobweb-phi} to obtain ${\cal A}_{1v}^{(N)}$ and  ${\cal A}_{1\phi}^{(N)}$ (yellow interval). Further updates over $N$ offer a slight reduction, with the limiting size shown in (c) and (d). Here, $\mathcal{A}_1^{(2)}: v_k\in [0.673,0.789];\; \phi_k\in [0.093,0.725]$ and $\mathcal{A}_1^{(6)}: v_k\in [0.638,0.803];\; \phi_k\in [0.088, 0.868]$ The stars with ($p_v,q_v$) and ($p_\phi,q_\phi$) in panels (c) and (d) indicate the min and max of the period-2 cycle of the WCS iterative procedure \eqref{WCScobweb-v}-\eqref{WCScobweb-phi}.} Panels (e) and (f) show the decrease of the  {attracting domain size} to a limiting size with the limiting size equal to $0.166$ and $0.780$ for $v$ and $\phi$, respectively.}
\label{Auximapsd26}
\end{figure}

 We also apply the auxiliary map method to estimate the attracting domain for non-trivial dynamics in Case PD and Case CD and observe the contraction property from Fig. \ref{Auximapsd30}-\ref{Auximapsd26} for Case PD and CD. 

Similar to the cobweb illustration of the updates in the Case FP, Fig.~\ref{Auximapsd30} and Fig.~\ref{Auximapsd26} illustrate the updates of the region ${\cal A}_1^{(N)}$  and  ${\cal M}_{\cal A}^{(N)}$ in Case PD and Case CD, respectively. The setup in Fig.~\ref{Auximapsd30} and Fig.~\ref{Auximapsd26} is the same as in Fig.~\ref{cobweb12d35}, with each row showing results from updates of ${\cal A}_1^{(N)}$. In Case PD, $N=1, N=2$, and $N=11$ are shown, while Case CD demonstrates the cobwebbing diagrams for $N=1$ and $N=6$.  Moreover, in contrast to the Case FP, where the limiting dynamics approaches a point for $N$  large,  for Cases PD and CD, the attracting domain follows from the attracting period-2 cycle, with a limiting size at a finite $N$ dictated by $|p_v-q_v|$ and $|p_{\phi}-q_{\phi}|$. The pairs of points ($p_v,q_v$) and ($p_\phi,q_\phi$) shown in Figs.~\ref{cobweb12d35},\ref{Auximapsd30},\ref{Auximapsd26} for the largest value of $N$ indicate the maximum $q_{\bullet}$ and minimum $p_{\bullet}$ of period-2 cycle for $v$ and $\phi$. Likewise, these values can be used to explicitly determine the size of the globally  {attracting domain} via {constructing second-iterate maps}, as discussed in the next section.

In Case PD, the limiting dynamics converge to an attracting period-2 cycle for both $v_k$ and $\phi_k$ when $N$ is large, with much of the size reduction of $\mathcal{A}_1^{(N)}$ occurring in the first two updates, as shown in Fig.~\ref{Auximapsd30_updates}. 
  Similar to Case PD, Fig.~\ref{Auximapsd26} shows that the limiting dynamics of Case CD for sufficiently large $N$ yield cycles that bound a relatively larger range of $v_k$ and $\phi_k$. Case CD does not allow using the generic cobwebbing via \eqref{auxv}-\eqref{auxphi} since the lower bounds $\xi^{(N)}_L(v_k)$ and $\eta^{(N)}_L(\phi_k)$ are not monotonically decreasing functions on any updated interval. 

\subsection{Second-iterate maps for the attracting domain} \label{S6.3GD}
The auxiliary map method developed in the previous subsection opens the door to explicitly characterizing the global dynamics of the composite map. Figures~\ref{cobweb12d35}-\ref{Auximapsd26} demonstrate that the final size of the attracting domain is always bounded by a period-$2$ cycle. In the scenarios where the WCS cobwebbing was applied (Fig.~\ref{cobweb12d35}(a) for the FP case, Fig.~\ref{Auximapsd30}(a) for the PD case, and Fig.~\ref{Auximapsd26}(a)-(d) for the CD case), this orbit is a $2$-cycle of the WCS cobwebbing iterative procedure \eqref{WCScobweb-v}-\eqref{WCScobweb-phi}. In the scenarios where the generic cobwebbing  \eqref{auxv}-\eqref{auxphi} yields the same result as WCS (Fig.~\ref{cobweb12d35}(b)-(f) for the FP case, Fig.~\ref{Auximapsd30}(b)-(f) for the PD case), this orbit is a $2$-cycle of the auxiliary map based on alternating upper and lower bound curves in \eqref{auxv}-\eqref{auxphi}.

Since these $2$-cycles bound a subset of the auxiliary map's phase space, their existence and global stability imply the existence of a globally stable  {attracting domain} for the trajectories of the composite map $\mathcal{M}$ \eqref {compositemapM}. The bounds on the  {attracting domain} are  indicated as $q_v,p_v,q_\phi,$ and $p_\phi$ in Figs.~\ref{cobweb12d35},  \ref{Auximapsd30}, and \ref{Auximapsd26}  for the largest value of $N$ shown.
When applicable, computing these values as the roots of $m$ iterations of the maps \eqref{auxv} and \eqref{auxphi} for appropriate $m$, we obtain their stability and thus bounds on the attracting domain for the dynamics.

First, to obtain the bounds on $v_k$ and $\phi_k$  used in the  $(N+1)^{\rm th}$ update, we consider the general second-iterate WCS maps for $v_{k+2}^{\max}$ and $\phi_{k+2}^{\max}$, defined via \eqref{WCScobweb-v} and \eqref{WCScobweb-phi}, respectively: 
\begin{align}\label{WCSsecond-v}
v_{k+2}^{\rm min} = v_{k+1}^{\rm max}(v_{k}^{\rm min}),\;\;\;{\rm where}\;\; v_{k}^{\rm min}=\min \limits_{v_{k-1} \in I_{(k-1)v}^{(N)}} \{ \xi_L^{(N)}(v_{k-1})\}\;{\rm and}\; v_{k+1}^{\rm max}=\max \limits_{v_k \in I_{kv}^{(N)}} \{ \xi_U^{(N)}(v_k)\},\\
\label{WCSsecond-phi}
\phi_{k+2}^{\rm min} = \phi_{k+1}^{\rm max}(\phi_{k}^{\rm min}),\;\;\;{\rm where}\;\; \phi_{k}^{\rm min}=\min \limits_{\phi_{k-1} \in I_{(k-1)\phi}^{(N)}} \{ \eta_L^{(N)}(\phi_{k-1})\}\;{\rm and}\; \phi_{k+1}^{\rm max}=\max \limits_{\phi_k \in I_{k\phi}^{(N)}} \{ \eta_U^{(N)}(\phi_k)\}.
\end{align}

In cases where the WCS cobwebbing via \eqref{WCScobweb-v} is equivalent to the generic cobwebbing using \eqref{auxv} as detailed in Subsection~\ref{S6.2}, the second-iterate WCS map \eqref{WCSsecond-v} for $v_{k+2}^{\rm min}$ transforms into the second-iterate map:
\begin{eqnarray}\label{auxv2nd}
v_{k+2}(v_k) & = \xi_{L}^{(N)}\big(\xi_{U}^{(N)}(v_k) 
\big).
\end{eqnarray}

The maps $\xi^{(N)}_{L/U}$ are written explicitly in terms of $f_1$ evaluated at $\phi^{\rm min/max}_{0}$. 
They do not cross each other, analogous to $f_1$  shown in Fig.~\ref{fig:illauximaps}(a). 
Then, we have the closed-form expression for {the first-iterate \eqref{auxv} and  second-iterate map \eqref{auxv2nd}, where the latter} for $v_{k+2}$ is a $9^{\rm th}$-order polynomial of the form 
\begin{align}
    v_{k+2}(v_k) &= f_1(f_1(v_k, \phi^{\rm max}_0), \phi^{\rm min}_0) \nonumber\\
    &= \alpha_0 + \alpha_1 v_k^1 + \alpha_2 v_k^2 + \alpha_3 v_k^3 + \alpha_4 v_k^4 + \alpha_5 v_k^5 + \alpha_6 v_k^6 + \alpha_7 v_k^7 + \alpha_8 v_k^8 + \alpha_9 v_k^9 \, .
    \label{seconditerativev}
\end{align}
Here, $\alpha_i,i=1,...,9$ are polynomials that depend on $d$ and on $\phi^{\rm min}_0$ and $\phi^{\rm max}_0$,  whose coefficients $b_0, b_1, ..., b_9$ are listed in Supplementary Section III.
The (stable) root $v_{k+2}=v_k=p_v$ of \eqref{seconditerativev} corresponds to the minimum on the limiting behavior of  $\xi_1^{(N)}$ \eqref{auxv}, with the maximum $q_v$ obtained by

\begin{eqnarray}
& v_k = p_v,    \qquad 
v_{k+1} = q_v = f_1(v_k, \phi^{\rm min}_0) = f_1(p_v, \phi^{\rm min}_0)  = \xi_{U}^{\rm (N)} (p_v),  \label{pvqv_modified} \\
\implies & v_{k+2} = p_v = f_1(v_{k+1},\phi^{\rm max}_0) = f_1(q_v,\phi^{\rm max}_0)  =f_1(f_1(p_v,\phi^{\rm min}_0),\phi^{\rm max}_0) = \xi_{L}^{\rm (N)}(p_v)\, . \nonumber 
\end{eqnarray}

These values $p_v$ and $q_v$ determine explicit bounds for the attracting domain for $v_k$ indicated by the red boxes for sufficiently large $N$ in Figs.~\ref{cobweb12d35}(e) and Fig.~\ref{Auximapsd30}(e) for the FP and PD cases that allow using \eqref{seconditerativev}. Note that deriving such tight bounds for global dynamics directly from the 2D composite map ${\cal M}_{\cal A}$ via a Lyapunov function or similar approaches without the constructive use of the explicit auxiliary maps 
seems elusive. 

 Similarly, the {periodic solutions} for $\phi_k$ are based on the definition of $\eta_1^{(N)}$ in \eqref{auxphi}. For the FP and PD cases, we consider 
\begin{eqnarray}\label{auxphi2nd}
\phi_{k+2}(\phi_k) = \eta_{L}^{(N)}\big(\eta_{U}^{(N)}(\phi_k) \big).
\end{eqnarray}
In contrast to \eqref{seconditerativev} for $v_k$, the family of curves $g_1(v_k,\phi_k)$,  in the definition of $\eta_{L/U}$ \eqref{auxphi} cross each other for different fixed $v_k\in [v^{\rm min}_0,v^{\rm max}_0]$, analogous to 
Fig.~\ref{fig:illauximaps}(b). Then, there is no closed-form expression for the first- and second-iterate maps $\phi_{k+1}$ and $\phi_{k+2}$, and  $\eta_{L/U}$ are determined numerically in \eqref{auxphi2nd}. 

\begin{figure}[htbp]
    \centering
\includegraphics[width=0.7\textwidth]{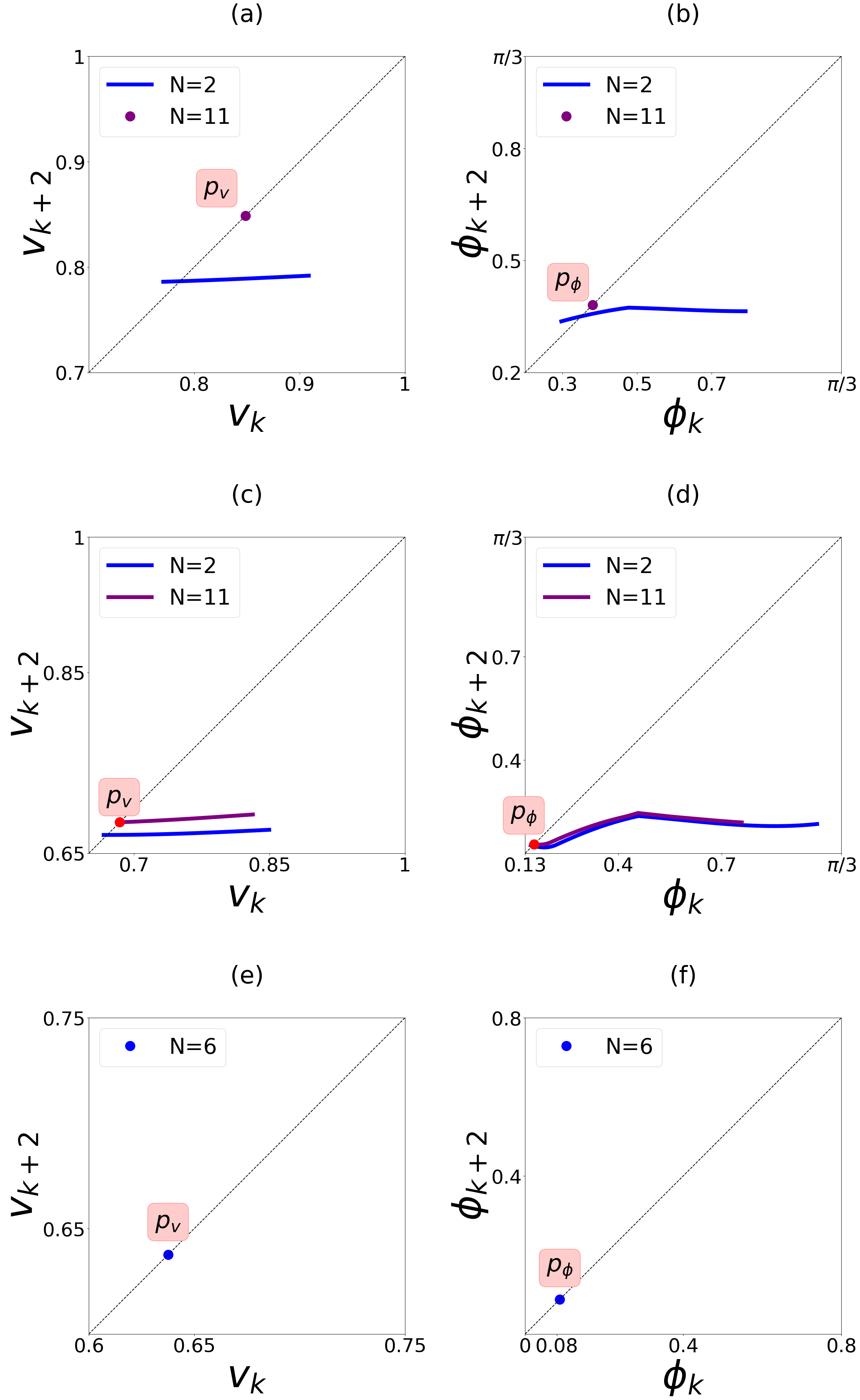}
\caption{Second-iterate maps for the $N$th update of the auxiliary maps. The curves are obtained from \eqref{auxv2nd} and \eqref{auxphi2nd}.
(a)-(b):  The FP case for $N= 2,11$. The blue curve corresponds to the $2$nd update. The purple curve for the $11$th update is invisible due to the strong contraction of the attracting domain after the updates. (c)-(d): The PD case  for $N= 2,11$. (e)-(f). The CD case for $N=6$. The purple (red) points $p_v$ and $p_{\phi}$ are fixed points of the second-iterate maps \eqref{auxv2nd} and \eqref{auxphi2nd} in the FP  (PD) cases and correspond to 2-cycles, defined in \eqref{pvqv_modified}-\eqref{pphiqphi_modified} and depicted in Fig.~\ref{cobweb12d35}(e)-(f) and 
Fig.~\ref{Auximapsd30}(e)-(f), respectively. The blue fixed points in (e)-(f) correspond to 2-cycles of the WCS auxiliary maps \eqref{WCSsecond-v}-\eqref{WCSsecond-phi}, also shown in Fig.~\ref{Auximapsd26}(c)-(d).}
\label{highorderiterate}
\end{figure}

 For the FP and PD cases, we calculate $p_\phi$ and $q_\phi$, which give the minimum and maximum of the limiting behavior shown by the red boxes in Fig.~\ref{cobweb12d35}(f) and  Fig.~\ref{Auximapsd30}(f) for sufficiently large $N$.  They are given by
\begin{eqnarray}
 & \phi_k = p_\phi, \qquad 
 \phi_{k+1} = q_\phi = \max_{v_k} g_1(v_k, \phi_k) = \max_{v_k}g_1(v_k, p_\phi) = \eta_{U}^{\rm (N)}(p_\phi),\label{pphiqphi_modified}\\
\implies & \phi_{k+2} = p_\phi = \min_{v_k} g_1(v_k, \phi_{k+1}) = \min_{v_k} g_1(v_k, q_\phi) =  \min_{v_k} g_1(v_k, \max_{v_k}g_1(v_k, p_\phi)) = \eta_{L}^{\rm (N)}(\eta_{U}^{\rm (N)}(p_\phi)).
\nonumber
\end{eqnarray}
 The curves obtained from applying the iterates given in \eqref{seconditerativev} and \eqref{auxphi2nd} are shown in Fig.~\ref{highorderiterate}(a-d) in the FP and PD cases. Panels (a)-(d) illustrate the stability of the fixed points $p_v$ and $p_\phi$ for the period-2 cycle. There, the curves show the limiting behavior of the second iterate of ${\cal M}_{\cal A}^{(N)}$, given by \eqref{auxv2nd} and \eqref{auxphi2nd}. They intersect the diagonals in the $v_{k+2}-v_k$ and $\phi_{k+2}-\phi_k$ phase planes with a slope less than unity.   Then, for sufficiently large $N$ we obtain the stable fixed points $p_v$ and $p_{\phi}$, likewise implying the stability of the fixed points  $q_v$ and $q_{\phi}$,  which all together provide the range of the  {attracting domain} for ${\cal M}_{\cal A}^{(N)}$ in Fig.~\ref{cobweb12d35} and Fig.~\ref{Auximapsd30}.  

In the CD case where the lower bound function $\xi_L^{(N)}(v_k)$ does not allow  using \eqref{seconditerativev} and \eqref{auxphi2nd} even for sufficiently large updates $N$, the WCS iterations  \eqref{WCSsecond-v} and \eqref{WCSsecond-phi} yield a $2$-cycle for $v_{k+2}^{\min}$ and $\phi_{k+2}^{\min}$, respectively:
\begin{align}\label{WCSsecond-CD-v}
v_{k+2}^{\rm min}=p_v=\min_{v_{k+1}} \{ f_1(v_{k+1},\phi^{\rm max}_0)\},\;
v_{k+1}^{\rm max}=q_v=\max_{v_k} \{ f_1(v_{k},\phi^{\rm min}_0)\},\\
\label{WCSsecond-CD-phi}
\phi_{k+2}^{\rm min} =p_{\phi}= \min_{\phi_{k+1}} \{ \eta_L^{(N)}(\phi_{k+1})\};\; \phi_{k+1}^{\rm max}=q_{\phi}=\max_{\phi_{k}} \{ \eta_U^{(N)}(\phi_{k})\}.
\end{align}
Figure~\ref{Auximapsd26}(c),(d) illustrates the 2-cycle for $v_{k}^{\min}$ and $\phi_{k}^{\min}$, respectively. Additionally, Fig.~\ref{highorderiterate}(e),(f) shows the fixed points $p_v$ and $p_{\phi}$ of the second-iterate    WCS map \eqref{WCSsecond-v} and \eqref{WCSsecond-phi} that provide the lower bound for the attracting domain whose upper bounds, $q_v$ and $q_{\phi},$ can similarly be identified from \eqref{WCSsecond-CD-v}-\eqref{WCSsecond-CD-phi}.

The following statement summarizes the results for the existence of a globally attracting domain for the auxiliary composite map ${\cal M}_{\cal A}^{(N)}$  highlighting cases where the generic and WCS iterations provide the same results. \\

\begin{statement1}\label{attracting-domain-bounds}
[Bounds for the attracting domain]. \\ 
1).  The generic \eqref{auxv}-\eqref{auxphi} and WCS
\eqref{WCScobweb-v} -\eqref{WCScobweb-phi} iterates yield the same upper and lower bounds functions $\xi_{U,L}^{(N)}$ and $\eta_{U,L}^{(N)}$  in the case that the auxiliary map $\mathcal{M}_{\cal A}^{(N)}$ in $A_1^{(N)}\in R^+_1$ satisfies the property on each $k^{\rm th}$ iterate,
\begin{equation}
\begin{array}{l}
{\rm argmax}_{v_{k-1}\in I^{(N)}_{(k-1)v}} \xi_{U}(v_{k-1}) \notin I^{(N)}_{kv}=[v_k^{\min},\;v_k^{\max}],\\
{\rm argmin}_{v_{k-1}\in I^{(N)}_{(k-1)v}} \xi_{L}(v_{k-1}) \notin I^{(N)}_{kv},
\end{array}\label{part1-v}
\end{equation}
\begin{equation}
\begin{array}{l}
{\rm argmax}_{\phi_{k-1}\in I^{(N)}_{(k-1)\phi}} \eta_{U}(\phi_{k-1}) \notin I^{(N)}_{k \phi}=[\phi_k^{\min},\;\phi_k^{\max}],\\
{\rm argmin}_{\phi_{k-1}\in I^{(N)}_{(k-1)\phi}} \eta_{L}(\phi_{k-1}) \notin I^{(N)}_{k\phi} \,. 
\end{array}\label{part1-phi} 
\end{equation}
Then, each successive application of either of these auxiliary maps yields new global maxima of $\xi_{U}^{(N)}$ and $\eta_{U}^{(N)}$ and global minima of $\xi_{L}^{(N)}$ and $\eta_{L}^{(N)}$ on the updated $v$- and $\phi$-intervals,  $I^{(N)}_{k\phi}$ and  $I^{(N)}_{k\phi}$.
Then, the auxiliary map $\mathcal{M}_{\cal A}^{(N)}$ can be defined by  \eqref{auxv}-\eqref{auxphi} and has a  stable $2$-cycle with alternating values $p_v$ and $q_v$ for $v_k$ and $p_{\phi}$ and $q_{\phi}$ for $\phi_k,$ given in \eqref{pvqv_modified}-\eqref{pphiqphi_modified}. These values determine the bounds for the attracting domain $ A_1^{(N+1)}=\{p_v<v_k<q_v,\;p_{\phi}<\phi_k<q_{\phi}\}$ in the composite 2D map ${\cal M}$ \eqref{auxM} in region $\mathcal{R}_1^+.$ \\
2). The auxiliary map $\mathcal{M}_{\cal A}^{(N)}$ is defined by the WCS map \eqref{WCScobweb-v}-\eqref{WCScobweb-phi}  if the  $k^{\rm th}$ iterate of the upper and lower bound functions $\xi_{U/L}^{(N)}$ and $\eta_{U/L}^{(N)}$ in $A_1^{(N)}\in R^+_1$ does not shrink (or expand) the $v$-interval and $\phi$-interval, i.e.,  $I^{(N)}_{kv} =I^{(N)}_{(k+1)v}$ and $I^{(N)}_{k\phi} =I^{(N)}_{(k+1)\phi}$. 
Again, the auxiliary map $\mathcal{M}_{\cal A}^{(N)}$ defined via the WCS map  has a  stable $2$-cycle with alternating values $p_v$ and $q_v$ for $v_k$ and $p_{\phi}$ and $q_{\phi}$ for $\phi_k,$ given in \eqref{WCSsecond-CD-v}-\eqref{WCSsecond-CD-phi}. These values determine the  bounds for the attracting domain $A_1^{(N+1)}$ in the composite 2D map ${\cal M}$ \eqref{auxM} in the region $\mathcal{R}_1^+$.
\end{statement1}

\begin{remark1}\label{proof} Parts 1) and 2) follow directly from the construction of the auxiliary map $\mathcal{M}_{\cal A}^{(N)}$, either via the generic   \eqref{auxv}-\eqref{auxphi}, or via WCS  iterates \eqref{WCScobweb-v}-\eqref{WCScobweb-phi}.  When \eqref{part1-v}-\eqref{part1-phi} holds,   e.g.,  when upper and lower bound functions $\xi_{U/L}^{(N)}$ and $\eta_{U/L}^{(N)}$ are monotonically decreasing in $A_1^{(N)}\in R^+_1$, the generic   \eqref{auxv}-\eqref{auxphi} suffices. 
In the more general setting of attracting ${\cal R}_1^+$, the most conservative one-iterate bound for ${\cal A}_1^{(N)}$ (see Remark \ref{WCS attracting-domain}) may be improved by further WCS iterates, i.e., updates via Steps 2-3 in Section \ref{S6.2}.  
In either case, the iterates converge to a stable 2-cycle,  alternating between $\xi_{U/L}^{(N)}$ and between $\eta_{U/L}^{(N)}$, preventing the emergence of higher-period orbits due either to the contraction condition  \eqref{part1-v}-\eqref{part1-phi}  or to the WCS formulation  \eqref{WCSsecond-CD-v}-\eqref{WCSsecond-CD-phi}.  
 \end{remark1}

Note that the FP and PD cases at $N\ge 2$ satisfy the condition \eqref{part1-v}-\eqref{part1-phi} of Statement \ref{attracting-domain-bounds} (Fig.~\ref{cobweb12d35}(c)-(f) for the FP case, Fig.~\ref{Auximapsd30}(c)-(f) for the PD case). In contrast, the upper and lower bound functions in the CD case do not obey \eqref{part1-v}-\eqref{part1-phi} (Fig.~\ref{Auximapsd26}), so the attracting domain in the CD case is determined by the conditions of Part 2) in Statement \ref{attracting-domain-bounds}.
 
As described in Section \ref{S6.1auxiliary}, one can apply the auxiliary approach for all regions ${\cal R}_j$ for $j=2,3,4,5$, which confirms the transient behavior for regions outside of ${\cal R}_1 $.   Combining this transient behavior with the results of this section, we have the complete confirmation of the bounds on the attracting domains for ${\cal M}$ for different $d$,   obtained via the limiting regions of the auxiliary map as applied in Sections \ref{S6.2}, \ref{S6.3GD}.

\section{Conclusion}\label{S7conclusion}
While studying VI systems through local stability analysis has gained significant momentum, understanding their global dynamics and bifurcations remains challenging due to the limited applicability of classical global stability methods developed for smooth dynamical systems. In particular, the engineering literature has focused on linear stability and bifurcations, yet global behavior is important in design.

In this paper, we propose a computer-assisted analysis based on reduced smooth maps for studying the global dynamics of the VI pair.  The framework is designed to be generic, ideally for application to other non-smooth dynamical systems. The global stability analysis is facilitated by an approximation of the exact map for the states at impact, specifically the relative impact velocity $\dot{Z}_k$ between the outer (the capsule) and the inner (the ball) masses and the impact phase $\psi_k$ relative to the forcing. The exact non-smooth maps for these quantities are given by complex coupled transcendental equations for $\dot{Z}_k$ and $\psi_k$. While the non-smooth dynamics present a challenge in using commonly defined maps,  they also provide an opportunity for designing a new approach for impacting systems. Specifically, we use short sequences of returns to one side of the capsule to define building blocks for the maps. The output of such a return map yields surfaces for  $\dot{Z}_{k+1}$ and $\psi_{k+1}$ in terms of $\dot{Z}_k$ and $\psi_k$. Return maps based on these building blocks give the foundation for dividing the state space into a few regions with potentially attracting or transient behavior, thus yielding valuable, distinguishing features that can be used for global stability results. Generating polynomial approximations of the exact return maps for $\dot{Z}_k$  and $\psi_k$ on each region in state space, we combine these to obtain a piecewise smooth approximate composite map to reconstruct the system's dynamics. This framework is computationally efficient. It reduces the main computation to constructing polynomial return maps for only short-time realizations of the impact pair over the space of initial conditions. The method calculates a single return, which is a sequence of only a few impacts. This requires limited computation, as compared to, e.g., computing basins of attraction or cell mapping \cite{Rounak_basins2022,Kroetz_basins2018,Wang_cellmap2022}, and in contrast to long-time simulations over the entire state space traditionally used in deriving flow-defined Poincar\'{e} maps for global dynamics of limit-cycle or chaotic systems. Yet, our approximate return maps can be viewed as geometrical models of VI pair systems, analogous to geometrical Lorenz maps used to analyze global dynamics and bifurcations in the chaotic Lorenz system \cite{afraimovic1977origin,robinson1989homoclinic,guckenheimer1983nonlinear} and its more analytically tractable piecewise smooth counterpart \cite{belykh2021sliding}.
Certain aspects of our computation-based analysis do not rely on finding polynomial approximations for the return maps; for example, the efficient comparison of the surfaces projected in the phase planes already identifies potential regions for attracting behavior, on which to focus the computer-assisted analysis. Then, we also pursue polynomial approximations, aiming for explicit expressions for the global analysis.

Anchored in relatively simple return maps, our framework is valuable for cobweb analysis in the phase planes of the state variables.  The relevant global analysis is facilitated by introducing 1D auxiliary maps based on the extreme bounds of the 2D maps in the regions with different types of dynamics. Repeated updates of these auxiliary maps within regions with attracting dynamics yield attraction basins for limit-cycle and chaotic dynamics. Thus, our computer-assisted method of reducing non-smooth systems into a composite piecewise smooth map provides a framework to study the global dynamics of non-smooth systems with impacts. Here, we have focused on parameter regions corresponding to energetically favorable states in VI pair-based energy harvesting systems,  so that the results are relevant for recent designs of VI-based energy harvesters \cite{Yurchenko2017} and nonlinear energy transfer \cite{Kumar2024}. While motivated by a specific vibro-impact energy harvester, nevertheless,  our approach uses generic return maps composed of short sequences of impacts that, in turn, decompose the full dynamics. Thus, the paradigm can be generalized for application in other non-smooth systems.   It may also be interesting to see if this approach,  motivated by a particular class of applied models, is relevant for 2D maps considered in generic mathematical settings \cite{SpecialJDE2023}.

Adapting these findings to realistic external environments remains critical for future exploration. Future work will focus on refining these theoretical frameworks and methodologies to effectively integrate vibro-impact systems into practical applications. This pursuit involves enhancing our understanding of the underlying dynamics and engineering solutions that can withstand and thrive in realistic external environments.

\begin{figure}[H]
    \centering
\includegraphics[width=0.6\textwidth]{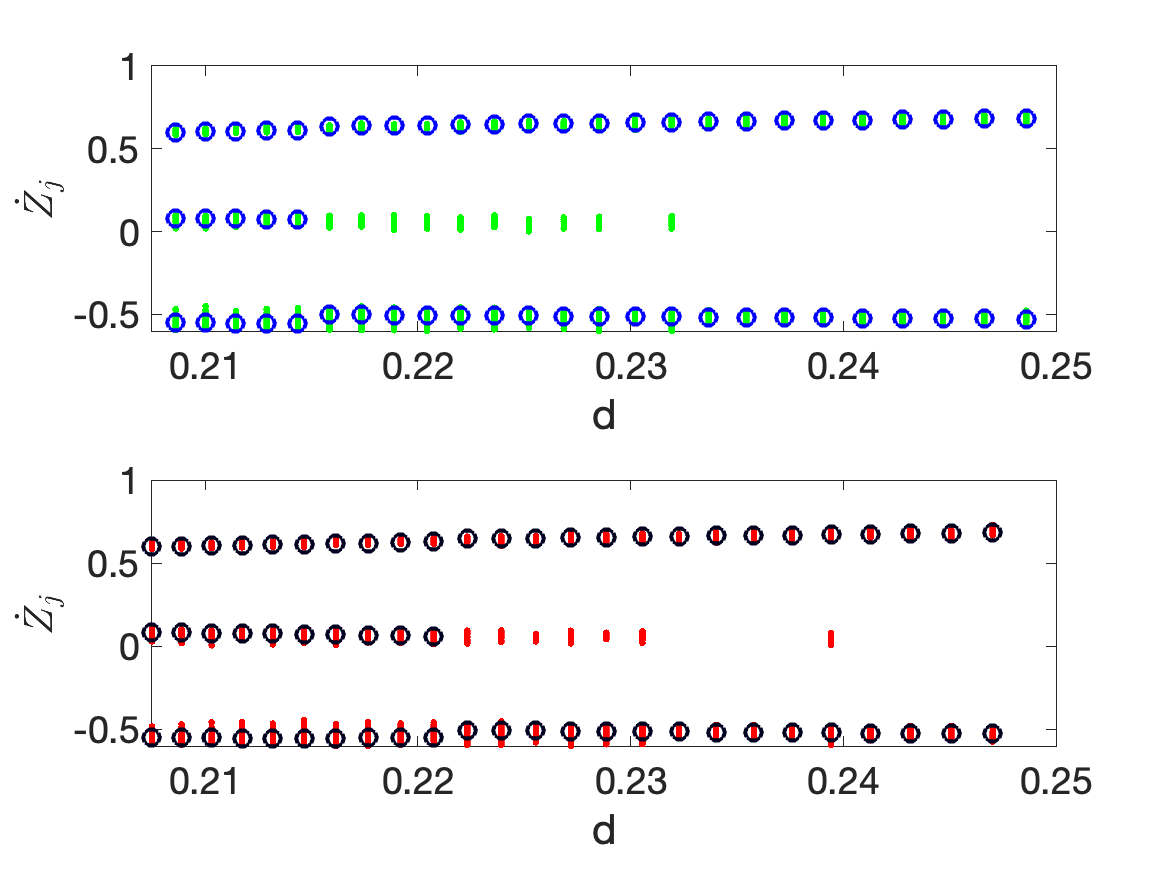}
\caption{Bifurcation diagrams for $\dot{Z}_j$ from \eqref{eq_impact} based on continuation-type methods for decreasing $d$ (top) and increasing $d$ (bottom).  Blue and black open circles correspond to deterministic forcing, and green and red dots correspond to additive noise forcing via an Ornstein-Uhlenbeck process $\zeta$, with limiting behavior $\zeta \sim N(0,0.002)$. Parameters: $r= 0.25$, $\beta = \pi/6$.}
\label{comp_stoch}
\end{figure}

One example of a realistic external setting is the consideration of the VI energy harvester,  illustrated in Fig. \ref{model:abc}(a), under stochastic external forcing.   Figure~\ref{comp_stoch} gives the bifurcation structure with two different types of periodic behavior for the system \eqref{cyl_dim}-\eqref{simpact_dim}, shown via the impact velocity $\dot{Z}_j$ vs. the non-dimensional capsule length parameter $d$. Both panels show deterministic (open circles) vs. stochastic (dots) results for $\dot{Z}_j$.   The top and bottom panels show bifurcation diagrams obtained via a continuation-type method for decreasing and increasing $d$, respectively.
Comparing these indicates bi-stability of two different periodic behaviors. For larger $d$, we observe  1:1 periodic behavior with alternating impacts on $\partial T$  with $\dot{Z}_j<0$ and $\partial B$ with  $\dot{Z}_j>0$ per forcing period. For smaller $d$, we observe 2:1 behavior  with two impacts on $\partial B$  followed by a single impact on $\partial T$ per forcing period. The bi-stability is apparent from the co-existence of branches for the 1:1 and 2:1 solutions in a range of $d$, approximately $ 0.221<d<.216$. At the same time, the stochastic results shown by the green and red points indicate the regular appearance of 2:1 behavior, even for larger values of $d$ beyond the region of bi-stability.
A preliminary analysis, based on the  algorithm from Section \ref{S4algorithm} with an augmented set of return maps analogous to \eqref{exactmap}, includes both ${\cal P}_{BTB}$ to capture  1:1 behavior and a new map ${\cal P}_{BBTB}$ to capture 2:1  behavior.  These maps capture the attraction to either  1:1 and 2:1 behaviors or both. Furthermore, this novel return map framework also provides critical information about the stochastic sensitivity of the 1:1 behavior. This information can be generated quickly since the surfaces for the maps are generated from short-time simulations. Furthermore, we can again compare the shape of these surfaces and their projections in the phase planes to focus on smaller regions with potentially attracting dynamics after eliminating larger transient regions. These results, together with the geometry of the surfaces of these maps, analogous to those shown in Fig. \ref{3d}, efficiently suggest how the noise can bias the dynamics towards  2:1 behavior when combined with the phase plane analysis. We leave the details of that analysis to future work, noting that the algorithm's combined flexibility and efficiency allow for a straightforward augmentation that includes new return maps representing the 2:1 behavior. Then, within the dynamical characterization of the state space provided by our algorithm, we can study non-smooth dynamics in a stochastic setting.

This paper has focused on the development of a novel return map formulation as the basis for a computer-assisted global analysis, obtaining explicit expressions wherever possible. There are a number of other features that we expect are valuable for future generalizations that we have not pursued here.  For example, we expect that more steps of the algorithm could be automated, such as integrating defined criteria to aid in partitioning and comparing approximations for different orders of polynomials for the composite map. 
Furthermore, while we have given the algorithm in terms of  2D maps for simplicity of exposition, we expect that the ideas of this approach can be adapted to higher dimensions.
In addition, if we relax the demand for a nearly explicit global analysis, we anticipate that accurate auxiliary maps that are purely computation-based could be used to approximate the attracting domain(s).   \\

\bibliographystyle{siam}
\bibliography{SIADS_computer_assisted_global_analysis}

\appendix
\section{Return Maps and Composite Map Construction}

\subsection{Division of state space for the return maps} \label{Appendix-statespace} 
\begin{figure}[htbp]
\centering
\includegraphics[width=0.8\textwidth]{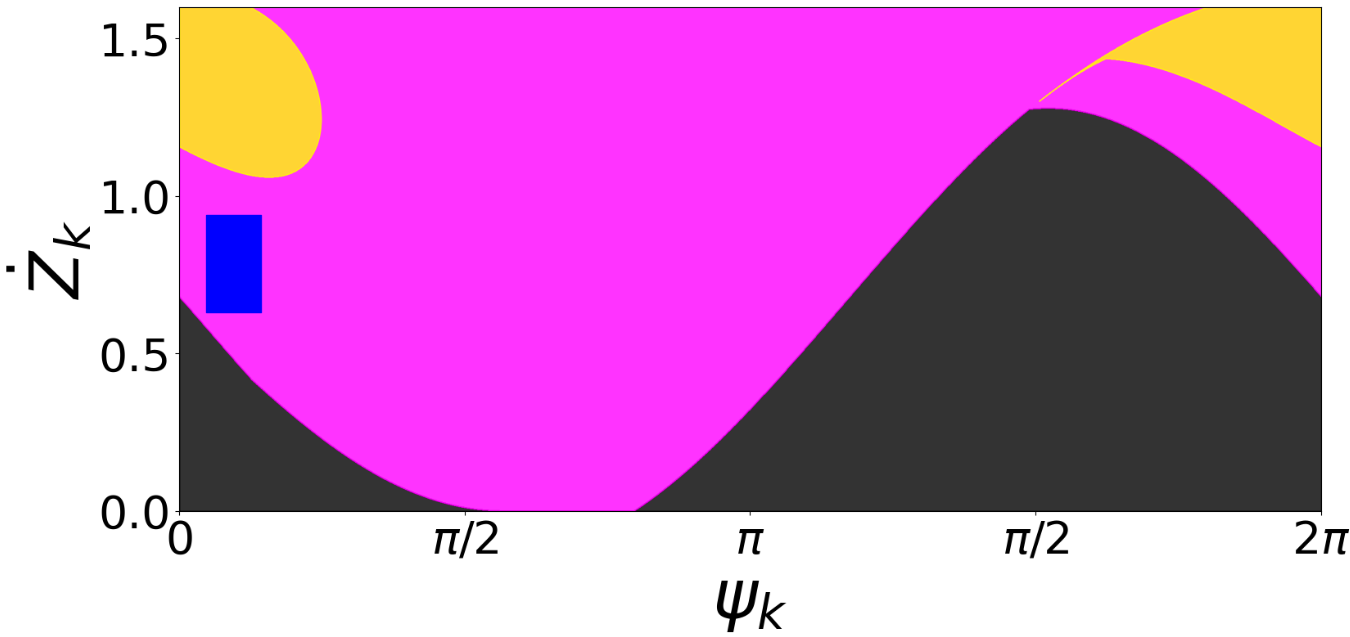}
\caption{Division of the $(\dot{Z}_k, \psi_k)$ state space, corresponding to exact return maps with BTB motion (blue and magenta regions), BB motion (black regions), and BTTB motion (yellow regions). Parameter: $d=0.26$.}
\label{ICPP02pi}
\end{figure}

We show the regions in the state space $(\dot{Z}_k,\psi_k)$ whose images correspond to BB, BTB, and BTTB motion, with $P_{BB}$ and $P_{BTB}$ as defined in \eqref{exactmap} in Section \ref{S3return-map}, and $P_{BTTB}$.  Figure \ref{ICPP02pi}  shows the full range of $\psi_k$, from 0 to $2\pi$, and a larger range of $\dot{Z}_k$ as compared to Fig. \ref{regions}. The region with $\phi_k>\pi$ is comprised of mostly BB motion and, as discussed in  Remark \ref{range} and shown in Fig. \ref{2Dpi-2pi}, is strongly transient. Likewise, the yellow regions, corresponding to BTTB motion, are strongly transient for $\beta>0$, which drives the motion away from multiple impacts on the top membrane $\partial T$. Therefore, we restrict our attention to the state space with range $\psi_k\in [0,\pi]$ and $\dot{Z}_k\leq \; 1.0$ (below the yellow regions) when constructing the composite map $\mathcal{M}$, with a focus on understanding the attracting domain and those regions in state space in close proximity to it. 

\subsection{Phase plane projection of the exact maps}\label{AppendixPrjt}

\begin{figure}[htbp]
\centering
\includegraphics[width=0.8\textwidth]{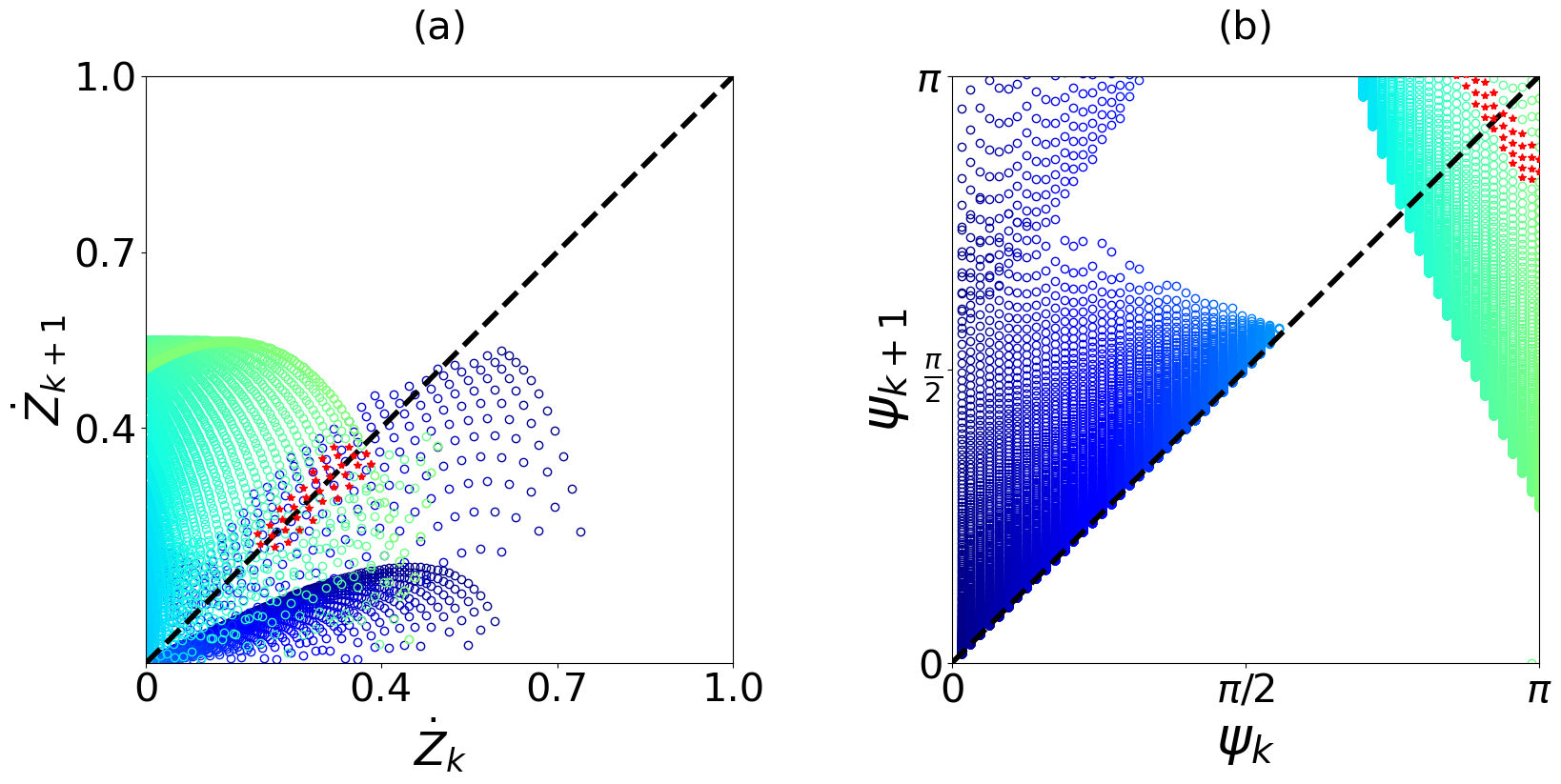}
\includegraphics[width=0.8\textwidth]{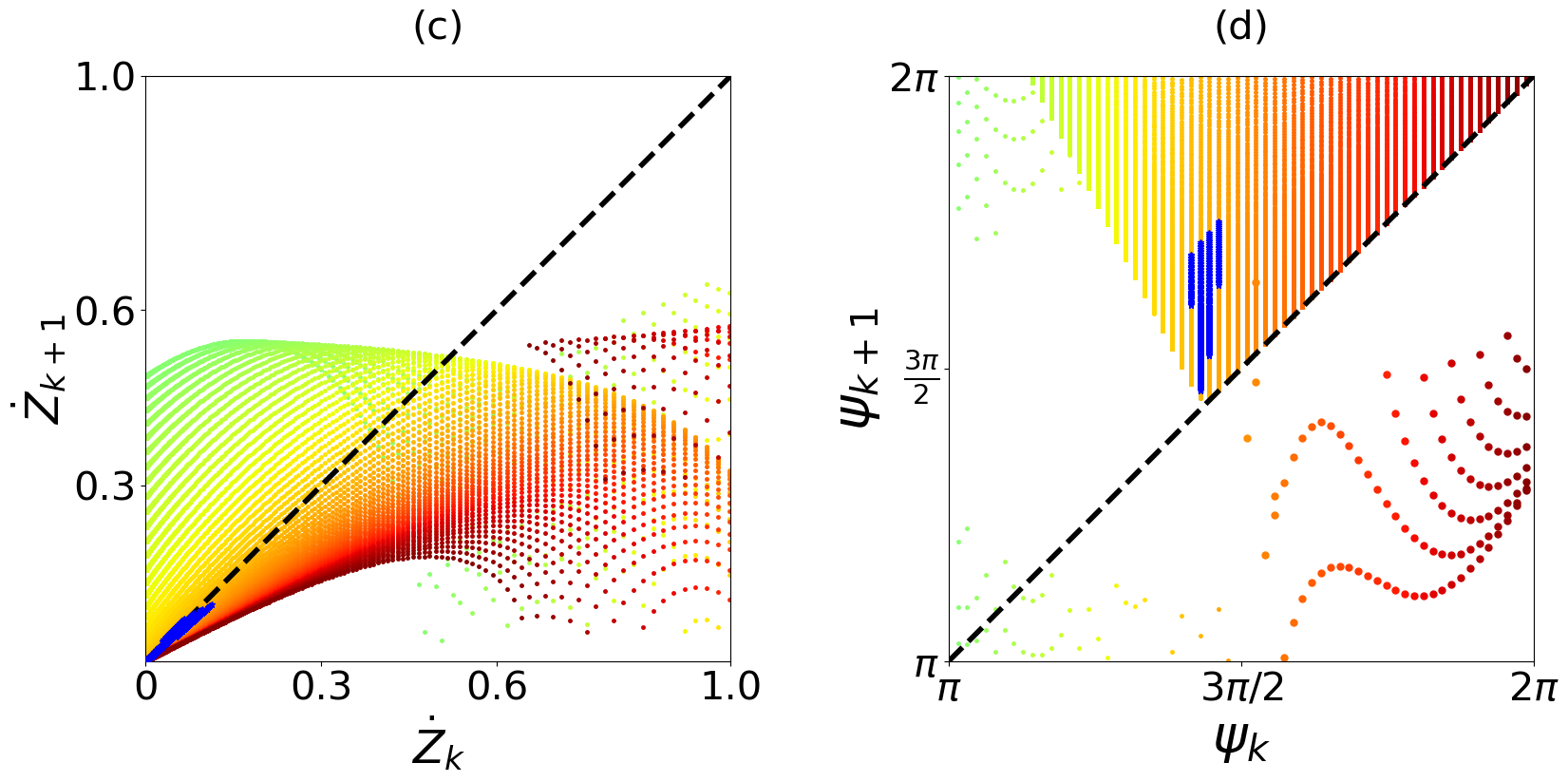}
\caption{ Similar to Fig. \ref{2Dpi-2pi}, we show the 2D projection of  the $P_{BB}$ maps (black surfaces in Fig.~\ref{3d}) on the phase plane $\dot{Z}_{k}-\dot{Z}_{k+1}$ and $\psi_{k}-\psi_{k+1}$ with $d=0.35$. Different colors correspond to the maps for different values of $\psi_k$. (a)-(b) shows results for initial condition $\psi_k\in [0, \pi]$, while  (c)-(d) shows $\psi_k\in [\pi, 2\pi]$. Stars showcase where both maps take values near the diagonals in both phase planes; red stars in (a)-(b) correspond to values on steep slopes of the surfaces, while blue star $\psi_k<3\pi/2$ in (c)-(d), also on surfaces with steep slopes, do not cross the diagonals. These properties indicate transient dynamics for these regions marked with stars.} 
\label{2D0-2pi}
\end{figure}

Figure \ref{2D0-2pi} shows the projections of   the exact maps, defined by \eqref{exactmap} in Section \ref{S3return-map}, on  the  $\dot{Z}_k - \dot{Z}_{k+1}$ and $\psi_{k} -\psi_{k+1}$ phase planes, as referenced in Remark \ref{range}. This 2-D projection of Fig. \ref{3d} gives separate views of the dynamics for $\dot{Z}_k$ and $\psi_k$ in their respective phase planes.
The points delineate curves for  $ \dot{Z}_{k+1}$ and $\psi_{k+1}$ in the image of the return map, some of which cross both diagonals in the $\dot{Z}_{k} - \dot{Z}_{k+1}$ and $\psi_{k }- \psi_{k+1}$ planes. The slopes of the curves that intercept the diagonals suggest that there is a smaller subregion of the state space $(\dot{Z}_k, \psi_k)$ that is attracting. 

\subsection{Comments on Region ${\cal R}_1$}\label{appendix-R1}
In the next six sections of the appendix, we further comment on the details of the algorithm implementation for the specific VI pair model, as discussed in Section \ref{S4.2}.

\begin{figure}[htbp]
\centering
\begin{subfigure}[b]{0.8\textwidth}
   \centering
\includegraphics[width=\textwidth]{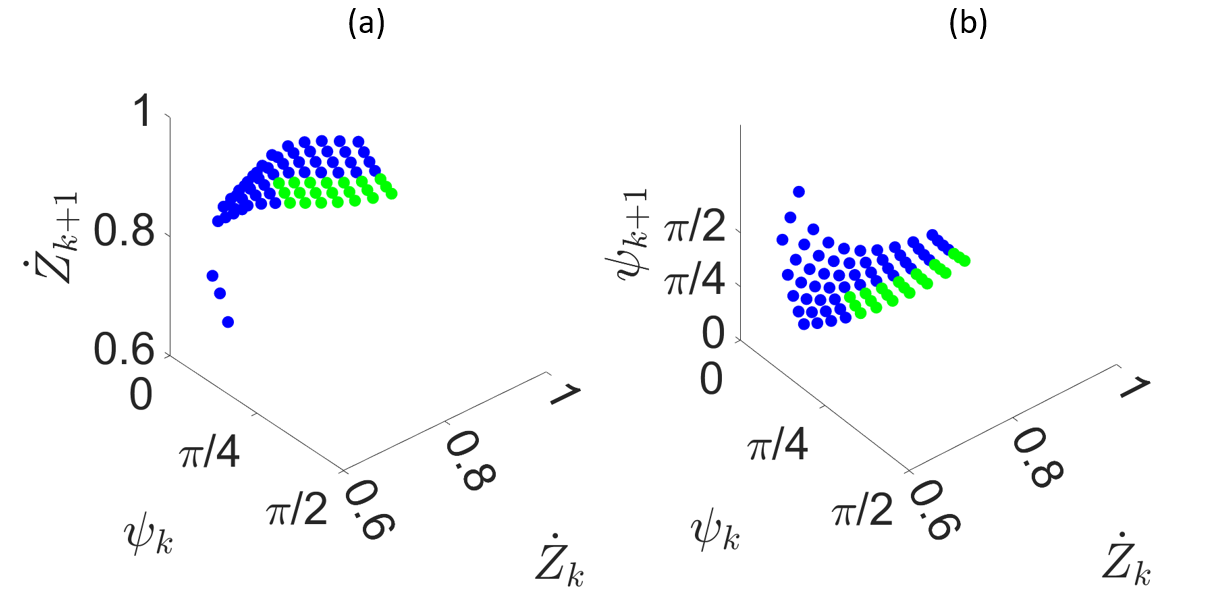}
\end{subfigure}
\vfill
\begin{subfigure}[b]{0.8\textwidth}
   \centering
\includegraphics[width=\textwidth]{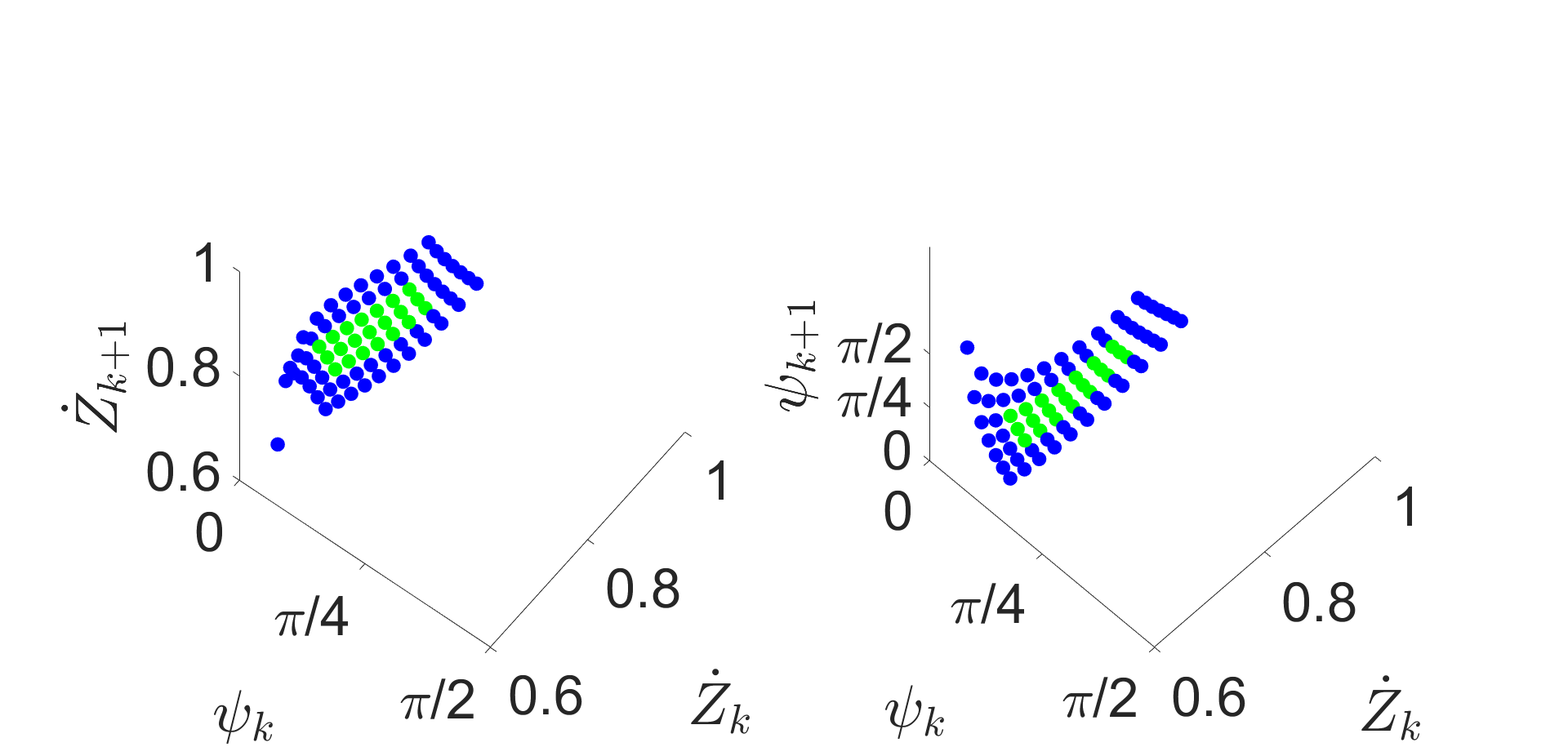}
\end{subfigure}
\vfill
\begin{subfigure}[b]{0.8\textwidth}
   \centering
\includegraphics[width=\textwidth]{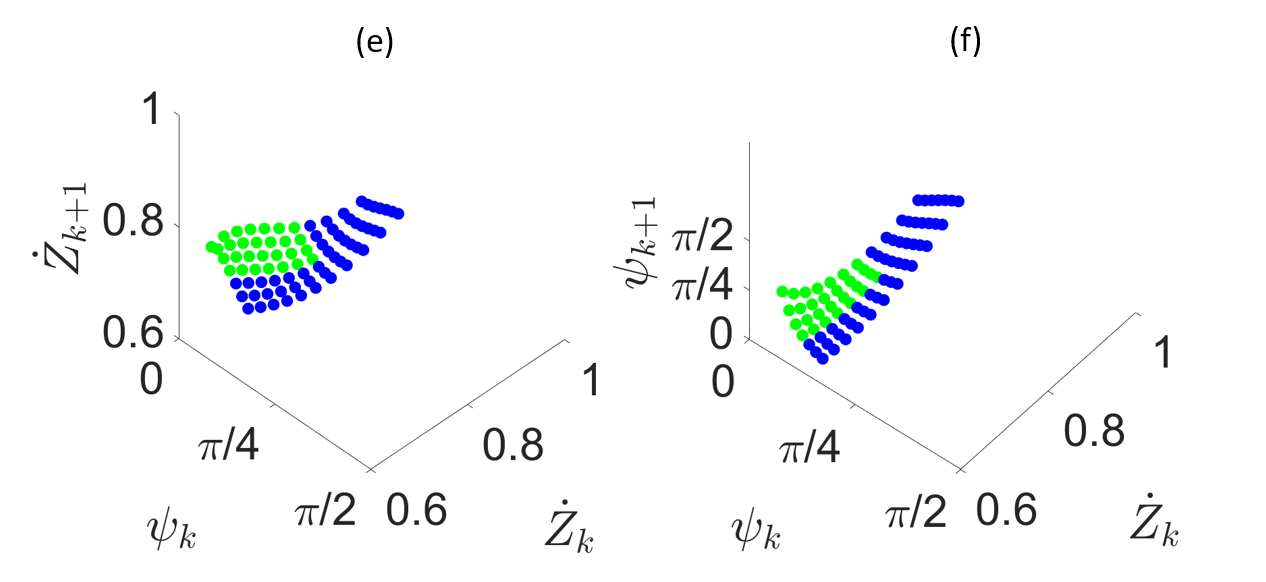}
\end{subfigure}
\caption{Illustration of the location change of the subregions filtered by \eqref{r1uniond}, as shown in green. The blue region surrounding it is the union of all such regions $\cup_{d\in[0.26,0.35]}\mathcal{R}_{1.2}$, as described in \eqref{r1uniond}. Parameters: (a)-(b): $d=0.35$; (c)-(d): $d=0.30$; (e)-(f): $d=0.26$.}
\label{R1compared}
\end{figure}

 In order to capture the full dynamics for all $d$ near the diagonals of both phase planes $\dot{Z}_{k} - \dot{Z}_{k+1}$  and  $\psi_k - \psi_{k+1}$,   we define region $\mathcal{R}_1$  as the union of the subregions  obtained using \eqref{r1uniond}.  Figure \ref{R1compared} illustrates the location of the subregion (green) based on the filter in \eqref{r1uniond} corresponding to one $d$ value. These are shown relative to the union of the  subregions over   all $d$ in the range of interest (blue). Through this definition, we can use the same map for $\mathcal{R}_1$ for all $d$ considered rather than finding different approximate maps for each $d$. 

 We have explored a range of $\delta$ values, $\delta=1.2, 1.3, 1.4$, which is the filter parameter in \eqref{r1uniond}.  In summary, a smaller $\delta$ yields a smaller ${\cal R}_1$ which allows a more accurate approximation of $f_1$ and $g_1$ to the surface of the exact map. On the other hand, a larger ${\cal R}_1$  can capture more dynamics near this region which is desirable. In that case, one can compensate for the increased error associated with larger $\delta$ by increasing the polynomial orders in the approximation. Here, we chose $\delta=1.2$ for the benefit of a simpler expression to construct the approximate map. 

\begin{figure}[htbp]
     \centering
     \begin{subfigure}[b]{\textwidth}
         \centering
    \includegraphics[width=\textwidth]{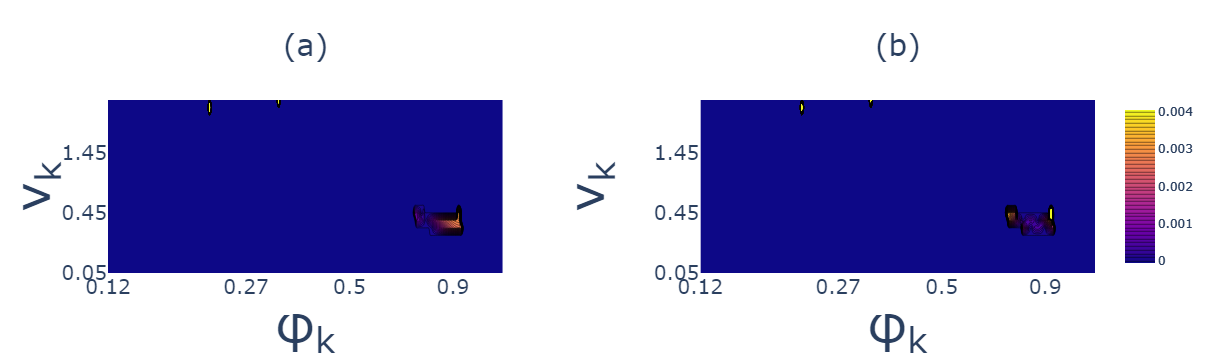}
     \end{subfigure}
     \vfill
    \begin{subfigure}[b]{\textwidth}
         \centering
    \includegraphics[width=\textwidth]{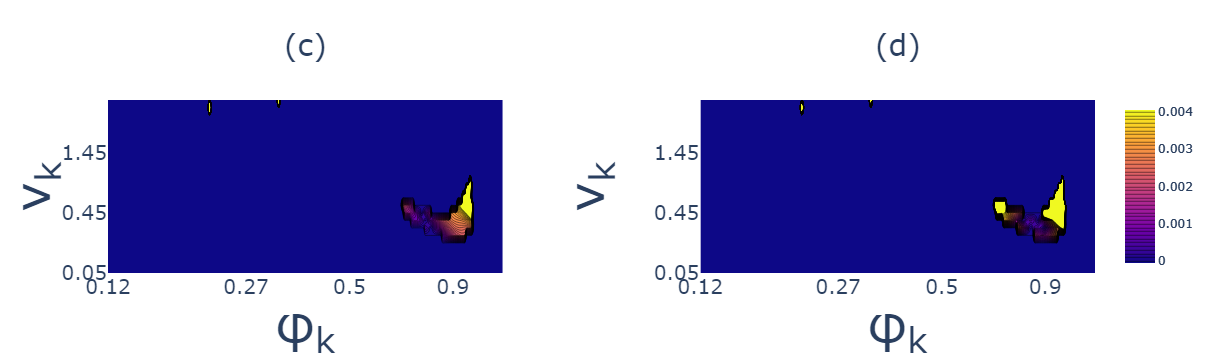}
     \end{subfigure}
     \vfill
    \begin{subfigure}[b]{\textwidth}
         \centering
    \includegraphics[width=\textwidth]{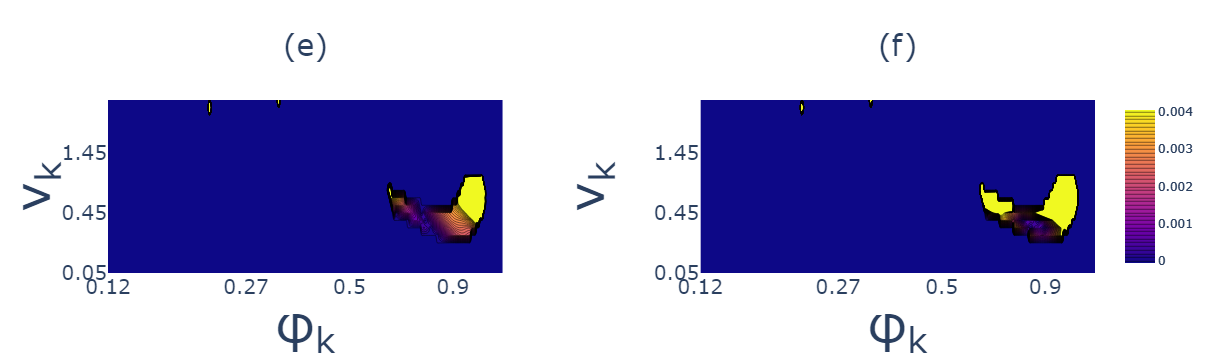}
     \end{subfigure}
\caption{ Heat maps corresponding to the approximation error in Region $\mathcal{R}_1$ with different $\delta$ in \eqref{r1uniond}. The approximation errors $\epsilon_v = |\dot{Z}_{k+1}-v_{k+1}|$ are shown in (a),(c),(e) and $\epsilon_{\phi}=|\psi_{k+1}-\phi_{k+1}|$ are shown in (b),(d),(f) for $(\dot{Z}_{k+1},\phi_{k+1})$ in the exact map and $(v_{k+1},\phi_{k+1})$ in the coupled 2-D approximate map \eqref{cubic2dv}-\eqref{cubic2dp} for ${\cal R}_1$. Note lighter colors indicate larger errors $\epsilon$. As $\delta$ increases, the size of $\mathcal{R}_1$ increases, which includes more variation that yields the larger approximation error. Parameters: $d=0.35$ in all panels; (a)-(b): $\delta=1.2$; (c)-(d): $\delta =1.3$; (e)-(f): $\delta =1.4$.}
\label{comparert}
\end{figure}

In considering the choice for the order of polynomials, we note that higher-order polynomials give more accurate approximations, but this will increase the complexity of the 2D map. Hence, we choose the lowest order polynomial such that the approximation can also reproduce similar dynamics to the exact map. In this case, the polynomial map is quadratic in $\phi_k$  and cubic in $v_k$. Specifically, the polynomials given in the map $(f_1(v_k,\phi_k) , g_1(v_k,\phi_k))$ \eqref{cubic2dv}-\eqref{cubic2dp} in  $ {\cal R}_{1.2}$ approximate the surface using the Matlab function \texttt{fit([x,y],z,fitType)} with argument \texttt{fitType} set to \texttt{"poly23"}. A detailed comparison between the order of the polynomials used in the approximation and the associated error is given in Table \ref{TableDeltaErr} and  Fig. \ref{comparert}.

Table \ref{TableDeltaErr} compares different types of approximation error statistics, $R^2$, and the Summation Squared Error (SSE), using different $\delta$ and different orders of polynomials. Figure \ref{comparert} indicates that a smaller $\delta$ gives a better approximation for a given polynomial order, as a larger $\delta$ includes more variability of the surfaces for ($\dot{Z}_{+1},\psi_{k+1})$. Table \ref{TableDeltaErr} shows that the combination of $\delta=1.2$ and the polynomial order poly23 gives the best result.

\vspace{5 mm}

\begin{table}[htbp]
    \centering
\begin{tabular}{ |p{0.8cm}|p{2cm}|p{2cm}|p{2cm}|p{2cm}|p{2cm}|  }
 \hline
 \multirow{2}{8cm}{$\delta$}& \multirow{2}{8cm}{Poly degree} &\multicolumn{2}{|c|}{$v_{k+1}$} &\multicolumn{2}{|c|}{$\phi_{k+1}$}  \\\cline{3-6}
 & &$R^2$&SSE&$R^2$&SSE\\
 \hline\hline
1.2   & poly23    & 0.9992 & $2.2705\times 10^{-5}$ & 0.9998 & $2.2181\times 10^{-5}$\\[0.1cm]
1.3   & poly23    & 0.99827 & 0.0025092 & 0.99984 & 0.0032939\\[0.1cm]
1.3   & poly33    & 0.99827 & 0.0025055 & 0.99994 & 0.0011577\\[0.1cm]
1.4   & poly23    & 0.99735 & 0.0055033 & 0.99981 & 0.0055713\\[0.1cm]
1.4   & poly33    & 0.99735 & 0.0054874 & 0.9999 & 0.0031359\\[0.1cm]
 \hline 
\end{tabular}
\caption{Comparison of the approximation error $R^2$ and SSE in ${\cal R}_1$ for different $\delta$ and different polynomial orders. Here, $R^2 = 1-\frac{SSE}{SST}$, where the Summation Squared Error and the Summation Squared Total are given by  $SSE = \sum^n_i (y_i-\hat{y}_i)^2$ and $SST = \sum^n_i (y_i-\overline{y})^2$, respectively. Here, $y_i$ is the exact value corresponding to  $\dot{Z}_{k+1}$ or $\psi_{k+1}$, and $\hat{y}_i$ is the estimation $v_{k+1}$ or $\phi_{k+1}$, and $\overline{y}$ is the average of all exact values $\overline{\dot{Z}_{k+1}}$ or $ \overline {\psi_{k+1}}$.}
\label{TableDeltaErr}
\end{table}

\subsection{ Comments on  Region ${\cal R}_2$} \label{appR2}
The surfaces generated over ${\cal R}_2$ correspond to the BTB behavior. As described in Remark \ref{separable}, we use separable maps to represent the dynamics of Region ${\cal R}_2$. Recall that the separable map takes the form of a single variable polynomial, e.g. $v_{k+1}=f_2(v_k)$ and $\phi_{k+1}=g_2(\phi_k)$ \eqref{r2map} in this case. Given the strongly transient nature of the dynamics in ${\cal R}_2$, also indicated by the steep surfaces shown in Fig. \ref{3d}, this 1-D approximation with separable maps is sufficient to represent the dynamics of ${\cal R}_2$.

\subsection{Comments on Region  ${\cal R}_4$} \label{AppR4}
Similar to Region ${\cal R}_2$, the surfaces over ${\cal R}_4$ also correspond to the BTB behavior. However, the surfaces in this region must be approximated separately because of its steep descending surfaces over smaller values of $\dot{Z}_k$, making it difficult to obtain a good approximation over the combined regions of ${\cal R}_2$ and ${\cal R}_4$. The approximate location of ${\cal R}_4$ is given by $\{ (\dot{Z},\psi_k): \dot{Z}_k<0.55,\;  1.1<\psi_k<2.5, \text{ and } \dot{Z}_k>0.63-0.53\psi_k\}$.

Similar to  $\mathcal{R}_2$, we use separable maps for the approximation in  ${\cal R}_4$, choosing two 1-D maps that represent the dynamics given by the surfaces for $\dot{Z}_{k+1}$ and $\psi_{k+1}$
\begin{align}\label{R4eqn}
   v_{k+1}(v_k) &= f_4(v_k)= b_{40} v_k^8 + b_{41} v_k^7 + b_{42} v_k^6 + b_{43} v_k^5 + b_{44} v_k^4 + b_{45} v_k^3 + b_{46} v_k^2 + b_{47} v_k + b_{48},\nonumber\\
   \phi_{k+1}(\phi_k) &= g_4(v_k)= a_{40} \phi_k^4 + a_{41} \phi_k^3 + a_{42} \phi_k^2 + a_{43} \phi_k + a_{44}.
\end{align}
The steep drop of the surface for smaller values of $\dot{Z}_{k+1}$, as shown in Fig.~\ref{r2approx}(f), indicates that the dynamics in $\mathcal{R}_4$ is also strongly transient. That is, at  the fixed point of $v_{k+1} = f_4(v_k)$ the  slope is $|f_4'(v_k)|> 1$, as shown in Fig. \ref{r2approx}(e).

\subsection{Comments on Region  ${\cal R}_3$}
\label{appendixR3}

The approximation for $\mathcal{R}_3$ covers  the surfaces in Fig.~\ref{3d} over the region $\{(\dot{Z}_k,\psi_k): 0<\dot{Z}_k<0.63-0.53\psi_k\}$ within the state space considered. The approximations for the  lower triangular surfaces in this region are given by 
\begin{align} \label{Apr3eqn}
v_{k+1}(v_k, \phi_k) &= f_3(v_k,\phi_k) = b_{300} + b_{301} \phi_k + b_{302} v_k +b_{303} \phi_k^2 + b_{304} \phi_k v_k + b_{305} v_k^2 + b_{306} \phi_k^3 +b_{307} \phi_k^2 v_k \nonumber\\ &+b_{308} \phi_k v_k^2 +b_{309} v_k^3
+ b_{310} \phi_k^3 v_k + b_{311} \phi_k^2 v_k^2 + b_{312} \phi_k v_k^3 
+ b_{313} v_k^4 +b_{314} \phi_k^3 v_k^2 \nonumber\\
&+b_{315} \phi_k^2 v_k^3 +b_{316} \phi_k v_k^4 +b_{317} v_k^5, \nonumber \\
\phi_{k+1}( v_k, \phi_k) &= g_3(v_k,\phi_k) =  a_{300} + a_{301} \phi_k + a_{302} v_k + a_{303} \phi_k^2 +a_{304} \phi_k v_k +a_{305} v_k^2 +a_{306} \phi_k^3 + a_{307} \phi_k^2 v_k \nonumber\\
&+ a_{308} \phi_k v_k^2 
+ a_{309} v_k^3 + a_{310} \phi_k^4 +a_{311} \phi_k^3 v_k +a_{312} \phi_k^2 v_k^2 
+a_{313} \phi_k v_k^3 +a_{314} v_k^4 + a_{315} \phi_k^4 v_k \nonumber\\
&+ a_{316} \phi_k^3 v_k^2 + a_{317} \phi_k^2 v_k^3 
+a_{318} \phi_k v_k^4 + a_{319} v_k^5.
\end{align}

As discussed in Section \ref{S4.1}, Iteration 1 steps iv) and vi), there is also a nearly vertical surface in this region, shown in Fig. \ref{3d}. It represents strongly transient dynamics corresponding to rapid transitions from BB to BTB behavior, so we treat this as immediately transient. As a result, we use the lower triangular surface to capture the dynamics of this region, taking the map \eqref{Apr3eqn} over all of ${\cal R}_3$. We find that these surfaces do not shift or change shape with  $d$ varying. Therefore, the coefficients in \eqref{Apr3eqn} are constant instead of being functions of $d$.

\subsection{Comments on Region  ${\cal R}_5$}\label{AppendixR5}

\begin{figure}[H]
\centering
\begin{subfigure}[t]{0.25\textwidth}
 \centering
 \caption{}
\includegraphics[width=\textwidth]{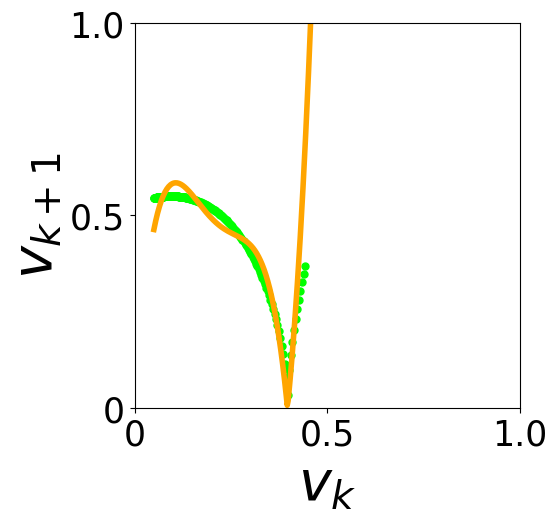}
\end{subfigure}
\hfill
\begin{subfigure}[t]{0.215\textwidth}
 \centering
 \caption{}
 \includegraphics[width=\textwidth]{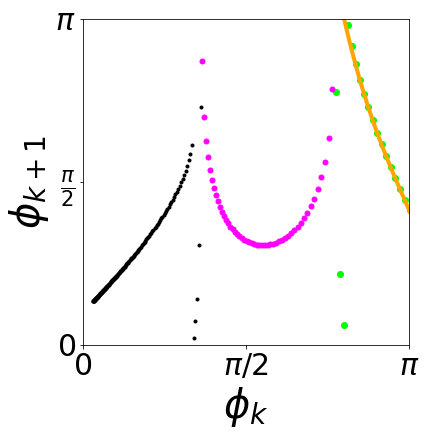}
\end{subfigure}
\hfill
\begin{subfigure}[t]{0.48\textwidth}
 \centering
 \caption{}
 \includegraphics[width=\textwidth]{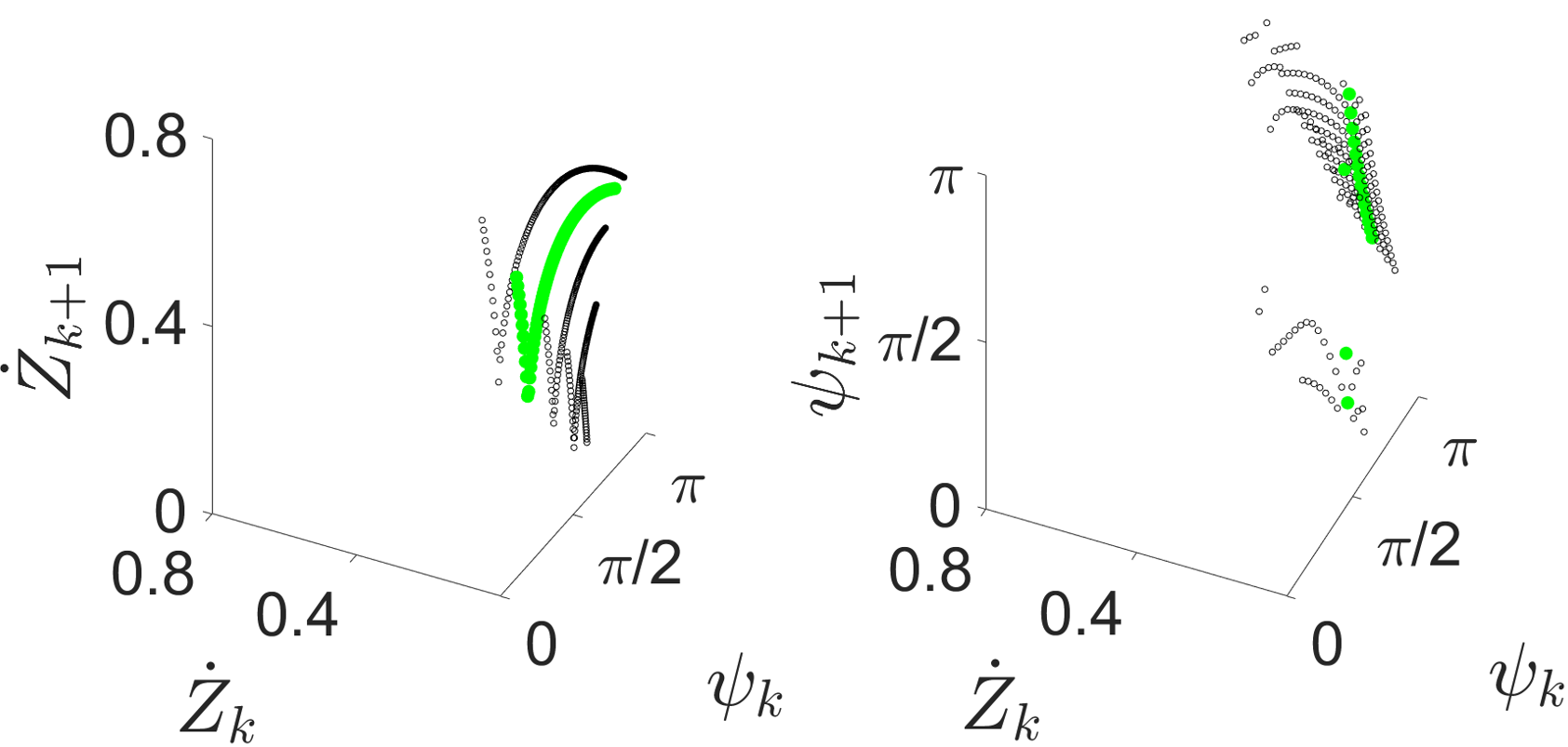}
\end{subfigure}
 \caption{Approximation of $(Z_{k+1},\psi_{k+1})$ in ${\cal R}_5$ for $d=0.35$, which has ranges $\dot{Z}_k<0.55$ and $2.5<\psi_k<\pi$. Panels (a),(b) compare the orange curves for the approximate separable map \eqref{R5eqn} with the green curves in the corresponding phase planes. In panel (c), the green curves are generated with the exact map \eqref{exactmap},  giving a separable representation of the variation of the surface for fixed  $\psi_k=3.05$ (left) and  $\dot{Z}_k = 0.12$ (right).}
 \label{r5approx}
\end{figure}

Region  ${\cal R}_5$ corresponds to smaller $\dot{Z}_k<0.55$, as in ${\cal R}_4$, and for larger $\psi$:  $2.5<\psi_k<\pi$. The dynamics in this region are BB motion instead of BTB motion, with the map $(f_5,g_5)$ based on a separable approximation as in ${\cal R}_2$ and ${\cal R}_4$. The green curves in Fig. \ref{r5approx}(a),(b)  capture the dynamics on the surfaces for $\dot{Z}_{k+1}$ and $\psi_{k+1}$, and are approximated with  orange curves  that give the separable maps
\begin{align}\label{R5eqn}
   v_{k+1}(v_k) &= f_5(v_k)= |b_{50} v_k^4 + b_{51} v_k^3 + b_{52} v_k^2 + b_{53} v_k + b_{54}|,\nonumber\\
   \phi_{k+1}(\phi_k) &= g_5(\phi_k)= a_{50} \phi_k^3 + a_{51} \phi_k^2 + a_{52} \phi_k + a_{53}.
\end{align}
 The coefficients $a_{5i}, b_{5i}, i=0,1,...,4$, are functions of $d$, with $a_{54}=0$ in $\phi_{k+1}$.

 Note there is a small nearly vertical area in the surface for $\psi_{k+1}$, similar to that observed in ${\cal R}_3$ mentioned in Appendix \ref{appendixR3}. As discussed in step vi) of Iteration 1 of the algorithm (Section \ref{S4algorithm}), we treat this as immediately transient, taking the map \eqref{R5eqn} over all of ${\cal R}_5$.

\subsection{The pseudocode used in the programming the composite map}\label{algo_code}

\smallskip

Here, we provide the pseudocode for the approximate composite map for $(v_n,\phi_n)$, as used in Figure \ref{test}, with references to the bounds and maps for each region ${\cal R}_n$. 

\vspace{0.5 mm}
{\bf Algorithm: Composite map for $(v_n,\phi_n)$}
\begin{algorithmic} \label{seudocode}
\If{$\phi_k>\pi$ OR $\phi_k<0$,}\\
 Reset as in  Section \ref{S4.2}, Iteration 1, step vi):  $\phi_{k+1}=1.2$ and $v_{k+1}=v_k$
 \smallskip
\ElsIf{$0.63\le v_k\le 0.94$ AND $0.15\le \phi_k\le 0.45$.} 
    \State Use Region $\mathcal{R}_1$ approximate maps \eqref{cubic2dv}-\eqref{cubic2dp}:
    \smallskip
\ElsIf{$v_k>0.63 - 0.53 \phi_k$ AND $v_k>0.55$ AND $(v_k, \phi_k)\notin \mathcal{R}_1$,}
    \State Use $\mathcal{R}_2$ approximate map \eqref{r2map}:
    \smallskip
\ElsIf{$v_k>0.63 - 0.53 \phi_k$ AND $1.1<\phi_k <2.5$ and $v_k<0.55$,}
    \State Use $\mathcal{R}_4$ approximate map \eqref{R4eqn}:
    \smallskip
\ElsIf{$2.5<\phi_k<\pi$ AND $v_k<0.55$,}
    \State Use $\mathcal{R}_5$ approximate map \eqref{R5eqn}: \smallskip
\ElsIf{$v_k<0.63 - 0.53 \phi_k$,}
    \State Use $\mathcal{R}_3$ approximate map \eqref{Apr3eqn}:
\EndIf
\end{algorithmic}

As discussed in Section \ref{S4.2}, Iteration 1, step iv), the reset value of $\phi_{k+1}=1.2$ can be identified based on the shape of the return maps, indicating that the system moves towards $\phi<\pi/2$. To allow all possible behaviors, a reset value is chosen  in a transient region. There could be other values or distributions of values that would give the same results. The remaining constants in the algorithms are not user defined, but follow from the definition of the maps.  $R_3$ and $R_5$ are obtained in the first definition of the return map (Section \ref{S4.2}, Iteration 1, step iii), $R_2$  is obtained from the part of the BTB region outside of $R_1$ and $R_4$, with the range of $R_4$ based on the slope of the surface in $R_2$. The bounds for $R_1$ depend on the choice of $\delta$ and the order of the polynomial approximation in $(f_1,g_1)$ as described in Appendix \ref{appendix-R1}.

\subsection{Navigation for Fig. \ref{glosmallic}}\label{NaviFig13}

We will use Case FP, shown in Fig. 13 (a)(b), to demonstrate how the orbits are drawn in Fig. 13. 

Step 1: In panel (a), the initial condition is $v_0=0.2$. This has possible images in Region $\mathcal{R}_3$ and $\mathcal{R}_4$ since maps for these three regions overlap. In panel (b), the initial condition is $\phi_0=0.1$. This has possible images in Region $\mathcal{R}_2$ and $\mathcal{R}_3$ since maps for these two regions overlap. Then, the step is taken using the map that is common to both of these, which is Region $\mathcal{R}_3$. Then this map gives us $(v_1, \phi_1) =  ( 0.093, 2.116 )$.

We repeat this process for each $k$. Here we provide the 
 next steps for both $v_k$ and $\phi_k$ until the system  reaches Region $\mathcal{R}_1$: 

Step 2:  We observe from (a) that $v_1=0.093$ has possible images in Region $\mathcal{R}_3$, $\mathcal{R}_4$ and $\mathcal{R}_5$ since maps for these three regions overlap. From (b), $\phi_1$ has possible images in Region $\mathcal{R}_2$ and $\mathcal{R}_4$. Since the region in common $\mathcal{R}_4$, we apply $(f_4, g_4)$ in this step. This gives us the output $(v_2, \phi_2) =  (0.799, 1.150)$.

Step 3: We observe from (a) that $v_2=0.799$ has possible images in Region $\mathcal{R}_1$ and $\mathcal{R}_2$, and $\phi_2=1.150$ has possible images in Region $\mathcal{R}_2$ and $\mathcal{R}_4$. Therefore, the region in common is $\mathcal{R}_2$. In this case, we apply maps $(f_2,g_2)$ to $(v_2, \phi_2)$ and have output $(v_3, \phi_3) =  (0.843, 0.298)$.

Step 4: We observe from (a) that $v_3=0.843$ has possible images in Region $\mathcal{R}_1$ and $\mathcal{R}_2$, and $\phi_3=0.298$ has possible images in Region $\mathcal{R}_1$, $\mathcal{R}_2$, and $\mathcal{R}_3$. In this case, since both $v_3$ and $\phi_3$ have reached $\mathcal{R}_1$, we apply maps $(f_1, g_1)$ to $(v_3, \phi_3)$. The output is $(v_4, \phi_4)=( 0.844, 0.396 )$, and we observe that it is still in the attraction region $\mathcal{R}_1$. 

From this step forward, we observe that the outputs remain in $\mathcal{R}_1$, and hence we repeatedly apply maps $(f_1, g_1)$.

\end{document}


\title{Supplementary material for ``Computer-assisted global analysis for vibro-impact dynamics: a reduced smooth maps approach''}
	\author{Lanjing Bao$^{1}$, Rachel Kuske$^{2}$, Daniil Yurchenko$^{3}$, and Igor Belykh$^{1}$}
	\address{$^1$Department of Mathematics and Statistics, Georgia State University, P.O. Box 4110, Atlanta, Georgia, 30302-410, USA,\\
		$^2$School of Mathematics, Georgia Institute of Technology,
			686 Cherry Street Atlanta, GA 30332-0160 USA,\\
   $^3$Institute for Sound and Vibration Research, University of Southampton, UK}
	\draft \maketitle

\section{Exact map derivations}
Section 2 of the main text summarizes the results in [2] for calculating the exact maps between two consecutive impacts. Ideally, we would like to derive closed-form solutions at the $(j+1)^{\rm th}$ impact $(\dot{Z}_{j+1}, t_{j+1})$ if given the $j^{\rm th}$ impact. The closed-form solutions enable us to study the global dynamics of the VI-EH system. However, we briefly showed in Section 2 of the main text that deriving the closed-form solutions is not feasible. 
In this section, we expand the solving process for the closed-form analytical expressions for the BTB motion. 

BTB motion is composed of maps $P_{TB} \circ P_{BT}$, where map $P_{TB}:(Z_j\in \partial T, \dot{Z}_j, t_j) \; \mapsto \; (Z_{j+1}\in \partial B, \dot{Z}_{j+1}, t_{j+1})$ and map $P_{BT}:(Z_j\in \partial B, \dot{Z}_j, t_j) \; \mapsto \; (Z_{j+1}\in \partial T, \dot{Z}_{j+1}, t_{j+1})$. The equations are 
\begin{align}
&P_{BT}: &\dot{Z}_{j+1} &= -r \dot{Z}_{j} +\Bar{g}\cdot (t_{j+1}-t_{j})+ F_1(t_{j+1})-F_1(t_{j}),\nonumber\\
&   &-\frac{d}{2} &=  \frac{d}{2} -r \dot{Z}_{j} \cdot (t_{j+1}-t_{j})+\frac{\Bar{g}}{2}\cdot (t_{j+1}-t_{j})^2 +F_2(t_{j+1}) - F_2(t_{j}) - F_1(t_{j})\cdot(t_{j+1}-t_{j})\\
&P_{TB}: &\dot{Z}_{j+2} &= -r \dot{Z}_{j+1} +\Bar{g}\cdot (t_{j+2}-t_{j+1})+ F_1(t_{j+2})-F_1(t_{j+1}),\nonumber\\
&   &\frac{d}{2} &=  -\frac{d}{2} -r \dot{Z}_{j+1} \cdot (t_{j+2}-t_{j+1})+\frac{\Bar{g}}{2}\cdot (t_{j+2}-t_{j+1})^2 +F_2(t_{j+2}) - F_2(t_{j+1}) - F_1(t_{j+1})\cdot(t_{j+2}-t_{j+1})
\end{align}
If given $f(t)=F(\pi t+\varphi) = \cos(\pi t+\varphi)$, then $F_1(t) = \frac{1}{\pi} \sin(\pi t+\varphi)$ and $F_2(t) = -\frac{1}{\pi^2}\cos(\pi t+\varphi)$. The maps $P_{BT}, P_{TB}$ can be written as follows: 
\begin{align}
    \dot{Z}_{j+1} &= -r \dot{Z}_{j} +\Bar{g}\cdot (t_{j+1}-t_{j}) +\frac{1}{\pi}\Big(\sin(\pi t_{j+1}+\varphi) - \sin(\pi t_j + \varphi)\Big),\label{p2zj1}\\
    -d &= -r \dot{Z}_{j} \cdot (t_{j+1}-t_{j})+\frac{\Bar{g}}{2}\cdot (t_{j+1}-t_{j})^2 
    -\frac{1}{\pi^2}\Big(\cos(\pi t_{j+1}+\varphi) - \cos(\pi t_j + \varphi)\Big) \label{p2map}\\
     & \qquad \qquad \qquad \qquad \qquad \qquad \qquad \qquad \qquad
    - \frac{1}{\pi}(t_{j+1}-t_{j})\sin(\pi t_j+\varphi),\nonumber\\
    \dot{Z}_{j+2} &= -r \dot{Z}_{j+1} +\Bar{g}\cdot (t_{j+2}-t_{j+1}) +\frac{1}{\pi}\Big(\sin(\pi t_{j+2}+\varphi) - \sin(\pi t_{j+1} + \varphi)\Big)\label{p3j2}\\
    d &= -r \dot{Z}_{j+1} \cdot (t_{j+2}-t_{j+1})+\frac{\Bar{g}}{2}\cdot (t_{j+2}-t_{j+1})^2 
    -\frac{1}{\pi^2}\Big(\cos(\pi t_{j+2}+\varphi) - \cos(\pi t_{j+1} + \varphi)\Big)\label{p3map}\\
     & \qquad \qquad \qquad \qquad \qquad \qquad \qquad \qquad \qquad
    - \frac{1}{\pi}(t_{j+2}-t_{j+1})\sin(\pi t_{j+1}+\varphi)\nonumber
\end{align}

First, unpack \eqref{p2map}:
\begin{align*}
    -d &= -r \dot{Z}_{j} t_{j+1} + r \dot{Z}_{j}t_{j}+\frac{\Bar{g}}{2} t_{j+1}^2 -\Bar{g} t_{j+1} t_{j} + \frac{\Bar{g}}{2} t_{j}^2 
    -\frac{1}{\pi^2}\cos(\pi t_{j+1}+\varphi) + \frac{1}{\pi^2} \cos(\pi t_j + \varphi) \\
    &  \qquad \qquad \qquad \qquad \qquad \qquad \qquad \qquad \qquad
    - \frac{1}{\pi}\sin(\pi t_j+\varphi) t_{j+1} + \frac{1}{\pi}\sin(\pi t_j+\varphi)t_{j}
\end{align*}

Sort all terms containing $t_{j+1}$ such that it's a quadratic equation on the LHS and cosine on the RHS:
\begin{align}\label{p2tj1}
    \frac{\Bar{g}}{2} t_{j+1}^2 - \Big( r\dot{Z}_j + \Bar{g} t_j + \frac{1}{\pi}\sin(\pi t_j+\varphi)\Big) t_{j+1} +\Big( d+ r \dot{Z}_{j}t_{j} + \frac{\Bar{g}}{2} t_{j}^2 + \frac{1}{\pi^2} \cos(\pi t_j + \varphi) + \frac{t_j}{\pi}\sin(\pi t_j+\varphi) \Big)&\nonumber\\
       = \frac{1}{\pi^2} \cos(\pi t_{j+1}+\varphi)\qquad \qquad \qquad \qquad \qquad \qquad&
\end{align}
Equation \eqref{p2tj1} has a solution if the quadratic function on the LHS and the cosine function on the RHS have an intersection. If $t_{j+1}$  has a closed-form solution, then the expression can be written as a function of $\dot{Z}_j$ and $t_j$:
\begin{align} \label{h1}
    t_{j+1} = h_1 (\dot{Z}_j, t_j)
\end{align}
Observe equation \eqref{p2zj1} that $\dot{Z}_{j+1}$ is a function of $\dot{Z}_j, t_{j+1}$, and $t_j$. Apply \eqref{h1} in \eqref{p2zj1}, we can rewrite \eqref{p2zj1} as a function of $\dot{Z}_j$ and $t_j$ only:
\begin{align}\label{h3}
    \dot{Z}_{j+1} &= h_2(\dot{Z}_j, t_{j+1}, t_j)\nonumber\\
    & = h_2(\dot{Z}_j, h_1 (\dot{Z}_j, t_j), t_j)\nonumber\\
    & = h_3 (\dot{Z}_j, t_j)
\end{align}
By unpacking and regrouping \eqref{p3map}, the equation can be written as a quadratic equation of $t_{j+2}$ on the LHS and cosine function of $t_{j+2}$ on the RHS as well. This means $t_{j+2}$ can be written as a function of $\dot{Z}_{j+1}$ and $t_{j+1}$. Apply \eqref{h1} and \eqref{h3}, we have
\begin{align}\label{h5}
    t_{j+2} &= h_4 (\dot{Z}_{j+1}, t_{j+1})\nonumber\\
    &= h_4 (h_3 (\dot{Z}_j, t_j), h_1 (\dot{Z}_j, t_j) )\nonumber\\
    &= h_5(\dot{Z}_j, t_j)
\end{align}
by definition of $\psi_j$, $\psi_{j+2} = \mod (\pi t_{j+2} + \varphi, 2\pi) = h_6(\dot{Z}_j, t_j)$.

Applying \eqref{h1}, \eqref{h3}, \eqref{h5} into \eqref{p3j2}, $\dot{Z}_{j+2}$ can be written as a function of $\dot{Z}_j, t_j$ as well:
\begin{align*}
    \dot{Z}_{j+2} &= h_7(\dot{Z}_{j+1}, t_{j+2}, t_{j+1})\nonumber\\
    &= h_7(h_3 (\dot{Z}_j, t_j), h_5(\dot{Z}_j, t_j), h_1 (\dot{Z}_j, t_j))\nonumber\\
    &= h_8 (\dot{Z}_j, t_j)
\end{align*}

Therefore, the exact map for BTB case $P_{TB}\circ P_{BT}$ can be written as  
\begin{align*}
    \dot{Z}_{j+2} &= h_8(\dot{Z}_j, t_j)\\
    \psi_{j+2} &= h_6(\dot{Z}_j, t_j)
\end{align*}
Solving \eqref{p2tj1} which involves finding $t_{j+1}$ that satisfies both the quadratic equation on the LHS and the cosine function on the RHS is the main obstacle in finding the closed form solution for $(\dot{Z}_{j+2}, t_{j+2})$. If we had the explicit expression for $t_{j+1}=h_1(\dot{Z}_j, t_j)$, the explicit expressions for $h_i(\dot{Z}_j, t_j), i=2,3,...,8$ would follow. However, it is not possible to write down an explicit expression for the solution to \eqref{p2tj1}. Therefore, we are not able to find the closed form expressions for $\dot{Z}_{j+2}$ and $\psi_{j+2}$.

\section{Coefficients for the Composite Map}
\label{coefficients}
In Section 4 of the main text, an algorithm for constructing the composite map is developed. The composite map combines the approximate return maps for each subregion $\mathcal{R}_i$ for $i=1,2,3,4,5$ in Fig. 3(b) of the main text. The approximate return maps are given in terms of the variables $(v_k,\phi_k)$ that denote the approximate relative impact velocity on $\partial B$ and the corresponding impact phase, respectively, at the $k^{\rm th}$ return to $\partial B$.
In this section, we give the specific coefficients of the approximate maps.

\subsection{Region $\mathcal{R}_1$}
The polynomial approximate map of $\mathcal{R}_1$:
\begin{equation} \label{region1_eqn}
    \begin{split}
            g_1(v_k, \phi_k) &= a_0+a_1 \phi_k + a_2 v_k + a_3 \phi_k^2 + a_4 \phi_k v_k + a_5 v_k^2  + a_6 \phi_k^2 v_k + a_7 \phi_k v_k^2 + a_8 v_k^3,\\
    f_1(v_k, \phi_k) &= b_0+b_1 \phi_k + b_2 v_k + b_3 \phi_k^2 + b_4 \phi_k v_k + b_5 v_k^2 + b_6 \phi_k^2 v_k + b_7 \phi_k v_k^2 + b_8 v_k^3,
    \end{split}
\end{equation}
where the coefficients are functions of $d$.
\begin{align*}
a_0 &= -1.499 d^2 + 18.39 d + 10.21,
&  b_0 &= -381.7 d^2 +210 d -27.04,  \\
a_1 &= -196.5 d^2 + 146.2 d -51.59 ,        
&  b_1 &= 777.6 d^2 -438.9 d +59.89 ,\\
a_2 &= 81.56 d^2 -47.97 d -28.93     , 
&  b_2 &= 1036 d^2 -567.3 d +76.38, \\
a_3 &= 257.1 d^2 -189.4 d + 45.34,
&  b_3 &= -487.7 d^2 +278.1 d -38.23, \\
a_4 &= 380.3 d^2 -321.8 d +104.5,
&  b_4 &= -1504 d^2 +860.1 d -121,\\
a_5 &= -218.2 d^2 +125.2 d +7.025,
&  b_5 &= -961.4 d^2 +535.5 d -75.31,\\
a_6 &= -361 d^2 +268.6 d -59.56,
&  b_6 &= 599.4 d^2 -345.7 d +48.48,\\
a_7 &= -84.22 d^2 +91.69 d -35.86,
&  b_7 &= 706.9 d^2 -413.1 d +60.34,\\
a_8 &= 167.4 d^2 -111.7 d +14.11,
&  b_8 &= 313.2 d^2 -180 d +26.68.
\end{align*}

\subsection{Region $\mathcal{R}_2$}
The polynomial approximate map of $\mathcal{R}_2$:
\begin{equation} \label{Supr2eqn}
    \begin{split}
    g_2(\phi_k) &= a_{20} \phi_k^5 + a_{21} \phi_k^4 + a_{22} \phi_k^3 + a_{23} \phi_k^2 + a_{24} \phi_k + a_{25},\\
    f_2(v_k) &= b_{20} v_k^5 + b_{21} v_k^4 + b_{22} v_k^3 + b_{23} v_k^2 + b_{24} v_k + b_{25},
    \end{split}
\end{equation}
where the coefficients are functions of $d$.
\begin{align*}
a_{20} &= -314721.3 d^5 +491841.99 d^4 -306600.36 d^3 +95280.8 d^2
      -14757.75 d  +910.99,\\
a_{21} &= 3254508.65 d^5 -5091024.58 d^4 +3176650.07 d^3 -988128.92 d^2
       +153191.85 d  -9465.61,\\
a_{22} &= -12716817.41 d^5 +19914833.37 d^4 -12439815.01 d^3
       +3873678.24 d^2   -601181.49 d +37187.01,\\
a_{23} &= 23067186.96 d^5  -36172597.397 d^4  +22625380.93 d^3
      -7054633.67 d^2  +1096273.23 d    -67901.79,\\
a_{24} &= -18352978.56 d^5  +28821687.94 d^4 -18053223.19 d^3
       +5636927.09 d^2   -877183.51 d  +54409.05,\\
a_{25} &= 4085295.603 d^5 -6431417.18 d^4  +4038349.16 d^3 -1263982.02 d^2
       +197165.56 d   -12259.01,
\end{align*}
\begin{align*}
    b_{20} &= -8423791.87 d^4 + 10162592.6 d^3 - 4551825.9 d^2 + 895903.8 d - 65176,  \\
    b_{21} &= 39053115.4 d^4 - 47089789.2 d^3 + 21089709.9 d^2 - 4152787.8 d + 302441.7, \\
    b_{22} &= -72167960.4 d^4 + 86937329.1 d^3 - 38914922.6 d^2 + 7662387.6 d - 558347.5, \\
    b_{23} &= 66515161.4 d^4 - 80025768.3 d^3 + 35789133.6 d^2 - 7043852.7 d + 513345.19, \\
    b_{24} &= -30587995.6 d^4 + 36746130.9 d^3 - 16414881.8 d^2 + 3228420.6 d - 235247.1,\\
    b_{25} &= 5609812.9 d^4 - 6728373.2 d^3 + 3001779.6 d^2 - 589859.9 d + 42967.7.
\end{align*}

\subsection{Region $\mathcal{R}_3$}
The polynomial approximate map of $\mathcal{R}_3$:
\begin{align} \label{Supr3eqn}
g_3( v_k, \phi_k) &= a_{300} + a_{301} \phi_k + a_{302} v_k + a_{303} \phi_k^2 +a_{304} \phi_k v_k +a_{305} v_k^2 +a_{306} \phi_k^3 + a_{307} \phi_k^2 v_k + a_{308} \phi_k v_k^2 \nonumber\\
&+ a_{309} v_k^3 + a_{310} \phi_k^4 +a_{311} \phi_k^3 v_k +a_{312} \phi_k^2 v_k^2 
+a_{313} \phi_k v_k^3 +a_{314} v_k^4 + a_{315} \phi_k^4 v_k \nonumber\\
&+ a_{316} \phi_k^3 v_k^2 + a_{317} \phi_k^2 v_k^3 
+a_{318} \phi_k v_k^4 + a_{319} v_k^5,\nonumber\\
f_3(v_k, \phi_k) &= b_{300} + b_{301} \phi_k + b_{302} v_k +b_{303} \phi_k^2 + b_{304} \phi_k v_k + b_{305} v_k^2 + b_{306} \phi_k^3 +b_{307} \phi_k^2 v_k +b_{308} \phi_k v_k^2 +b_{309} v_k^3\nonumber\\
&+ b_{310} \phi_k^3 v_k + b_{311} \phi_k^2 v_k^2 + b_{312} \phi_k v_k^3 
+ b_{313} v_k^4 +b_{314} \phi_k^3 v_k^2 +b_{315} \phi_k^2 v_k^3 +b_{316} \phi_k v_k^4 +b_{317} v_k^5,
\end{align}
where the coefficients are constants.
\begin{align*}
a_{300} &= -4.708\cdot 10^{-5}, & a_{310} &= 0.02214 ,& b_{300} &= 3.311 \cdot 10^{-5}, & b_{310} &= 0.09745,\\
a_{301} &= 0.99, & a_{311} &= -45.86, & b_{301} &= 0.0002375, & b_{311} &= 16.8,\\
a_{302} &= 3.456, & a_{312} &= -235.6, & b_{302} &= 0.4358 ,& b_{312} &= 44.31,\\
a_{303} &= 0.0483, & a_{313} &= -323.4, & b_{303} &= -0.0001751, & b_{313} &= 29.58,\\
a_{304} &= -11.35, & a_{314} &= -148, & b_{304} &= 0.268, & b_{314} &= -8.853,\\
a_{305} &= -13.29, & a_{315} &= 18.32, & b_{305} &= 1.895, & b_{315} &= -38.48,\\
a_{306} &= -0.06063, & a_{316} &= 143, & b_{306} &= 8.499\cdot 10^{-6}, & b_{316} &= -51.93,\\
a_{307} &= 39.33, & a_{317} &= 331.9, & b_{307} &= -0.3043, & b_{317} &= -24.49,\\
a_{308} &= 111.7, & a_{318} &= 308.4, & b_{308} &= -10.54, & \\
a_{309} &= 70.26, & a_{319} &= 114.2, & b_{309} &= -12.81. & \\
\end{align*}

\subsection{Region $\mathcal{R}_4$}
The polynomial approximate map of $\mathcal{R}_4$:

\begin{equation}\label{r4eqn}
    \begin{split}
    g_4(\phi_k) &= a_{40} \phi_k^4 + a_{41} \phi_k^3 + a_{42} \phi_k^2 + a_{43} \phi_k + a_{44},\\
    f_4(v_k) &= b_{40} v_k^8 + b_{41} v_k^7 + b_{42} v_k^6 + b_{43} v_k^5 + b_{44} v_k^4 + b_{45} v_k^3 + b_{46} v_k^2 + b_{47} v_k + b_{48},
    \end{split}
\end{equation}
where the coefficients are functions of $d$.
\begin{align*}
a_{40} &= -25564661 d^5 + 38856593 d^4 - 23532532 d^3 + 7099885 d^2 - 1067289 d + 63961,\\
a_{41} &= 187346514 d^5 - 284624988 d^4 + 172304032 d^3 - 51964934 d^2 + 7808829 d - 467815,\\
a_{42} &= -508479594 d^5 + 772240827 d^4 - 467346559 d^3 + 140905675 d^2 - 21168395 d + 1267853,\\
a_{43} &= 605074088 d^5 - 918738962 d^4 + 555892718 d^3 - 167571538 d^2 + 25170155 d - 1507297,\\
a_{44} &= -267117434 d^5 + 405554166 d^4 - 245366553 d^3 + 73959909 d^2 - 11108535 d + 665192,
\end{align*}
\begin{align*}
    b_{40} &= -33678323446 d^4 + 39732483684 d^3 - 17535685854 d^2 + 3431234055 d - 251148526,\\
    b_{41} &= 83698133214 d^4 - 98744923307 d^3 + 43580936553 d^2 - 8527653549 d + 624188381,\\
    b_{42} &= -87552753895 d^4 + 103292995807 d^3 - 45588589509 d^2 + 8920591453 d - 652957616,\\
    b_{43} &= 50107657144 d^4 - 59115843570 d^3 + 26090961385 d^2 - 5105403132 d + 373701618,\\
    b_{44} &= -17068153916 d^4 + 20136255271 d^3 - 8887112566 d^2 + 1738997288 d - 127289906,\\
    b_{45} &= 3522641275 d^4 - 4155665872 d^3 + 1834037891 d^2 - 358869361 d + 26267846,\\
    b_{46} &= -427643189 d^4 + 504439390 d^3 - 222608632 d^2 + 43555467 d - 3187931,\\
    b_{47} &= 27780980 d^4 - 32762498 d^3 + 14455481 d^2 - 2827938 d + 206958,\\
    b_{48} &= -737999 d^4 + 869871 d^3 - 383650 d^2 + 75032 d - 5489.
\end{align*}

\subsection{Region $\mathcal{R}_5$}
The polynomial approximate map of $\mathcal{R}_5$:

\begin{equation}\label{r5eqn}
    \begin{split}
    g_5(\phi_k) &= a_{50} \phi_k^3 + a_{51} \phi_k^2 + a_{52} \phi_k + a_{53},\\
    f_5(v_k) &= |b_{50} v_k^4 + b_{51} v_k^3 + b_{52} v_k^2 + b_{53} v_k + b_{54}|,
    \end{split}
\end{equation}
where the coefficients are functions of $d$.
\begin{align*}
a_{50} &= 2064.98 d^4 - 1231.18 d^3 - 75.3328 d^2 + 138.871 d - 19.476,\\
a_{51} &= -19752.202 d^4 + 12355.348 d^3 + 119.244 d^2 - 1133.696 d + 166.629,\\
a_{52} &= 61428.79 d^4 - 39362.33 d^3 + 662.836 d^2 + 3177.32 d - 485.139, \\
a_{53} &= -62366.62 d^4 + 40245.2 d^3 - 1078.83 d^2 - 3068.1 d + 482.49 ,
\end{align*}
\begin{align*}
    b_{50} &= -3327935009 d^4 + 4251589868 d^3 - 2036587076 d^2 + 433686951 d - 34659098 ,\\
    b_{51} &= 49128168 d^4 - 628668996 d^3 + 301980243 d^2 - 64578564 d + 5193379 ,\\
    b_{52} &= -24532591 d^4 + 31293322 d^3 - 15007193 d^2 + 3211068 d - 259329,\\
    b_{53} &= 438110 d^4 - 552384 d^3 + 262235 d^2 - 55690.3 d + 4496.38,\\
    b_{54} &=  -5.8882 d^4 + 7.2206 d^3 - 3.2965 d^2 + 0.6646 d - 0.0499.
\end{align*}

\section{The coefficients for the $v$ second-iterate map }
In Section 6 of the main text, a method using the auxiliary map approach is developed. In regions where the approximate maps are two-dimensional, we decouple the 2D systems into the 1D auxiliary maps. We construct a new composite map $\mathcal{M}_{\mathcal{A}}^{(N)}$, defined in (6.3), which assists us in identifying the global attracting region of the system. In the meantime, a higher-order iterate map is derived using $\mathcal{M}_{\mathcal{A}}^{(N)}$ to show the global dynamics and to pinpoint the location of the global absorbing domain. 

In the case of $v_k$, the second iterate map, equation (6.10) in the main text, has a closed-form expression:

\begin{align}
    v_{k+2}(v_k; \phi_{\rm min}, \phi_{\rm max}) &= f_n(f_n(v_k, \phi_{\rm max}), \phi_{\rm min}) \nonumber\\
    &= \alpha_0 + \alpha_1 v_k^1 + \alpha_2 v_k^2 + \alpha_3 v_k^3 + \alpha_4 v_k^4 + \alpha_5 v_k^5 + \alpha_6 v_k^6 + \alpha_7 v_k^7 + \alpha_8 v_k^8 + \alpha_9 v_k^9,
    \label{seconditerativev}
\end{align}

where $\alpha_i,i=1,...,9$ can be calculated if given the parameters $d, \phi_{\rm min}, \phi_{\rm max}$ since $b_1,...,b_8$ are functions of $d$. Each coefficient $\alpha_i$ are polynomials with combinations $b_0, b_1, ..., b_9$, $\phi_{\rm min}$, and $\phi_{\rm max}$: 

\begin{align*}
    \alpha_0 &= b_0 + b_0 b_2 + b_0^2 b_5 + b_0^3 b_8 + b_1 \phi_{\rm min} + b_0 b_4 \phi_{\rm min} + b_0^2 b_7 \phi_{\rm min} + 
 b_3 \phi_{\rm min}^2 + b_0 b_6 \phi_{\rm min}^2 + b_1 b_2 \phi_{\rm max} + 2 b_0 b_1 b_5 \phi_{\rm max} \\
 &+ 3 b_0^2 b_1 b_8 \phi_{\rm max} + 
 b_1 b_4 \phi_{\rm min} \phi_{\rm max} + 2 b_0 b_1 b_7 \phi_{\rm min} \phi_{\rm max} + b_1 b_6 \phi_{\rm min}^2 \phi_{\rm max} + b_2 b_3 \phi_{\rm max}^2 + 
 b_1^2 b_5 \phi_{\rm max}^2 + 2 b_0 b_3 b_5 \phi_{\rm max}^2 \\
 &+ 3 b_0 b_1^2 b_8 \phi_{\rm max}^2 + 
 3 b_0^2 b_3 b_8 \phi_{\rm max}^2 + b_3 b_4 \phi_{\rm min} \phi_{\rm max}^2 + b_1^2 b_7 \phi_{\rm min} \phi_{\rm max}^2 + 
 2 b_0 b_3 b_7 \phi_{\rm min} \phi_{\rm max}^2 + b_3 b_6 \phi_{\rm min}^2 \phi_{\rm max}^2 \\
 &+ 2 b_1 b_3 b_5 \phi_{\rm max}^3 + 
 b_1^3 b_8 \phi_{\rm max}^3 + 6 b_0 b_1 b_3 b_8 \phi_{\rm max}^3 + 2 b_1 b_3 b_7 \phi_{\rm min} \phi_{\rm max}^3 + 
 b_3^2 b_5 \phi_{\rm max}^4 + 3 b_1^2 b_3 b_8 \phi_{\rm max}^4 + 3 b_0 b_3^2 b_8 \phi_{\rm max}^4 \\
 &+ 
 b_3^2 b_7 \phi_{\rm min} \phi_{\rm max}^4 + 3 b_1 b_3^2 b_8 \phi_{\rm max}^5 + b_3^3 b_8 \phi_{\rm max}^6 +3 b_3^2 b_6 b_8 \phi_{\rm max}^6   \\
     \alpha_1 &= b_2^2 + 2 b_0 b_2 b_5 + 3 b_0^2 b_2 b_8 + b_2 b_4 \phi_{\rm min} + 2 b_0 b_2 b_7 \phi_{\rm min} + 
 b_2 b_6 \phi_{\rm min}^2 + b_2 b_4 \phi_{\rm max} + 2 b_1 b_2 b_5 \phi_{\rm max} + 2 b_0 b_4 b_5 \phi_{\rm max} \\
 &+ 
 6 b_0 b_1 b_2 b_8 \phi_{\rm max} + 3 b_0^2 b_4 b_8 \phi_{\rm max} + b_4^2 \phi_{\rm min} \phi_{\rm max} + 2 b_1 b_2 b_7 \phi_{\rm min} \phi_{\rm max} + 
 2 b_0 b_4 b_7 \phi_{\rm min} \phi_{\rm max} + b_4 b_6 \phi_{\rm min}^2 \phi_{\rm max} \\
 &+ 2 b_2 b_3 b_5 \phi_{\rm max}^2 + 
 2 b_1 b_4 b_5 \phi_{\rm max}^2 + b_2 b_6 \phi_{\rm max}^2 + 2 b_0 b_5 b_6 \phi_{\rm max}^2 + 3 b_1^2 b_2 b_8 \phi_{\rm max}^2 + 
 6 b_0 b_2 b_3 b_8 \phi_{\rm max}^2 + 6 b_0 b_1 b_4 b_8 \phi_{\rm max}^2 \\
 &+ 3 b_0^2 b_6 b_8 \phi_{\rm max}^2 + 
 b_4 b_6 \phi_{\rm min} \phi_{\rm max}^2 + 2 b_2 b_3 b_7 \phi_{\rm min} \phi_{\rm max}^2 + 2 b_1 b_4 b_7 \phi_{\rm min} \phi_{\rm max}^2 + 
 2 b_0 b_6 b_7 \phi_{\rm min} \phi_{\rm max}^2 + b_6^2 \phi_{\rm min}^2 \phi_{\rm max}^2 \\
 &+ 2 b_3 b_4 b_5 \phi_{\rm max}^3 + 
 2 b_1 b_5 b_6 \phi_{\rm max}^3 + 6 b_1 b_2 b_3 b_8 \phi_{\rm max}^3 + 3 b_1^2 b_4 b_8 \phi_{\rm max}^3 + 
 6 b_0 b_3 b_4 b_8 \phi_{\rm max}^3 + 6 b_0 b_1 b_6 b_8 \phi_{\rm max}^3 \\
 &+ 2 b_3 b_4 b_7 \phi_{\rm min} \phi_{\rm max}^3 + 
 2 b_1 b_6 b_7 \phi_{\rm min} \phi_{\rm max}^3 + 2 b_3 b_5 b_6 \phi_{\rm max}^4 + 3 b_2 b_3^2 b_8 \phi_{\rm max}^4 + 
 6 b_1 b_3 b_4 b_8 \phi_{\rm max}^4 + 3 b_1^2 b_6 b_8 \phi_{\rm max}^4 \\
 &+ 6 b_0 b_3 b_6 b_8 \phi_{\rm max}^4 + 
 2 b_3 b_6 b_7 \phi_{\rm min} \phi_{\rm max}^4 + 3 b_3^2 b_4 b_8 \phi_{\rm max}^5 + 6 b_1 b_3 b_6 b_8 \phi_{\rm max}^5 + 
 3 b_3^2 b_6 b_8 \phi_{\rm max}^6   
 \end{align*}
\begin{align*}
    \alpha_2 &= b_2 b_5 + b_2^2 b_5 + 2 b_0 b_5^2 + 3 b_0 b_2^2 b_8 + 3 b_0^2 b_5 b_8 + b_4 b_5 \phi_{\rm min} +
  b_2^2 b_7 \phi_{\rm min} + 2 b_0 b_5 b_7 \phi_{\rm min} + b_5 b_6 \phi_{\rm min}^2 + 2 b_2 b_4 b_5 \phi_{\rm max} \\
 &+ 
 2 b_1 b_5^2 \phi_{\rm max} + b_2 b_7 \phi_{\rm max} + 2 b_0 b_5 b_7 \phi_{\rm max} + 3 b_1 b_2^2 b_8 \phi_{\rm max} + 
 6 b_0 b_2 b_4 b_8 \phi_{\rm max} + 6 b_0 b_1 b_5 b_8 \phi_{\rm max} + 3 b_0^2 b_7 b_8 \phi_{\rm max} \\
 &+ 
 b_4 b_7 \phi_{\rm min} \phi_{\rm max} + 2 b_2 b_4 b_7 \phi_{\rm min} \phi_{\rm max} + 2 b_1 b_5 b_7 \phi_{\rm min} \phi_{\rm max} + 2 b_0 b_7^2 \phi_{\rm min} \phi_{\rm max} +
  b_6 b_7 \phi_{\rm min}^2 \phi_{\rm max} + b_4^2 b_5 \phi_{\rm max}^2 + 2 b_3 b_5^2 \phi_{\rm max}^2 \\
 &+ 2 b_2 b_5 b_6 \phi_{\rm max}^2 + 
 2 b_1 b_5 b_7 \phi_{\rm max}^2 + 3 b_2^2 b_3 b_8 \phi_{\rm max}^2 + 6 b_1 b_2 b_4 b_8 \phi_{\rm max}^2 + 
 3 b_0 b_4^2 b_8 \phi_{\rm max}^2 + 3 b_1^2 b_5 b_8 \phi_{\rm max}^2 + 6 b_0 b_3 b_5 b_8 \phi_{\rm max}^2 \\
 &+ 
 6 b_0 b_2 b_6 b_8 \phi_{\rm max}^2 + 6 b_0 b_1 b_7 b_8 \phi_{\rm max}^2 + b_4^2 b_7 \phi_{\rm min} \phi_{\rm max}^2 + 
 2 b_3 b_5 b_7 \phi_{\rm min} \phi_{\rm max}^2 + 2 b_2 b_6 b_7 \phi_{\rm min} \phi_{\rm max}^2 + 2 b_1 b_7^2 \phi_{\rm min} \phi_{\rm max}^2 \\
 &+ 
 2 b_4 b_5 b_6 \phi_{\rm max}^3 + 2 b_3 b_5 b_7 \phi_{\rm max}^3 + 6 b_2 b_3 b_4 b_8 \phi_{\rm max}^3 + 
 3 b_1 b_4^2 b_8 \phi_{\rm max}^3 + 6 b_1 b_3 b_5 b_8 \phi_{\rm max}^3 + 6 b_1 b_2 b_6 b_8 \phi_{\rm max}^3 + 
 6 b_0 b_4 b_6 b_8 \phi_{\rm max}^3 \\
 &+ 3 b_1^2 b_7 b_8 \phi_{\rm max}^3 + 6 b_0 b_3 b_7 b_8 \phi_{\rm max}^3 + 
 2 b_4 b_6 b_7 \phi_{\rm min} \phi_{\rm max}^3 + 2 b_3 b_7^2 \phi_{\rm min} \phi_{\rm max}^3 + b_5 b_6^2 \phi_{\rm max}^4 + 
 3 b_3 b_4^2 b_8 \phi_{\rm max}^4 + 3 b_3^2 b_5 b_8 \phi_{\rm max}^4 \\
 &+ 6 b_2 b_3 b_6 b_8 \phi_{\rm max}^4 + 
 6 b_1 b_4 b_6 b_8 \phi_{\rm max}^4 + 3 b_0 b_6^2 b_8 \phi_{\rm max}^4 + 6 b_1 b_3 b_7 b_8 \phi_{\rm max}^4 + 
 b_6^2 b_7 \phi_{\rm min} \phi_{\rm max}^4 + 6 b_3 b_4 b_6 b_8 \phi_{\rm max}^5 \\
 &+ 3 b_1 b_6^2 b_8 \phi_{\rm max}^5 + 
 3 b_3^2 b_7 b_8 \phi_{\rm max}^5 + 3 b_3 b_6^2 b_8 \phi_{\rm max}^6   \\
    \alpha_3 &= 2 b_2 b_5^2 + b_2 b_8 + b_2^3 b_8 + 2 b_0 b_5 b_8 + 6 b_0 b_2 b_5 b_8 + 
 3 b_0^2 b_8^2 + 2 b_2 b_5 b_7 \phi_{\rm min} + b_4 b_8 \phi_{\rm min} + 2 b_0 b_7 b_8 \phi_{\rm min} + 
 b_6 b_8 \phi_{\rm min}^2 \\
 &+ 2 b_4 b_5^2 \phi_{\rm max} + 2 b_2 b_5 b_7 \phi_{\rm max} + 3 b_2^2 b_4 b_8 \phi_{\rm max} + 
 2 b_1 b_5 b_8 \phi_{\rm max} + 6 b_1 b_2 b_5 b_8 \phi_{\rm max} + 6 b_0 b_4 b_5 b_8 \phi_{\rm max} + 
 6 b_0 b_2 b_7 b_8 \phi_{\rm max} \\
 &+ 6 b_0 b_1 b_8^2 \phi_{\rm max} + 2 b_4 b_5 b_7 \phi_{\rm min} \phi_{\rm max} + 
 2 b_2 b_7^2 \phi_{\rm min} \phi_{\rm max} + 2 b_1 b_7 b_8 \phi_{\rm min} \phi_{\rm max} + 2 b_5^2 b_6 \phi_{\rm max}^2 + 
 2 b_4 b_5 b_7 \phi_{\rm max}^2 + 3 b_2 b_4^2 b_8 \phi_{\rm max}^2 \\
 &+ 2 b_3 b_5 b_8 \phi_{\rm max}^2 + 
 6 b_2 b_3 b_5 b_8 \phi_{\rm max}^2 + 6 b_1 b_4 b_5 b_8 \phi_{\rm max}^2 + 3 b_2^2 b_6 b_8 \phi_{\rm max}^2 + 
 6 b_0 b_5 b_6 b_8 \phi_{\rm max}^2 + 6 b_1 b_2 b_7 b_8 \phi_{\rm max}^2 + 6 b_0 b_4 b_7 b_8 \phi_{\rm max}^2 \\
 &+ 
 3 b_1^2 b_8^2 \phi_{\rm max}^2 + 6 b_0 b_3 b_8^2 \phi_{\rm max}^2 + 2 b_5 b_6 b_7 \phi_{\rm min} \phi_{\rm max}^2 + 
 2 b_4 b_7^2 \phi_{\rm min} \phi_{\rm max}^2 + 2 b_3 b_7 b_8 \phi_{\rm min} \phi_{\rm max}^2 + 2 b_5 b_6 b_7 p_1^3 + 
 b_4^3 b_8 \phi_{\rm max}^3 \\
 &+ 6 b_3 b_4 b_5 b_8 \phi_{\rm max}^3 + 6 b_2 b_4 b_6 b_8 \phi_{\rm max}^3 + 
 6 b_1 b_5 b_6 b_8 \phi_{\rm max}^3 + 6 b_2 b_3 b_7 b_8 \phi_{\rm max}^3 + 6 b_1 b_4 b_7 b_8 \phi_{\rm max}^3 + 
 6 b_0 b_6 b_7 b_8 \phi_{\rm max}^3 \\
 &+ 6 b_1 b_3 b_8^2 \phi_{\rm max}^3 + 2 b_6 b_7^2 \phi_{\rm min} \phi_{\rm max}^3 + 
 3 b_4^2 b_6 b_8 \phi_{\rm max}^4 + 6 b_3 b_5 b_6 b_8 \phi_{\rm max}^4 + 3 b_2 b_6^2 b_8 \phi_{\rm max}^4 + 
 6 b_3 b_4 b_7 b_8 \phi_{\rm max}^4 \\
 &+ 6 b_1 b_6 b_7 b_8 \phi_{\rm max}^4 + 3 b_3^2 b_8^2 \phi_{\rm max}^4 + 
 3 b_4 b_6^2 b_8 \phi_{\rm max}^5 + 6 b_3 b_6 b_7 b_8 \phi_{\rm max}^5 + b_6^3 b_8 \phi_{\rm max}^6  \\
    \alpha_4 &=  b_5^3 + 2 b_2 b_5 b_8 + 3 b_2^2 b_5 b_8 + 3 b_0 b_5^2 b_8 + 6 b_0 b_2 b_8^2 + 
 b_5^2 b_7 \phi_{\rm min} + 2 b_2 b_7 b_8 \phi_{\rm min} + 2 b_5^2 b_7 \phi_{\rm max} + 2 b_4 b_5 b_8 \phi_{\rm max} \\
 &+ 
 6 b_2 b_4 b_5 b_8 \phi_{\rm max} + 3 b_1 b_5^2 b_8 \phi_{\rm max} + 3 b_2^2 b_7 b_8 \phi_{\rm max} + 
 6 b_0 b_5 b_7 b_8 \phi_{\rm max} + 6 b_1 b_2 b_8^2 \phi_{\rm max} + 6 b_0 b_4 b_8^2 \phi_{\rm max} + 
 2 b_5 b_7^2 \phi_{\rm min} \phi_{\rm max} \\
 &+ 2 b_4 b_7 b_8 \phi_{\rm min} \phi_{\rm max} + b_5 b_7^2 \phi_{\rm max}^2 + 
 3 b_4^2 b_5 b_8 \phi_{\rm max}^2 + 3 b_3 b_5^2 b_8 \phi_{\rm max}^2 + 2 b_5 b_6 b_8 \phi_{\rm max}^2 + 
 6 b_2 b_5 b_6 b_8 \phi_{\rm max}^2 + 6 b_2 b_4 b_7 b_8 \phi_{\rm max}^2 \\
 &+ 6 b_1 b_5 b_7 b_8 \phi_{\rm max}^2 + 
 3 b_0 b_7^2 b_8 \phi_{\rm max}^2 + 6 b_2 b_3 b_8^2 \phi_{\rm max}^2 + 6 b_1 b_4 b_8^2 \phi_{\rm max}^2 + 
 6 b_0 b_6 b_8^2 \phi_{\rm max}^2 + b_7^3 \phi_{\rm min} \phi_{\rm max}^2 + 2 b_6 b_7 b_8 \phi_{\rm min} \phi_{\rm max}^2 \\
 &+ 
 6 b_4 b_5 b_6 b_8 \phi_{\rm max}^3 + 3 b_4^2 b_7 b_8 \phi_{\rm max}^3 + 6 b_3 b_5 b_7 b_8 \phi_{\rm max}^3 + 
 6 b_2 b_6 b_7 b_8 \phi_{\rm max}^3 + 3 b_1 b_7^2 b_8 \phi_{\rm max}^3 + 6 b_3 b_4 b_8^2 \phi_{\rm max}^3 + 
 6 b_1 b_6 b_8^2 \phi_{\rm max}^3 \\
 &+ 3 b_5 b_6^2 b_8 \phi_{\rm max}^4 + 6 b_4 b_6 b_7 b_8 \phi_{\rm max}^4 + 
 3 b_3 b_7^2 b_8 \phi_{\rm max}^4 + 6 b_3 b_6 b_8^2 \phi_{\rm max}^4 + 3 b_6^2 b_7 b_8 \phi_{\rm max}^5  \\
    \alpha_5 &= 2 b_5^2 b_8 + 3 b_2 b_5^2 b_8 + 3 b_2^2 b_8^2 + 6 b_0 b_5 b_8^2 + 
 2 b_5 b_7 b_8 \phi_{\rm min} + 3 b_4 b_5^2 b_8 \phi_{\rm max} + 2 b_5 b_7 b_8 \phi_{\rm max} + 6 b_2 b_5 b_7 b_8 \phi_{\rm max} + 
 6 b_2 b_4 b_8^2 \phi_{\rm max} \\
 &+ 6 b_1 b_5 b_8^2 \phi_{\rm max} + 6 b_0 b_7 b_8^2 \phi_{\rm max} + 
 2 b_7^2 b_8 \phi_{\rm min} \phi_{\rm max} + 3 b_5^2 b_6 b_8 \phi_{\rm max}^2 + 6 b_4 b_5 b_7 b_8 \phi_{\rm max}^2 + 
 3 b_2 b_7^2 b_8 \phi_{\rm max}^2 + 3 b_4^2 b_8^2 \phi_{\rm max}^2 \\
 &+ 6 b_3 b_5 b_8^2 \phi_{\rm max}^2 + 
 6 b_2 b_6 b_8^2 \phi_{\rm max}^2 + 6 b_1 b_7 b_8^2 \phi_{\rm max}^2 + 6 b_5 b_6 b_7 b_8 \phi_{\rm max}^3 + 
 3 b_4 b_7^2 b_8 \phi_{\rm max}^3 + 6 b_4 b_6 b_8^2 \phi_{\rm max}^3 + 6 b_3 b_7 b_8^2 \phi_{\rm max}^3 \\
 &+ 
 3 b_6 b_7^2 b_8 \phi_{\rm max}^4 + 3 b_6^2 b_8^2 \phi_{\rm max}^4   
    \alpha_6 &= b_5^3 b_8 + b_5 b_8^2 + 6 b_2 b_5 b_8^2 + 3 b_0 b_8^3 + b_7 b_8^2 \phi_{\rm min} + 
 3 b_5^2 b_7 b_8 \phi_{\rm max} + 6 b_4 b_5 b_8^2 \phi_{\rm max} + 6 b_2 b_7 b_8^2 \phi_{\rm max} + 3 b_1 b_8^3 \phi_{\rm max} \\
 &+ 
 3 b_5 b_7^2 b_8 \phi_{\rm max}^2 + 6 b_5 b_6 b_8^2 \phi_{\rm max}^2 + 6 b_4 b_7 b_8^2 \phi_{\rm max}^2 + 
 3 b_3 b_8^3 \phi_{\rm max}^2 + b_7^3 b_8 \phi_{\rm max}^3 + 6 b_6 b_7 b_8^2 \phi_{\rm max}^3   \\
    \alpha_7 &=  3 b_5^2 b_8^2 + 3 b_2 b_8^3 + 6 b_5 b_7 b_8^2 \phi_{\rm max} + 3 b_4 b_8^3 \phi_{\rm max} + 
 3 b_7^2 b_8^2 \phi_{\rm max}^2 + 3 b_6 b_8^3 \phi_{\rm max}^2  \\
    \alpha_8 &= 3 b_5 b_8^3 + 3 b_7 b_8^3 \phi_{\rm max}   \\
    \alpha_9 &= b_8^4  
\end{align*}